\documentclass[11pt,a4paper]{amsart}
\usepackage[colorlinks=true,linkcolor=blue,citecolor=red,urlcolor=blue]{hyperref}
 \usepackage[left=1.2in,right=1.2in,bottom=1in,top=1in]{geometry}

\usepackage{tikz,amsthm,amsmath,amstext,amssymb,amscd,epsfig,euscript, mathrsfs, dsfont,pspicture,multicol, mathtools,graphpap,graphics,graphicx,times,enumerate,subfig,sidecap,wrapfig,color}
 \usepackage{url}
\makeatletter
\def\@tocline#1#2#3#4#5#6#7{\relax
  \ifnum #1>\c@tocdepth 
  \else
    \par \addpenalty\@secpenalty\addvspace{#2}%
    \begingroup \hyphenpenalty\@M
    \@ifempty{#4}{%
      \@tempdima\csname r@tocindent\number#1\endcsname\relax
    }{%
      \@tempdima#4\relax
    }%
    \parindent\z@ \leftskip#3\relax \advance\leftskip\@tempdima\relax
    \rightskip\@pnumwidth plus4em \parfillskip-\@pnumwidth
    #5\leavevmode\hskip-\@tempdima
      \ifcase #1
       \or\or \hskip 1em \or \hskip 2em \else \hskip 3em \fi%
      #6\nobreak\relax
    \dotfill\hbox to\@pnumwidth{\@tocpagenum{#7}}\par
    \nobreak
    \endgroup
  \fi}
\makeatother

\def \rr {{\mathbb R}}
\def \cc {{\mathbb C}}

\DeclareGraphicsExtensions{.png}

\def\et{\tilde{\eta}}
\def\ept{\tilde{\epsilon}}
\def\mt{\tilde{\mu}}
\def\ut{\tilde{\upsilon}}
\def\dom{\mathrm{dom}}
\def\CS{\mathit{C}_{\Sigma}}
\def\aN{\alpha\cdot\mathit{N}}
\def\l{\lambda}
\def\ll{\l^{\prime}}

\def\ep{\epsilon}
\def\u{\upsilon}
\def\g{\gamma}
\def\t{\tau}
\def\k{\kappa}
\def\e{\eta}

\def\s{\sigma}
\def\m{\mu}

\def\G{\Gamma}

\def\O{\Omega}
\def\o{\omega}
\def\S{\Sigma}

\def \cc {{\mathbb{C}}}

\newcommand{\bela}{\begin{equation} \label}
\newcommand{\eeq}{\end{equation}}
\newcommand{\ba}{\begin{array}}
\newcommand{\ea}{\end{array}}

\newtheorem{definition}{Definition}[section]
\newtheorem{theorem}{Theorem}[section]
\newtheorem{proposition}{Proposition}[section]
\newtheorem{lemma}{Lemma}[section]
\newtheorem{corollary}{Corollary}[section]
\newtheorem{remark}{Remark}[section]

\catcode`\@=12

\begin{document}


\title[On Dirac operators with $\delta$-Interactions]{Spectral Analysis of Dirac Operators with delta interactions supported on the boundaries of rough domains}

\author[Badreddine Benhellal]{Badreddine Benhellal}
\address{Departamento de Matem\'aticas, Universidad del Pa\' is Vasco, Barrio Sarriena s/n 48940 Leioa, SPAIN\\ and
 Universit\' e de Bordeaux, IMB, UMR 5251, 33405 Talence Cedex, FRANCE.}
\email{benhellal.badreddine@ehu.eus and badreddine.benhellal@u-bordeaux.fr}


\maketitle

\begin{abstract}
 Given  an open set $\O\subset\rr^3$. We deal with the spectral study of Dirac operators of the form $H_{a,\t}=H+A_{a,\t}\delta_{\partial\O}$, where $H$ is the free Dirac operator in $\rr^3$, $A_{a,\t}$ is a bounded invertible, self-adjoint operator in $\mathit{L}^{2}(\partial\O)^4$, depending on parameters $(a,\t)\in\rr\times\rr^n$, $n\geqslant1$.   We investigate the self-adjointness and the related spectral properties of $H_{a,\t}$, such as the phenomenon of confinement and the Sobolev regularity of the domain  in different situations.
 Our set of techniques, which is based on fundamental solutions and layer potentials, allows us to tackle the above problems  under mild geometric measure theoretic assumptions on $\O$. 
 
\end{abstract}

\section{Intoduction} 
Given $m>0$, the free Dirac operator in $\rr^3$ is defined by $H=- i \alpha \cdot\nabla + m\beta$, where  $\alpha= (\alpha_1,\alpha_2,\alpha_3)$ and $\alpha \cdot\nabla =\sum_{k=1}^{3}\alpha_k\partial_k $. Here each $\alpha_j$'s and $\beta$ are  $4\times4$ Hermitian and unitary matrices defined by
\begin{align}\label{les matrice}
\alpha_k=\begin{pmatrix}
0 & \sigma_k\\
\sigma_k & 0
\end{pmatrix} \text{ for } k=1,2,3,
\quad
\beta=\begin{pmatrix}
\mathit{I}_2 & 0\\
0 & -\mathit{I}_2
\end{pmatrix},\quad 
 \gamma_5 :=-i\alpha_1\alpha_2\alpha_3=\begin{pmatrix}
0& \mathit{I}_2 \\
\mathit{I}_2  & 0
\end{pmatrix}, 
\end{align}  
where $\mathit{I}_n$ is the $n\times n$ unit matrix, and  $\sigma = (\sigma_1,\sigma_2,\s_3)$ are the Pauli matrices given by
\begin{align}\label{Pauli}
\sigma_1=\begin{pmatrix}
0 & 1\\
1 & 0
\end{pmatrix},\quad \sigma_2=
\begin{pmatrix}
0 & -i\\
i & 0
\end{pmatrix} ,
\quad
\sigma_3=\begin{pmatrix}
1 & 0\\
0 & -1
\end{pmatrix}.
\end{align} 
 We consider Dirac operators acting on $\mathit{L}^2(\rr^3)^4$, defined by the following (differential) expression
\begin{align}\label{Premieredef}
H_{a,\t} =  H + A_{a,\t}\delta_{\partial\O}:=H + V_{a,\t},
\end{align} 
where $A_{a,\t}$ is a bounded invertible, self-adjoint operator in $\mathit{L}^{2}(\partial\O)^4$, which depends on parameters $(a,\t)\in\rr\times\rr^n$, $n\geqslant1$, and the $\delta$-potential is the Dirac distribution supported on $\partial\O$. From the physical point of view, the operator $H_{a,\t}$ describes the dynamics of the massive relativistic particles of spin-$1/2$ (i.e $H_{a,\t}$ is the Hamiltonian), in the external potential  $V_{a,\t}$. Due to their physical interest and the mathematical challenge that goes behind them, Dirac operator of the form \eqref{Premieredef} have been the subject of several mathematical studies. In particular, the case of bounded smooth domains, when $V_{a,\t}$ is a combination of \textit{electrostatic}, \textit{Lorentz scalar} and \textit{magnetic} $\delta$-potentials, which are given respectively by
\begin{align*}
V_\ep=\ep I_4\delta_{\partial\O},\quad  V_\m= \mu\beta\delta_{\partial\O}\quad  V_\e=\eta(\alpha\cdot N)\delta_{\partial\O},  \quad \ep,\m,\e\in\rr,
\end{align*} 
where $\mathit{N}$ is the (geometric measure theoretic) outward unit normal to $\O$; see e.g. \cite{AMV1,Albert,BEHL1,AP1,AP2,VR}. As it was observed in several works, for some particular values of the parameters the phenomenon of confinement arises, cf. \cite{AMV2,HOP,BEHL2}, and for some others values (which are called critical combinations of coupling constants) there is a loss of regularity in the operator domain of the underlying Dirac operators, cf. \cite{OBV,BH,BHOP,BB,CLMT}. More recently,  the following potentials
\begin{align*} 
V_{\ept}=  \ept\g_5\delta_{\partial\O}, \quad V_{\u}=  i\u\beta \left(\alpha\cdot N\right)\delta_{\partial\O},  \quad \ept,\u\in\rr,
\end{align*} 
(called here the \textit{modified electrostatic}  and the \textit{anomalous magnetic} $\delta$-interactions, respectively)  have been introduced and studied by the author in \cite{BB} (possibly for unbounded $\mathit{C}^2$-smooth domains), see also   \cite{CLMT} where the 2D analog of $V_{\u}$ was considered. In \cite{BB,CLMT} it was proved that the coupling $H+V_{\u}$ gathers the properties of confinement and the loss of regularity in the domain of the operator.  We mention that  when $A_{a,\t}$ is a $2\times 2$ matrix (i.e in the 2D case), a systematic spectral study was carried out in \cite{CLMT}. In there the authors considered the case of  $\delta$-interactions supported on a $\mathit{C}^\infty$-smooth closed curves. In other words, they recovered the most general local interactions given as $2\times 2$ Hermitian matrix (which can be written using the 2D analog of $V_\ep,V_\m, V_\e$ and $V_\u$). Unlike dimension two, the three-dimensional case is richer. Indeed, it turns out that in dimension three there are two other remarkable potentials, which are given by
\begin{align*} 
V_{\mt}=  i\mt\g_5\beta\delta_{\partial\O}, \quad V_{\ut}=  \ut\g_5\beta \left(\alpha\cdot N\right)\delta_{\partial\O},  \quad \mt,\ut\in\rr,
\end{align*} 
and that share almost the same properties as the Lorentz scalar and the anomalous magnetic $\delta$-interactions, respectively. For that reason, we call them the \textit{modified Lorentz scalar} and the  \textit{modified anomalous magnetic} $\delta$-interactions, respectively.  In order to study systematically the coupling of $H$ with the above potential we introduce the $\beta$ and the $\g$ transformations of the electrostatic and the magnetic $\delta$-interactions (cf. Subsection \ref{subb3.4}). 
 In fact, all the previous potentials can be obtained as the multiplication  (which keeps the symmetry) by $\beta$ and/or $\g_5$ of the electrostatic and the magnetic $\delta$-interactions.  With the help of these transformations, we classify in a  simple way the potentials which induce the confinement as well as their critical parameters.

A very important issue that arises  when we study such a coupling problems is the regularity of the surface $\partial\O$. In fact, to our knowledge,  all the works which deal with the three dimensional Dirac operators coupled with $\delta$-shell interactions have been done for $\O$ at least  $\mathit{C}^2$-smooth. In 2D,  the analogue of the coupling $(H+V_\m)$  was studied in \cite{PV},  when  $|\m|<2$, and  $\partial\O$ a closed curve with finitely many corners and otherwise smooth.  We also mention the works \cite{CL,LoicO}, where the authors studied the Dirac operator with infinite mass boundary conditions on sectors and wedges, respectively. It is worth mentioning that in the papers previously cited,  the techniques used depend significantly on the canonical identification $\rr^2\simeq \cc$, and the nature of the problem, for example in  \cite{CL,LoicO} the radial symmetry of the domain plays a main role. Clearly, if $\O$ is not sufficiently smooth (Lipschitz for example), then the strategies   of \cite{CL,LoicO,PV}  cannot be used to prove the self-adjointness of $H_\t$, especially in 3D.  Thus, we  ask the following questions:
\begin{itemize}
  \item[(Q1)] Until what extent the results on self-adjointness of $H_{a,\t}$ (at least when $A_{a,\t}$ is one of the usual $\delta$-interactions)  also hold for Lipschitz domains?
\end{itemize}
\begin{itemize}
  \item[(Q2)]   If $H_{a,\t}$ is self-adjoint, can we characterize its (essential/discrete) spectrum?
\end{itemize}
\begin{itemize}
  \item[(Q3)]  In the non-critical case, is it possible to characterize the dependence of the Sobolev regularity of the domain of $H_{a,\t}$ through the regularity of $\O$?
\end{itemize}
The main goal of the present paper, is to study questions  $(\mathrm{Q1})$,$(\mathrm{Q2})$ and $(\mathrm{Q3})$ under the weakest geometric assumptions on $\O$ (i.e possibly for non-Lipschitz domains).    However, we are not going to do a whole catalog, we will take the following operator
\begin{align*}
H_{\k} =  H + V_{\k}=H+(\ep I_4 +\mu\beta + \eta(\alpha\cdot N))\delta_{\S},\quad  \k:=(\ep,\m,\e)\in\rr^3,
\end{align*}
as a reference,  which was already studied by the author in \cite{BB}, for a critical and non-critical parameter when $\O_+$ is a $\mathit{C}^2$-smooth domain. Once we understand how to answer these questions for the operator $H_\k$, we will look at 
\begin{align}\label{famille1}
\begin{split}
H_{\mt} &=  H + V_{\ut}=H+i\mt\g_5\beta\delta_{\S},\quad \mt \in\rr,\\
H_{\ut} &=  H + V_{\ut}=H + i\ut\g_5\beta(\aN)\delta_{\S},  \quad \ut\in\rr,
\end{split}
\end{align}
and also   to the families of Dirac operators given by:
\begin{align}\label{famille2}
\begin{split}
 (-m,m)&\ni a\longmapsto H_{a,\l}= H+ \l \CS^{a}\delta_{\S},\quad\l\in\rr\setminus\{0\},\\
  (-m,m)&\ni a\longmapsto H_{a,\ll}= H+ \ll( \aN) \mathit{C}^{a}_{\S}( \aN)\delta_{\S},\quad \ll\in\rr\setminus\{0\},
  \end{split}
\end{align}
where $ \mathit{C}^{a}_{\S}$ is the Cauchy operator associated to $(H-a)$,  see Subsection \ref{Subb2.2}.  We will see that $ H_{\mt} $, $H_{\ut}$, $H_{a,\l}$ and $H_{a,\ll}$ induce confinement, with some new boundaries conditions.

Let us expose our set of techniques and summarize the main results concerning the operator $H_\k$.  Following the strategy in \cite{AMV1},  we define the Dirac operator $H_{a,\t}$ on the domain
 \begin{align*}
\mathrm{dom}(H_{a,\t})=\left\{ u+\Phi[g]: u\in\mathit{H}^1(\rr^3)^4, g\in\mathit{L}^2(\partial\O)^4, u_{|\partial\O}=-\Lambda_{\t,+}[g]\right\}, 
\end{align*}
where $\Phi$ is an appropriate fundamental solution of the unperturbed operator $H$, and  $\Lambda_{\t,\pm}$ are bounded linear operators  acting on $\mathit{L}^2(\partial\O)^4$  (see \ref{Lambda}).  So, unlike most existing works, we are not going to treat the  $\delta$-interactions as a transmission problem. There are several reasons prompt us to choose to work with the strategy of \cite{AMV1}, we cite here: 
\begin{itemize}
  \item One can give meaning to the trace on $\partial\O$ of functions in $\mathit{H}^1(\rr^3)$, for a very large class of surfaces (Ahlfors-David regular surfaces). If $\O$ is a non Lipschitz domain, then this becomes a real obstacle if we work with transmission conditions, see the discussion in the beginning of Subsection \ref{Subb3.3}.
  \item For non-critical combinations of coupling constants, the main result of \cite{AMV1} gives us a powerful tool  to prove  the self-adjointness  of $H_{a,\t}$. In fact, if $\Lambda_{\t,+}$ is Fredholm then $H_{a,\t}$ is self-adjoint, see Theorem \ref{luis}. 
  \item In order to study the Sobolev regularity of $\dom(H_{a,\t})$, it is sufficient to study the regularity properties of $\Phi[g]$.
\end{itemize}
 Fredholm's property  and (or) invertibility  of the boundary integral operator is one of the  important tools for the analysis of strongly elliptic boundary value problems;  such techniques have been exploited since a long time to solve for example the Dirichlet or Neumann problem on Lipschitz domains;  cf. \cite{JK},\cite{Verc} and \cite{DK}.  As we will see throughout the paper, the anticommutators $\{ \aN, \mathit{C}_{\partial\O}^a\}$ play a central role in our study. In fact, as it was observed in \cite{AMV1,AMV2,AMV3,BB}, the compactness on $\mathit{L}^2(\partial\O)^4$ of $\{ \aN, \mathit{C}_{\partial\O}\}$ (where  $\mathit{C}_{\partial\O}^0= \mathit{C}_{\partial\O}$) implies that $\Lambda_{\k,+}$ is a  Fredholm operator.   So we can naturally ask the following question:
\begin{itemize}
  \item[(Q4)] Given a bounded Lipschitz domain $\O$, what is the necessary regularity on $\partial\O$  so that the anticommutator $\lbrace \alpha\cdot \mathit{N}, \mathit{C}_{\partial\O}\rbrace$ gives rise to a  compact operator on  $\mathit{L}^2(\partial\O)^4$?
\end{itemize}
One of the main results of this article is the answer to this question, see Theorem \ref{qestion3}. Looking closely at the anticommutator $\lbrace \alpha\cdot \mathit{N}, \mathit{C}_{\partial\O}\rbrace$, we observe that it involves a matrix version of the principal value of the harmonic double-layer $K$, its adjoint $K^{\ast}$ and the commutators $[\mathit{N}_k, R_j]$, where $R_j$ are the Riesz transforms (see \eqref{pv double layer} for the definition). Hence the situation is more clear. In fact, from the harmonic analysis and geometric measure theory point of view, it is shown  that the boundedness of Riesz transforms characterizes the uniform rectifiability of $\partial\O$; cf. \cite{NVT}, for example. In addition,  functional analytic properties of the Riesz transforms (such as the identity $\sum^{3}_{j=1}R_j^2=-\mathit{I}$) and the analogue version of the strongly singular part of $ (\alpha\cdot \mathit{N}) \mathit{C}_{\partial\O}$ in the Clifford algebra $\mathcal{C}l_3$, i.e. the Cauchy-Clifford operator (especially its self-adjointness and compactness character) are strongly related to the regularity and geometric properties of the domain $\O$, for more details we refer to \cite{HMMPT} and \cite{HTM}. The most important fact which allow us to establish some results for the Lipschitz class, is that the compactness of $K$, $K^{\ast}$ and $[\mathit{N}_j, R_j]$  characterizes  the class of regular SKT (aka Semmes-Kenig-Toro) domains; see \cite{HTM}. However, regular SKT domains are not necessarily Lipschitz domains and vice versa.  So, to stay in the context of  compact Lipschitz domains,  we shall suppose that  $\O$  satisfies the following property: 
\begin{itemize}
\item[(H1)]  $\O$  is a bounded Lipschitz  domain with normal $\mathit{N}\in\mathrm{VMO}(\partial\O,\mathrm{dS})^3$.
  \end{itemize}
  This assumption characterizes the intersection of the Lipschitz class with the regular SKT class. Moreover,  the hypothesis $(\mathrm{H1})$ is the answer to  the question $(\mathrm{Q4})$. In fact, we prove the following:
  \begin{align}\label{VOMcha}
 \O \text{ satisfies the assumption }  (\mathrm{H1}) \Longleftrightarrow \lbrace \alpha\cdot \mathit{N}, \mathit{C}_{\partial\O}\rbrace \text{ is compact in }  \mathit{L}^2(\partial\O)^4.
\end{align}
   see Theorem \ref{qestion3}. Once we have established that and hence proved the compactness of $\lbrace \alpha\cdot \mathit{N}, \mathit{C}_{\partial\O}\rbrace$, the self-adjointness of $H_{\k}$ will be an easy consequence of \cite[Theorem 2.11]{AMV1}. Moreover,  we prove almost all the spectral properties as in the $\mathit{C}^2$-smooth case, see Theorem \ref{main4}.  Another geometric type result that we establish in this article, is a characterization of the class of regular SKT domains via the compactness  of the anticommutator $\lbrace \alpha\cdot \mathit{N}, \mathit{C}_{\partial\O}\rbrace $ in $ \mathit{L}^2(\partial\O)^4$, see Proposition \ref{remarkSKT}. More precisely,  using the material provided in \cite{HTM}, we show that if $\O$ is a two-sided NTA domain with a compact Ahlfors-David regular boundary, then it holds that 
 \begin{align}\label{SKTcha}
 \O \text{ is a regular SKT domain }  \Longleftrightarrow \lbrace \alpha\cdot \mathit{N}, \mathit{C}_{\partial\O}\rbrace \text{ is    compact in }  \mathit{L}^2(\partial\O)^4.
\end{align}
   At this stage, beyond the two classes of domains which are characterized by \eqref{VOMcha} and \eqref{SKTcha}, the compactness arguments mentioned previously are no longer valid. So, in order to go further in our study we change the strategy and we turn to the invertibility arguments, which are rather valid for a large class of domains. Indeed, we investigate the case of bounded uniformly rectifiable domains (see Section \ref{sec1} for the definitions). In one direction, making the assumption that  $0<|\ep^2-\m^2|<1/\| \mathit{C}_{\partial\O}\|^2_{\mathit{L}^2(\partial\O)^4\rightarrow\mathit{L}^2(\partial\O)^4} $, we then  show that  $(H+V_\ep+V_\m)$ is self-adjoint, cf. Theorem \ref{th5.6}. In another direction, assuming that 
   \begin{align*}
   \m^2>\ep^2\quad  \text{ or } \quad 16\| W\|^2_{\mathit{L}^2(\S)^2\rightarrow\mathit{L}^2(\S)^2}<\ep^2-\m^2<1/\| W\|^2_{\mathit{L}^2(\S)^2\rightarrow\mathit{L}^2(\S)^2},
   \end{align*}
    (here $W$ is the stroungly singular part of $\mathit{C}_{\partial\O}$ defined in \eqref{opeWW}), we also prove the self-adjointness of $(H+V_\ep+V_\m)$. In particular, if $\O$ is Lipschitz, we then recover the same spectral properties as in the case of the assumption $(\mathrm{H1})$.  Moreover, we show that $\mathcal{H}_{\ep,\m}$ generates confinement when $\ep^2-\m^2=-4$, cf. Theorem \ref{th5.7} and Proposition \ref{Un confinement}.

Having established the above results, and in order to give an answer to question  $(\mathrm{Q3})$ , we consider  the class of H\"{o}lder's domains $\mathit{C}^{1,\o}$, with $\o\in(0,1)$, and we prove that the functions in $\dom(H_{\k})$ have the $\mathit{H}^s$-Sobolev regularity, with $s>1/2$. In particular, we show that if $\o\in(1/2,1)$, then $\mathrm{dom}(H_{\k})\subset\mathit{H}^1(\rr^3\setminus\partial\O)^4$, cf. Theorem \ref{self gamma}. 

Finally, these compactness and invertibility arguments are then applied to the spectral study of the Dirac operators given by \eqref{famille1} and \eqref{famille2}.


\textbf{Structure of the paper.} The paper is organised as follows. In the nest section, we set up the necessary definition and give some auxiliary results related to the integral operators associated to $H$. Section \ref{sec3} and \ref{sec5} are the heart of the paper and  contain our most important contributions. In Section \ref{sec3}, we provide the necessary materials to tackle all the problems under consideration.  Section \ref{sec5}, is divided into four subsections as follows. Subsection \ref{sub5.1}, is devoted  to the study of $H_{\k}$, when $\O$ is a Lipschitz domain with normal in VMO. In there, the spectral properties of $H_\k$ are given for non-critical parameters, and we provide all the result about the compactness. Subsequently,  in subsection \ref{sub5.3}, we study the Sobolev regularity of $\dom(H_{\k})$, for  H\"{o}lder's domains $\mathit{C}^{1,\o}$. After this,   we consider in subsection \ref{sub5.2} the couplings  $(H+V_\ep+V_\m)$ and $(H+V_\e)$, in the case of bounded, uniformly rectifiable domains. Finally,   the spectral study of $H_{\mt}$ and $H_{\ut}$ is performed  in Subsection \ref{sub4.4}. The last  section of the paper is devoted to the spectral study of the families of Dirac operator $H_{a,\l}$ and $H_{a,\ll}$.

\section{Preliminaries} \label{sec1}
\setcounter{equation}{0}
 We begin by setting notations, and recalling some definitions   and results from geometric measure theory which are used in the current work. 
 \subsection{Notation and Definitions} We use the following notations:
\begin{itemize}
 \item For a Hilbert space $\mathfrak{h}$, we denote by $\mathcal{B}(\mathfrak{h})$ (respectively $\mathcal{K}(\mathfrak{h})$) the space of bounded  (respectively  bounded  and compact), everywhere defined linear operators in $\mathfrak{h}$.  If $T$ is a closed operator in  $\mathfrak{h}$ then its spectrum, essential spectrum, and discrete spectrum are denoted by $\mathrm{Sp}(T)$,  $\mathrm{Sp}_{\mathrm{ess}}(T)$, and $\mathrm{Sp}_{\mathrm{disc}}(T)$, respectively.
 \item For $A,\, B\in \mathcal{B}(\mathfrak{h})$,  we denote by $[A, B]$ (respectively $\{A, B\}$)  the usual commutator (respectively anticommutator) bracket.  
  \item We use the letter $C$ to denote harmless positive constant, not necessarily the same at each occurrence.
  \item We use the notation $\rr^3_\pm=\lbrace x\in\rr^3: \pm x_3>0\rbrace$ for the upper (respectively the lower) half space. Also, the upper/lower  complex half plane is denoted by $\cc_\pm$.
  \item By $\mathcal{H}^2$ we denote the $2$-dimensional Hausdorff measure, and  we let  $\mathrm{dS} =\mathrm{d}\mathcal{H}^2\big|_{\mathrm{E}}$ to be the surface measure on a closed set $\mathrm{E}\subset\rr^3$ of dimension $2$.
  \item We denote by $\mathrm{diam}(E)$ the  diameter of $E$, that is  $\mathrm{diam}(E):= \sup_{x,y\in E}|x-y|$.
  \item We denote by $B(x,r)$ the Euclidean ball of radius $r$ centred at $x\in\rr^3$.
  \item For a Borel set $B$ of $\rr^3$, the Lebesgue measure of $B$ is denoted by $|B|$.
  \item For a Borel measure $\mu$, and a Borel set $B$ with $0<\mu(B)<\infty$, we set  
  $$\oint_{B} U\mathrm{d}\mu:=\mu(B)^{-1}\int_{B} U\mathrm{d}\mu.$$
\end{itemize}

Throughout the current  paper $\O$ will be a connected open (proper) subset of $\rr^3$, and $\O^{\mathrm{c}}$ is the complement of $\O$.

\begin{definition}[Ahlfors-David regular] We say that a set $E\subset\rr^3$ is $2$-dimensional Ahlfors-David regular, or simply ADR, if it is closed and there is some uniform constant $C$ such that 
\begin{align}
   \frac{1}{C}r^2\leqslant \mathcal{H}^2( B(x,r)\cap E)\leqslant Cr^2,\quad \forall x\in E, r\in (0,\mathrm{diam}(E)).
\end{align}
\end{definition}
\begin{definition}[Uniformly rectifiable domains]\label{UR} We say that a compact set $E\subset\rr^3$ is uniformly rectifiable provided that it is ADR and the following holds. There exist $\rho$, $M\in(0,\infty)$ (called the UR character of $E$) such that for each $x\in E$, $r\in(0,1]$, there is a Lipschitz map $\phi: B_r\rightarrow\rr^3$ (where $B_r$ is a ball of radius $r$ in $\rr^2$) with Lipschitz constant $L_\phi\leqslant M$, such that 
\begin{align}
  \mathcal{H}^2( E\cap B(x,r)\cap \phi(B_r))\geqslant \rho r^2.
\end{align}
 A nonempty, proper and bounded open subset $\O\subset\rr^3$ is called uniformly rectifiable, or simply UR,  provided that $\partial\O$ is uniformly rectifiable and also $ \mathcal{H}^2(\partial\O\setminus\partial_{\ast}\O)=0$, where $\partial_{\ast}\O$ is the measure theoretic boundary of $\O$ defined by
 \begin{align}
 \partial_{\ast}\O:=\left\{ x\in\partial\O: \limsup\limits_{r\rightarrow0}\frac{|B(x,r)\cap\O|}{r^{3}}>0, \,\limsup\limits_{r\rightarrow0}\frac{|B(x,r)\cap\O^{\mathrm{c}}|}{r^{3}}>0\right\}.
 \end{align}
\end{definition}  
\begin{definition}[BMO and VMO] For an ADR set $E\subset\rr^3$,  $\mathrm{BMO}(E,\mathrm{dS})$ stands for the space of functions with bounded mean oscillation,  relative to the surface measure $\mathrm{dS}$. We denote by $\mathrm{VMO}(E,\mathrm{dS})$ the Sarason space of functions with vanishing mean oscillation on $E$, i.e the closure of
the set of bounded uniformly continuous functions defined on $E$ in $\mathrm{BMO}(E,\mathrm{dS})$.
\end{definition}
From now, unless stated otherwise, we always assume that $\O$ is a UR domain with $\partial\O=\partial\overline{\O}$, and we set 
\begin{align}
 \O_+=\O \text{ and }  \O_-:=\rr^3\setminus\overline{\O}.
 \end{align}
 As a consequence,  $\O_-$ is also a UR domain with the same ADR boundary as $\O_+$ (cf. \cite[Proposition 3.10]{HTM}), which we denote by $\S:=\partial\O_+=\partial\O_-$.
\begin{definition} Let $\O_\pm$ be as above. Fix $a>0$ and let $x\in\S$, then the nontangential approach regions of opening $a$ at the point $x$ is defined by
\begin{align}\label{ntregion}
\Gamma^{\O_\pm}(x)=\Gamma_a^{\O_\pm}(x)=\{ y\in\O_\pm: |x-y|<(1+a)\mathrm{dist}(y,\S)\}.
\end{align}
If $x\in\S$ and $U:\O_\pm\rightarrow\cc^4$, then 
\begin{align}\label{ntTrace}
U_\pm(x):= \lim\limits_{\Gamma^{\O_\pm}(x) \ni y\xrightarrow[]{}x}U(y),
\end{align}
is the nontangential limit of $U$ with respect to $\O_\pm$ at $x$.  
We also define the nontangential maximal function of $U$ on $\S$ by
\begin{align}\label{ntMaximali}
\mathcal{N}^{\O_\pm}[U](x)=\mathcal{N}_a^{\O_\pm}[U](x)=\sup\{ |U(y)|: y\in\Gamma^{\O_\pm}(x)\},\quad x\in\S,
\end{align}
with the convention that $\mathcal{N}^{\O_\pm}[U](x)=0$ when $\Gamma^{\O_\pm}(x)=\emptyset$.
Given  $g\in\mathit{L}^2(\S)$, we define the Hardy-Littlewood maximal operator by 
\begin{align}\label{Maximalfunction}
M^{\S}g(x)=\sup_{r>0}\oint_{B(x,r)\cap\S} |g(y)|\mathrm{dS},\quad x\in\S.
\end{align}
Then, by \cite[p. 624]{CW}, there is $C>0$ such that 
\begin{align}\label{FeW}
\| M^{\S}g\|_{\mathit{L}^2(\S)} \leqslant C \| g\|_{\mathit{L}^2(\S)}.
\end{align}
\end{definition}

\subsubsection{Sobolev and Besov spaces}\label{BesovSobolev}

Throughout the paper,  we shall  work on the Hilbert space $\mathit{L}^2(\rr^{3})^4$ with respect to the Lebesgue measure. $\mathcal{D}(\O_\pm)^4$ denotes the usual space of indefinitely differentiable functions with compact support, and $\mathcal{D}^{\prime}(\O_\pm)^4$ is the  space of distributions defined as the dual space of $\mathcal{D}(\O_\pm)^4$. By $\mathit{H}^s(\rr^{3})^4$  we denote the first order $\mathit{L}^2$-based Sobolev space over $\rr^3$. The Sobolev space $\mathit{H}^1(\O_\pm)^4$ is defined as follows:
\begin{align}
\mathit{H}^1(\O_\pm)^4=\{ \varphi\in\mathit{L}^2(\O_\pm)^4: \text{ there exists } \tilde{\varphi}\in\mathit{H}^1(\rr^{3})^4 \text{ such that }  \tilde{\varphi}|_{\O_\pm} =\varphi\}.
\end{align}
By $\mathit{L}^2(\S, \mathrm{d}S)^4:=\mathit{L}^2(\S)^4$ we denote the usual $\mathit{L}^2$-space over $\S$.   If $\S$ is Lipschitz, we then define the Sobolev spaces  $\mathit{H}^s(\S)^4$, $s\in[-1,1]$, using local coordinates representation on the surface $\S$; cf. \cite{Mc}. Also, if $\S$ is of class $\mathit{C}^{0,\g}$ and $s\in(0,\g)$, then one can define equivalently the Sobolev space  $\mathit{H}^{s}(\S)^4$ (see \cite[Chapter 4]{HW} for example) using the following norm
 \begin{align}
 \|g \|^2_{\mathit{H}^{s}(\S)^4}:= \int_{\S}|g(x)|^2dS(x) + \int_{\S}\int_{\S}\frac{|g(x)-g(y)|^2}{|x-y|^{2(1+s)}}\mathrm{dS}(y) \mathrm{dS}(x).
  \end{align}
Again if $\S$ is Lipschitz, then we denote by $t_{\S}:\mathit{H}^1(\O_\pm)^4\rightarrow \mathit{H}^{1/2}(\S)^4$  the classical trace operator, and by   $E_{\O_\pm}:  \mathit{H}^{1/2}(\S)^4\rightarrow  \mathit{H}^1(\O_\pm)^4$ the extension operator, i.e $t_{\S}E_{\O_\pm}$ is the identity operator.  For a function $u\in\mathit{H}^1(\rr^3)^4$, with a slight abuse of terminology we will refer to $t_{\S}u$ as the restriction of $u$ on $\S$. 

We shall use a trace theorem for  functions in the Sobolev space $\mathit{H}^1(\rr^3)$ in the case of ADR surfaces.  Assume that $\S$ is ADR, then the Besov space $\mathit{B}^2_{1/2}(\S)$ (see \cite[Chapter V]{JW} for example), consists of all functions $g\in\mathit{L}^2(\S)$ for which 
 \begin{align}\label{Besov}
 \int\int_{|x-y|<1}\frac{|g(x)-g(y)|^2}{|x-y|^{3}}\mathrm{dS}(y) \mathrm{dS}(x)<\infty.
 \end{align}
The Besov space is equipped with the norm
 \begin{align}\label{Besov norm}
 \|g \|^2_{\mathit{B}^2_{1/2}(\S)}:= \int_{\S}|g(x)|^2dS(x) + \int\int_{|x-y|<1}\frac{|g(x)-g(y)|^2}{|x-y|^{3}}\mathrm{dS}(y) \mathrm{dS}(x).
 \end{align}
 Given $U\in\mathit{H}^1(\rr^3)$,  set
 \begin{align}\label{TTrace}
 T_{\S}u(x):=\lim\limits_{r\searrow 0}\oint_{B(x,r)}U(y)\mathrm{d}y,
 \end{align}
 at every point $x\in\S$ where the limit exists. Then, we have the following trace theorem, for  the proof we refer to \cite[Theorem 1 and Example 1]{AJ} and \cite[Theorem 1, p.182 ]{JW}. 
 \begin{proposition}\label{Les traces}  Let $\S$ be as above. Then the trace operator $t_\S$ (which until now was defined on $\mathcal{D}(\rr^3)$) extend to a bounded linear operator $T_\S$ from $\mathit{H}^1(\rr^3)$ to $\mathit{B}^2_{1/2}(\S)$ (where $T_\S$ is given by \eqref{TTrace}) with a bounded linear inverse operator $\mathcal{E}$ from $\mathit{B}^2_{1/2}(\S)$ to $\mathit{H}^1(\rr^3)$. In other words,  $\mathit{B}^2_{1/2}(\S)$ is the trace to $\S$ of $\mathit{H}^1(\rr^3)$ and $T_\S\mathcal{E}$ is the identity operator.
  \end{proposition}
  \begin{remark}\label{rktrace}We often use the fact that the trace operator $T_\S$ coincide with $t_\S$, when $\O_+$ is a Lipschitz domain (i.e $\mathit{H}^{1/2}(\S)$ is the trace to $\S$ of $\mathit{H}^1(\rr^3)$).
\end{remark}

\subsection{Integral operators associated to $H$, Hardy spaces and Calder\'on's decomposition}\label{Subb2.2} 

In this part, we recall some known results on the integral operators associated to the fundamental solution of the Dirac operator. Recall that $(H,\mathit{H}^1(\rr^3)^4)$ is self-adjoint (cf.  \cite{Tha}) and its spectrum is given by
\begin{align*}
  \mathrm{Sp}(H) =\mathrm{Sp}_{\mathrm{ess}}(H)=(-\infty,-m]\cup [m,+\infty).
\end{align*}
 Given $z\in\cc\setminus\left((-\infty,-m]\cup[m,\infty)\right)$  with the convention that $\mathrm{Im}\sqrt{z^2-m^2}>0$, then the fundamental solution of $(H-z)$ is given by
\begin{align*}\label{}
  \phi^z(x)=\frac{e^{i\sqrt{z^2-m^2}|x|}}{4\pi|x|}\left(z +m\beta+( 1-i\sqrt{z^2-m^2}|x|)i\alpha\cdot\frac{x}{|x|^2}\right), \quad \text{for all } x\in\rr^3\setminus\{0\}.
\end{align*}  
Let $g \in\mathit{L}^2(\S)^4$, we define the following operators
\begin{align}\label{Phi}
\begin{split}
 \Phi^{z}[g](x) &=\int_\S \phi^{z}(x-y)g(y)\mathrm{dS}(y), \quad\text{ for all } x\in\rr^3\setminus\S,\\
  \mathit{C}^{z}_{\S}[g](x)&=  \lim\limits_{\rho\searrow 0}\int_{|x-y|>\rho}\phi^{z}(x-y)g(y)\mathrm{dS}(y), \quad \text{ for all } x\in\S,\\
  \mathit{C}^{z}_{\pm}[g](x)&=  \lim\limits_{\Gamma^{\O_\pm}(x) \ni y\xrightarrow[]{}x}\Phi^{z}[g](y), \quad \text{ for all } x\in\S.
\end{split}
\end{align}
Thanks to \cite[Lemma 2.1]{AMV1}, we know that $\Phi^z: \mathit{L}^2(\S)^4  \longrightarrow \mathit{L}^2(\rr^3)^4$ is well defined and bounded.

 Next, we denote by $\tilde{\phi}$ the fundamental solution of the massless Dirac operator $-i\sigma\cdot\nabla$, that is 
\begin{align}
\tilde{\phi}(x)= i\sigma\cdot\frac{x}{|x|^3}, \quad\text{ for all } x\in\rr^3\setminus\{0\},
\end{align}
and we define the bounded operator $\tilde{\Phi}: \mathit{L}^2(\S)^2  \longrightarrow \mathit{L}^2(\rr^3)^2$  as follows
\begin{align*}
  \tilde{\Phi}[h](x) =\int_\S \tilde{\phi}(x-y)h(y)\mathrm{dS}(y), \quad \text{for all } x\in\rr^3\setminus\S \text{ and } \forall h\in\mathit{L}^2(\S)^2.
\end{align*}
 Also, for  $x\in\S$  and  $h\in\mathit{L}^2(\S)^2$, we set  
\begin{align}\label{opeWW}
\begin{split}
W_{\pm}[h](x)&=\lim\limits_{\Gamma^{\O_\pm}(x) \ni y\xrightarrow[]{}x} \tilde{\Phi}[h](y),\\
 W[h](x)&=  \lim\limits_{\rho\searrow 0} \int_{|x-y|>\rho} \tilde{\phi} (x-y)h(y)\mathrm{dS}(y).
\end{split}
\end{align} 
Then, the following lemma gives us the relations between the operators defined above, and gathers their important properties.  We mention that when $\O_+$ is a bounded Lipschitz domain these results are well known, in this case we refer to \cite[Lemma 3.3]{AMV1} for example.  In the case of UR domains the lemma is somehow contained in \cite{HTM},  but for the convenience of the reader   we give here the main ideas to establish it.

\begin{lemma}\label{lemme 2.1} Given $z\in\cc\setminus\left((-\infty,-m]\cup[m,\infty)\right)$. Let $\mathit{C}^{z}_{\S}$, $\mathit{C}^{z}_{\pm}$, $W_\pm$ and $W$ be as above. Then $ \mathit{C}^{z}_{\S}[g](x)$, $\mathit{C}^{z}_{\pm}[g](x)$,  $W_\pm[h](x)$ and $W[h](x)$ exist for $\mathrm{dS}$-a.e. $x\in\S$,  $\mathit{C}^{z}_{\S}, \, \mathit{C}^{z}_{\pm}\in\mathcal{B}( \mathit{L}^2(\S)^4)$ and $W, \,  W_{\pm}\in\mathcal{B}( \mathit{L}^2(\S)^2)$. Furthermore, the following hold true:
\begin{itemize}
  \item[(i)]  $W_{\pm}=  \mp\frac{i}{2}(\sigma\cdot\mathit{N})+ W$.
     \item[(ii)]   $\mathit{C}^{z}_{\pm}= \mp\frac{i}{2}(\alpha\cdot\mathit{N}) + \mathit{C}^{z}_{\S} $.
  \item[(iii)] $((\sigma\cdot\mathit{N})W)^2=(W (\sigma\cdot\mathit{N}) )^2 =-\frac{1}{4}\mathit{I}_2$. In particular, we have $\Vert W\Vert\geqslant \frac{1}{2}$. 
  \item[(iv)] $( (\alpha\cdot\mathit{N})\mathit{C}^{z}_{\S} )^2 =(\mathit{C}^{z}_{\S} (\alpha\cdot\mathit{N}) )^2 =-\frac{1}{4}\mathit{I}_4$. In particular, we have $\Vert \mathit{C}^{z}_{\S}\Vert\geqslant \frac{1}{2}$.
\end{itemize}
\end{lemma}
\textbf{Proof.}  Given $f\in \mathit{L}^2(\S)$, thanks to \cite[Proposition 3.30]{HTM} we know that for each $j\in\{1,2,3\}$, the limit
\begin{align}
    \lim\limits_{\rho\searrow 0}\int_{|x-y|>\rho}\frac{x_j-y_j}{4\pi |x-y|^3}f(y)\mathrm{dS}(y),
\end{align}
exist at almost every $x\in\S$. Moreover, it holds that 
\begin{align}\label{jusmprelation}
   \lim\limits_{\Gamma^{\O_\pm}(x)  \ni w\xrightarrow[]{}x}\int\frac{w_j-y_j}{4\pi |w-y|^3}f(y)\mathrm{dS}(y)= \mp\frac{1}{2}\mathit{N}_j(x)f(x) +  \lim\limits_{\rho\searrow 0}\int_{|x-y|>\rho}\frac{x_j-y_j}{4\pi |x-y|^3}f(y)\mathrm{dS}(y).
\end{align}
Thus, working component by component it follows that $W_{\pm}[h](x)$ and $W[h](x)$ exist for $\mathrm{dS}$-a.e. $x\in\S$, and $W,\,W_\pm\in\mathcal{B}( \mathit{L}^2(\S)^2)$. Item $(\mathrm{i})$ follows by applying  the jump relation \eqref{jusmprelation} to the functions $\sigma_j h$, $j=1,2,3$.

 Now, we are going to show $(\mathrm{ii})$  and complete the proof of the first statement. For that, fix
  $z\in\cc\setminus\left((-\infty,-m]\cup[m,\infty)\right)$ and set 
\begin{align}\label{derivsingl}
  k(x):=\phi^z(x)-i(\alpha\cdot\frac{x}{|x|^3}), \quad \text{for all } x\in\rr^3\setminus\{0\}.
\end{align}
Then there is a constant $C$ such that $|k(\omega,y)|\leqslant C/|\omega-y|^{3/2}:=\tilde{k}(\omega,y)$, for $\omega,y\in\overline{\O}_{+}$. Define
\begin{align}\label{Tilde}
T[g](x)=\int_{\S}\tilde{k}(x,y)g(y)\mathrm{dS}(y).
\end{align}
Clearly, $T$ is bounded in $\mathit{L}^2(\S)^4$. Now, recall the definition of $\Gamma^{\O_\pm}(x)$ from \eqref{ntregion}. Let $x\in\S$ and $\omega\in\Gamma^{\O_+}(x)$, then 
\begin{align}
\left| \int_{B(\omega,2|x-\omega|)\cap\S}k(\omega,y)g(y)\mathrm{dS}(y)\right| \leqslant \int_{B(\omega,2|x-\omega|)\cap\S}C\left(\frac{1+a}{|x-\omega|}\right)^{3/2}|g(y)|\mathrm{dS}(y).
\end{align}
Using the ADR property of $\S$ (more precisely, use the inequality $\mathcal{H}^2( B(x,r)\cap \S)\leqslant Cr^2$), it follows that there is $C_1$ depending only on the ADR  constant of $\S$ such that 
\begin{align}\label{KK1}
\left| \int_{B(\omega,2|x-\omega|)\cap\S}k(\omega,y)g(y)\mathrm{dS}(y)\right| \leqslant C_1  |x-\omega|^{1/2}M^{\S}g(x), 
\end{align}
where $M^{\S}$ is the  Hardy-Littlewood maximal operator defined by \eqref{Maximalfunction}. Now, let $y\in\S\setminus B(x,2|x-\omega|)$, then $|\omega -y|\leqslant 2|x-y|$, and thus $|k(\omega,y)|\leqslant \tilde{k}(\omega,y)\leqslant 2^3\tilde{k}(x,y)$. Therefore, we get 
\begin{align}\label{KK2}
\left| \int_{\S\setminus B(\omega,2|x-\omega|)}k(\omega,y)g(y)\mathrm{dS}(y)\right| \leqslant 2^3 T[|g|](x). 
\end{align}
Thus,  \eqref{KK1}, \eqref{KK2}  and the dominate convergence theorem yield that 
\begin{align}\label{limnont+}
\lim\limits_{\Gamma^{\O_+}(x) \ni \omega\xrightarrow[]{}x}\int_{\S}k(\omega,y)g(y)\mathrm{dS}(y)= \int_{\S}k(x,y)g(y)\mathrm{dS}(y),
\end{align}
holds for all $g\in\mathit{L}^{2}(\S)^4$ and $\mathrm{dS}$-a.e. $x\in\S$. Similarly, one can show that  
\begin{align}\label{limnont-}
\lim\limits_{\Gamma^{\O_-}(x) \ni \omega\xrightarrow[]{}x}\int_{\S}k(\omega,y)g(y)\mathrm{dS}(y)= \int_{\S}k(x,y)g(y)\mathrm{dS}(y),
\end{align}
holds for all $g\in\mathit{L}^{2}(\S)^4$ and for $\mathrm{dS}$-a.e. $x\in\S$. Thus, given any $g\in\mathit{L}^2(\S)^4$, it follows from  the above considerations and \eqref{derivsingl} that  $ \mathit{C}^{z}_{\S}[g](x)$ and $\mathit{C}^{z}_{\pm}[g](x)$ exist for $\mathrm{dS}$-a.e. $x\in\S$, and $\mathit{C}^{z}_{\S}, \, \mathit{C}^{z}_{\pm}\in \mathcal{B}( \mathit{L}^2(\S)^4)$. Now, using  \eqref{limnont+}, \eqref{limnont-} and $(\mathrm{i})$ (i.e working component by component) we easily get  $(\mathrm{ii})$. 

 Finally the proof of $(\mathrm{iii})$ and  $(\mathrm{iv})$  is a relatively straightforward modification of the technique used in the proof of \cite[Lemma 3.3]{AMV1}$(\mathrm{ii})$. Indeed,  by \cite[p. 2659]{HTM} it follows that 
   \begin{align}\label{besoinHTM}
   \begin{split}
\|\mathcal{N}_a^{\O_\pm}[\tilde{\Phi}[h]]\|_{\mathit{L}^2(\S)^2} &\leqslant C \| h\|_{\mathit{L}^2(\S)^2},\\
\|\mathcal{N}_a^{\O_\pm}[\Phi^z[g]]\|_{\mathit{L}^2(\S)^4}& \leqslant C \| g\|_{\mathit{L}^2(\S)^4},
\end{split}
\end{align}
for some $C>0$ depending only on $a$ as well as the ADR  and the UR constants of $\S$. Now, observe that 
\begin{align}
(-i\sigma\cdot\nabla)\tilde{\Phi}[h]=0 \text{ and } (H-z)\Phi^{z}[g]=0 \text{ in } \Omega_\pm.
\end{align}
Then, by \cite[Theorem 4.49]{HTM} it holds that 
\begin{align}\label{reproducingformula}
\begin{split}
\tilde{\Phi}[h]&= \int_{\S}\tilde{\phi}(x-y)(\pm i\sigma\cdot\mathit{N}(y))h(y)\mathrm{dS}(y), \quad x\in\O_\pm,\\
\Phi^z[g]&= \int_{\S}\phi^z(x-y)(\pm i\alpha\cdot\mathit{N}(y))g(y)\mathrm{dS}(y), \quad x\in\O_\pm.
\end{split}
\end{align}
Although \cite[Theorem 4.49]{HTM} was stated in the case of tow-sided NTA domains (cf. Definition \ref{NTAD}) it also holds for UR domains by the discussion on  \cite[p. 2758]{HTM}.  Now, given $x\in\O_+$, $h\in\mathit{L}^2(\S)^2$ and $g\in\mathit{L}^2(\S)^4$. Then, $(\mathrm{i})$ (respectively   $(\mathrm{ii})$) and \eqref{reproducingformula} yield that 
\begin{align}\label{KK4}
\begin{split}
\tilde{\Phi}[( i\sigma\cdot\mathit{N})h](x)&= \tilde{\Phi}[( i\sigma\cdot\mathit{N})W_+( i\sigma\cdot\mathit{N})h](x), \\
\Phi^z[( i\alpha\cdot\mathit{N})g](x)&=\Phi^z[( i\alpha\cdot\mathit{N})\mathit{C}^z_+( i\sigma\cdot\mathit{N})g](x).
\end{split}
\end{align}
Therefore, $(\mathrm{iii})$  (respectively $(\mathrm{iv})$) follows by taking the nontangential limit in \eqref{KK4}) and using $(\mathrm{i})$ (respectively  $(\mathrm{ii})$).  This completes the proof of the lemma.
  \qed
  \newline
\begin{remark}\label{remark2.1} Note that since $\overline{\phi^{z}(y-x)}=\phi^{\overline{z}}(x-y)$, it follows that $(\mathit{C}^{z}_{\S})^{\ast}= \mathit{C}^{\overline{z}}_{\S}$ in $\mathit{L}^2(\S)^4$. In particular,     $\mathit{C}^{z}_{\S}$ ($W$)  is self-adjoint operators in $\mathit{L}^2(\S)^4$ (respectively in $\mathit{L}^2(\S)^2$), for all $z\in(-m,m)$.
\end{remark}
In order to understand better Lemma \ref{lemme 2.1}, we need to investigate the following class of domains.
\begin{definition}[two-sided NTA domains]\label{NTAD} Following \cite{JK1}, we say that  a nonempty, proper open set $\O$ of $\rr^3$  is an  NTA (non-tangentially accessible) domain if $\O$ satisfies both the two-sided Corkscrew and Harnack Chain conditions\footnote{Generally speaking,  the Corkscrew condition is a quantitative, scale invariant version of
openness, and the Harnack Chain condition is a scale invariant version of path connectedness.} (see \cite{JK1},\cite{HTM} or \cite[Appendix]{HMMPT} for the precise definition). Furthermore, we say that $\O$ is a two-sided NTA domain if both $\O$ and $\rr^3\setminus\overline{\O}$ are non-tangentially
accessible domains.
\end{definition}

 Assume that $\O_+$ is a two-sided NTA domain with an ADR boundary\footnote{In the literature, a two-sided NTA domain whose boundary is ADR  is often referred to as  a 2-sided Chord arc domain.} (which makes it a UR domain). Following \cite{HMMPT}, we define the Hardy spaces $\mathbb{H}^2_{z}(\O_\pm)^4$ by
\begin{align*}
\mathbb{H}^2_{z}(\O_+)^4=\left\{ u:\O_+\rightarrow \cc^4:\mathcal{N}[u]\in  \mathit{L}^2(\S)^4 \text{ and }(\mathcal{H}-z)u=0\right\},
\end{align*}
and 
\begin{align*}
\mathbb{H}^2_{z}(\O_-)^4=\big\{ u:\O_-\rightarrow \cc^4:\mathcal{N}[u]\in  &\mathit{L}^2(\S)^4,\, (\mathcal{H}-z)u=0 \\
&\text{ and } u(x)=\mathcal{O}(|x|^{-2}) \text{ as } |x|\rightarrow \infty\big\}.
\end{align*}
Then, from \eqref{besoinHTM} it follow that 
\begin{align*}
\Phi^z\downharpoonright_{\O_\pm}\in\mathbb{H}^2_z(\O_\pm)^4.
\end{align*}
 Now, the boundary Hardy spaces are defined as follows 
\begin{align*}
\mathbb{H}^2_{z,\pm}(\S)^4=\left\{ u\downharpoonright_{\S}: u\in\mathbb{H}^2(\O_\pm)^4 \right\},
\end{align*}
where the boundary trace is taken in a nontangential pointwise sense. Then, we have the following proposition, for the proof we refer to \cite[Subsection 2.3]{HMMPT}.
\begin{proposition}\label{Calderonprojectors} Let $ z\in\cc\setminus\left((-\infty,-m]\cup[m,\infty)\right)$, then the following decomposition  holds
\begin{align*}
\mathit{L}^2(\S)^4&=\mathbb{H}^2_{z,+}(\S)^4 \oplus\mathbb{H}^2_{z,-}(\S)^4.
\end{align*}
Moreover, it holds that 
\begin{align*}
\mathrm{Rn}\left(\frac{1}{2}+i(\alpha\cdot\mathit{N}) \mathit{C}^{z}_{\S} \right)&=\mathbb{H}^2_{z,+}(\S)^4= \mathrm{Kr}\left(-\frac{1}{2}+i(\alpha\cdot\mathit{N}) \mathit{C}^{z}_{\S}\right) ,\\
\mathrm{Rn}\left(-\frac{1}{2}+i(\alpha\cdot\mathit{N}) \mathit{C}^{z}_{\S} \right)&=\mathbb{H}^2_{z,-}(\S)^4= \mathrm{Kr}\left(+\frac{1}{2}+i(\alpha\cdot\mathit{N}) \mathit{C}^{z}_{\S}\right).
\end{align*}
In other words,  $\left(\frac{1}{2}\pm i(\alpha\cdot\mathit{N}) \mathit{C}^{z}_{\S} \right)$ is the Calder\'on's projector associated to  $\mathbb{H}^2_{z,\pm}(\S)^4$.
\end{proposition}
\begin{remark}\label{remarproj}From the above proposition, we conclude that $\mathit{C}^z_\pm[i(\alpha\cdot\mathit{N})g]\in\mathbb{H}^2_{z,\pm}(\S)^4$, for all $g\in\mathit{L}^2(\S)^4$. Thus, $\mathit{C}^z_\pm$ are projectors. Note that $\left(1/2\pm i \mathit{C}^{z}_{\S}(\alpha\cdot\mathit{N}) \right)$ are also projectors. This observation, is the main idea behind the Dirac operators considered in Section \ref{Sec 8}.
\end{remark}

 \section{Main Tools}\label{Main tools}\label{sec3}
 \setcounter{equation}{0}

 In this section, we gather the main tools to tackle the different problems that we are going to consider.  
Before going any further, let us define the general form of the operators we are interested on and explain the philosophy of our technique. 

 Let $A_\t:  \mathit{L}^2(\S)^4\longrightarrow\mathit{L}^2(\S)^4$ be a bounded invertible, and self-adjoint operator depending on a parameter  $\t\in\rr^n$ with $n\in\mathbb{N}^{\ast}$. We assume that the inverse of $A_\t$ is given explicitly by
 \begin{align}
 A_\t^{-1}= \frac{1}{\mathrm{sgn}(\t)}\tilde{A}_\t,
 \end{align}
 where  $ \mathrm{sgn}(\t)$ is a real number, defined by $ \mathrm{sgn}(\t)\mathit{I}_4= \tilde{A}_\t  A_\t= A_\t\tilde{A}_\t$. Except for some particular situations, we always deal with the case $ \mathrm{sgn}(\t)\neq0$. In the case $ \mathrm{sgn}(\t)=0$, we need to change slightly our definition to treat such a situation.

 Next, we define the operators $\Lambda^{z}_{\t,\pm}$ as follows:
\begin{align}\label{Lambda}
\Lambda^{z}_{\t,\pm}=A_\t^{-1}\pm\mathit{C}^{z}_{\S},\quad\forall z\in\cc\setminus\left((-\infty,-m]\cup[m,\infty)\right).
\end{align}
Clearly,  $\Lambda^{z}_{\t,\pm}$ are bounded (and self-adjoint for $z\in(-m,m)$) from  $\mathit{L}^{2}(\S)^4$ onto itself. We wish to note that for technical reasons, some times one need to change slightly the definition of the operator $\Lambda^{z}_{\t,-}$.  In the sequel, we shall write $\Phi$, $ \mathit{C}_{\S} $,  $\mathit{C}_{\pm}$ and $\Lambda^{}_{\t,\pm}$ instead of $\Phi^0$, $ \mathit{C}^{0}_{\S} $,  $ \mathit{C}^{0}_{\pm}$ and $\Lambda^{0}_{\t,\pm}$.

Now, we define the perturbed Dirac operator $H_{\t}$ acting in $ \mathit{L}^2(\rr^3)^4$,  by  
 \begin{align}
 H_{\t}=  H + V_{\t}= H+ A_\t\delta_{\S},
  \end{align}
on the domain
 \begin{align}\label{dom1}
\mathrm{dom}(H_{\t})=\left\{ \varphi=u+\Phi[g]: u\in\mathit{H}^1(\rr^3)^4, g\in\mathit{L}^2(\S)^4, T_{\S}u=-\Lambda_{\t,+}[g]\right\},
\end{align}
where $T_{\S}$ is the trace operator defined in Proposition \ref{Les traces}, and 
\begin{align}
V_\t(\varphi)=\frac{1}{2}A_\t(\varphi_+ +\varphi_-)\delta_{\S},
\end{align}
with $\varphi_\pm= T_\S u + \mathit{C}^{}_\pm [g]$. Thus, $H_{\t}$ acts in the sens of distributions as $H_{\t}(\varphi)= H(u)$, for all $\varphi=u+\Phi^{}[g]\in\mathrm{dom}(H_{\t})$. 

For some particular values of $ \mathrm{sgn}(\t)$, two interesting phenomena appear in the spectral study of the Dirac operator $H_\t$. The first one is the confinement phenomenon,  assuming that $H_\t$ is essentially self-adjoint, this means that for any datum $\varphi_0\in\mathrm{dom}(\overline{H_\t})$ with support in $\O_\pm$, the unique solution $\varphi\in\mathit{C}^1(\rr,\mathit{L}^2(\rr^3)^4)$ of the following Cauchy problem
 \begin{equation}\label{param}
   \left\{
\begin{aligned}
i\partial_t\varphi(t,x)&= \overline{\mathcal{H}}\varphi(t,x),\\
\varphi(0,x)&=\varphi_0(x),
\end{aligned}
  \right.
\end{equation}
remains for all times supported in $\O_\pm$. More concretely, if one let 
$$\mathit{L}^2(\rr^3)^p\cong\mathit{L}^2(\Omega_+)^4\oplus\mathit{L}^2(\Omega_-)^4.$$
Then $\overline{H_\t}$ decouples as follows 
\begin{align}\label{CONF}
\overline{H_\t}=H_\t^{\Omega_+}\oplus H_\t^{\Omega_-},
\end{align}
 where $H_\t^{\Omega_\pm}$ are self-adjoint Dirac operators  acting in $\O_\pm$ with some boundary conditions. Moreover,  the propagator satisfy
\begin{align*}
e^{-it\overline{H}_\t}=e^{-itH_\t^{\Omega_+}}\oplus e^{-itH_\t^{\Omega_-}}.
\end{align*}
In this article, we always use the characterization \eqref{CONF}, in such a case we say that $\overline{H}_\t$ (or $V_\t$) generates confinement or equivalently $\S$ is impenetrable.

The second one is called critical combinations of the coupling constants, it results in the loss of the Sobolev regularity of functions in the domain of $H_\t$  for smooth domain $\O_+$, i.e $\dom(\overline{H_\t})\not\subset\mathit{H}^s(\rr^3\setminus\S)^4$, for all $s>0$. To our knowledge, there is no fixed definition for such a case, because the spectral study of $H_\t$ depends significantly on the smoothness of the domain $\O_+$. As it was remarked by the author in \cite{BB}, it seems that the $\mathit{C}^2$-smoothness condition on $\O_+$ is necessary to prove the self-adjointness in $\mathit{L}^2(\rr^3)^4$ of $H_\t$ in such a case. Thus, for our applications we fix the definition of the critical combinations of the coupling constants as follows:
\begin{definition}[Critical parameters]\label{defcritic} Let $\O_+$ be  a bounded $\mathit{C}^2$-smooth domain and  let $H_\t$ be as in \eqref{dom1}. We say that the parameter  $\t\in\rr^n$, $n\in\mathbb{N}^{\ast}$, is critical if  $ \Lambda_{\t,+}\in\mathcal{K}(\mathit{L}^2(\S)^4)$ or $\Lambda_{\t,\mp} \Lambda_{\t,\pm}\in\mathcal{K}(\mathit{L}^2(\S)^4)$.
\end{definition}
In what follows, we shall use the phrase  \textit{"characteristics of the model"}, or simply  \textit{"characteristics of the potential"},  to refer to the above phenomena.


\subsection{Non-critical parameters}
As we mentioned in the introduction,   the self-adjointness (in the non-critical cases) of the Dirac operator $ H_{\t}$ will be derived using the main result of \cite{AMV1}. However, the way that \cite[Theorem 2.11]{AMV1} was stated does not take into account the case when $\O_+$ is UR. Therefore, a few comments on how to extend it should be in order. In fact, using Proposition \ref{Les traces} instead of \cite[Proposition 2.6]{AMV1}, and taking into account  \cite[Lemma 2.8 and Lemma 2.10]{AMV1}  and \cite[Remark 2.12]{AMV1},  the main result of \cite{AMV1} (more precisely  \cite[Theorem 2.11 $(\mathrm{iii})$]{AMV1}) reads as follows: 
  \begin{theorem}\label{luis} Let $\O_+$ and $\S$ be as above. Let $\Lambda: \mathit{L}^2(\S)^4\longrightarrow\mathit{L}^2(\S)^4$ be a bounded linear, self-adjoint operator. Define $ T= H +V$, with 
   \begin{align}
\mathrm{dom}(T)=\left\{ u+\Phi[g]: u\in\mathit{H}^1(\rr^3)^4, g\in\mathit{L}^2(\S)^4 \text{ and } T_{\S}u=\Lambda[g]\right\},
\end{align}
and $T(u+\Phi[g])=H(u)$, i.e $V(u+\Phi[g])=-g$,  for all $u+\Phi[g] \in \mathrm{dom}(T)$. If $\Lambda$ is Fredholm, then $(T,\mathrm{dom}(T))$ is self-adjoint.
  \end{theorem}
  We omit the proof of this theorem, since it is exactly the same as \cite[Theorem 2.11]{AMV1}.  Now, for our consideration if we let
  \begin{align}
(T,\mathrm{dom}(T))=(H_\k,\mathrm{dom}(H_\t)),\quad V=V_\t,\quad  \Lambda=-\Lambda_{\t,+}.
\end{align}
 Then,  for all $\varphi=u+\Phi[g]\in\mathrm{dom}(T)$, it holds that 
\begin{align*}
V(\varphi)=A_\t(t_{\S}u +\mathit{C}_{\S}[g])=A_\t(-\Lambda_{\t,+}+\mathit{C}_{\S})[g]=-(A_\t)( A_{\t}^{-1})[g]=-g.
 \end{align*}
 where in the first equality Lemma \ref{lemme 2.1} was used. Therefore, $V(\varphi)=-g$ and hence $ T(\varphi)=H(u)$ holds in the sense of distribution for all $\varphi=u+\Phi[g]\in\mathrm{dom}(T)$.  Hence, if $\Lambda_{\t,+}$  is Fredholm, we then fall under the conditions of Theorem \ref{luis}, which means that $(H_\t,\mathrm{dom}(H_\t))$ is self-adjoint. 
 
To  study  the spectral properties of $H_\t$ (for non-critical parameter) we shall restrict ourselves to the case of Lipschitz domains. The following proposition gives us a criterion for the existence of eigenvalues in the gap $(-m,m)$,  and a Krein-type resolvent
formula for $H_\t$. 
\begin{proposition}\label{surlespectre} Let $H_\t$ be as in \eqref{dom1} with a non-critical parameter  $\t\in\rr^n$. The following hold:
\begin{itemize}
 \item [(i)] Given $z\in(-m,m)$, then $\mathrm{Kr}(H_{\t}-a)\neq0$ $\Longleftrightarrow$ $\mathrm{Kr}( \Lambda^{z}_{\t,+})\neq0$. 
\item [(ii)] Given $z\in\cc\setminus\rr$ such that $\Lambda^{z}_{\t,+}$ is Fredholm. Then $\Lambda^{z}_{\t,+}$ is invertible in $\mathit{L}^2(\S)^4$ and it holds that 
 \begin{align}\label{Kresolvent}
    (H_{\t} -z)^{-1}(v)= (H -z)^{-1}(v) - \Phi^{z}(\Lambda^{z}_{\t,+} )^{-1} t_\S (H -z)^{-1}(v), \quad\forall\, v\in\mathit{L}^2(\rr^3)^4.
\end{align}
In particular, we have 
\begin{align}\label{spectreesse}
\mathrm{Sp}_{\mathrm{ess}}(H_{\t})=(-\infty,-m]\cup [m,+\infty). 
\end{align}
\end{itemize}
\end{proposition}
\textbf{Proof.}  The proof of item $(\mathrm{i})$ follows in the same way as in \cite[Proposition 4.1]{BB} (see also \cite[Proposition 3.1]{AMV2}).  Let us show $(\mathrm{ii})$. Fix $z\in\cc\setminus\rr$ such  that $\Lambda^{z}_{\t,+}$ is Fredholm. Then, from $(\mathrm{i})$ and the fact  $H_{\t}$ is self-adjoint it is clear that $\mathrm{Kr}(\Lambda^{z}_{+})=0 $ and $\mathrm{Rn}(\Lambda^{z}_{\t,+})= \mathit{L}^2(\S)^4$, as otherwise $z$ will be a non-real eigenvalue of $H_{\t}$. Hence,  we conclude that $\Lambda^{z}_{\t,+}: \mathit{L}^2(\S)^4\longrightarrow\mathit{L}^2(\S)^4$ is bijective and thus \eqref{Kresolvent} makes sense. Now given $v\in\mathit{L}^2(\rr^3)^4$,  we set 
\begin{align*}
 \varphi= (H -z)^{-1}(v) - \Phi^{z}(\Lambda^{z}_{\t,+} )^{-1} t_\S (H -z)^{-1}(v).
\end{align*}
To prove item $(\mathrm{ii})$, it remains to show that $\varphi\in\mathrm{dom}(H_\t)$. For this, remark  that $ \varphi= u+\Phi[g]$ where 
\begin{align*}
u=& (H -z)^{-1}(v) - (\Phi^{z}-\Phi)(\Lambda^{z}_{\t,+} )^{-1} t_\S (H -z)^{-1}(v),\\
  g=&-(\Lambda^{z}_{\t,+} )^{-1} t_\S (H -z)^{-1}(v).
\end{align*} 
Note that  $(\Lambda^{z}_{\t,+} )^{-1} t_\S (H -z)^{-1} $ is a bounded, and compact operator from $\mathit{L}^{2}(\rr^3)^4$ to $\mathit{L}^{2}(\S)^4$ and $(H -z)u= v+z \Phi[g]\in\mathit{L}^{2}(\rr^3)^4$. Consequently, we get that $g\in\mathit{L}^{2}(\S)^4 $ and $u\in\mathit{H}^{1}(\rr^3)^4$. Moreover, using Lemma \ref{lemme 2.1}$(\mathrm{ii})$, we obtain 
\begin{align*}
t_{\S}u&= \left( t_\S (H -z)^{-1} - ( \mathit{C}^z_{\S} - \mathit{C}_{\S})(\Lambda^{z}_{\t,+} )^{-1} t_\S (H -z)^{-1}\right)(v)\\
&= \left( t_\S (H -z)^{-1}- ( \Lambda^{z}_{\t,+}  -\Lambda_{\t,+} )(\Lambda^{z}_{\t,+} )^{-1} t_\S (H -z)^{-1}\right)[v]=-\Lambda_{\t,+}[g].
\end{align*}
Thus $\varphi\in\mathrm{dom}(H_\t)$, which yields $(\mathrm{ii})$. Since the Sobolev embedding $\mathit{H}^{1/2}(\S)^4\hookrightarrow\mathit{L}^{2}(\S)^4$ is compact, it follows that $\Phi^{z}(\Lambda^{z}_{\t,+} )^{-1} t_\S (H -z)^{-1}\in\mathcal{K}(\mathit{L}^{2}(\rr^3)^4)$. Therefore, we deduce by Weyl’s theorem  that $\mathrm{Sp}_{\mathrm{ess}}(H_{\t})$ is given by \eqref{spectreesse}. This completes the proof of the proposition.\qed
 \subsection{Critical parameters} In order to avoid repetition, in this part we explain our strategy  to prove the self-adjointness in the critical case, which is a generalisation of the technique developed in \cite{BB}. Here we assume that $\O_+$ is a bounded $\mathit{C}^2$-smooth domain.  First, we record some  known results that will be important in the proof of the self-adjointness of $H_\t$ for critical parameters.
  \begin{proposition}\label{extension}(\cite[Proposition 3.1]{BB}) Let $\Phi^z$ and $\mathit{C}^z_{\S}$ be as in Lemma \ref{lemme 2.1}. Then, the following hold true:
\begin{itemize}
 \item[(i)] The trace operator $t_{\S}$ (which until now was defined on  $\mathit{H}^{1}(\Omega_\pm)^4$) has a unique extension to a bounded linear operator from $\mathit{L}^{2}(\Omega_\pm)^4$ to  $\mathit{H}^{-1/2}(\S)^4$.
  \item[(ii)] The operator $\Phi^{z}$ admits a continuous extension from $\mathit{H}^{-1/2}(\S)^4$ to $\mathit{L}^{2}(\rr^3)^4$, which we still denote $\Phi^{z}$.
  \item[(iii)] The operator $\mathit{C}^{z}_{\S}$ admits a continuous extension  $\tilde{\mathit{C}^{z}_{\S}}:\mathit{H}^{-1/2}(\S)^4\rightarrow\mathit{H}^{-1/2}(\S)^4$. Moreover, we have
  \begin{align}\label{dual}
  t_\S\Phi^{z}_{\O_\pm}[h]:=\tilde{\mathit{C}^{z}_{\pm}}[h]&= (\mp\frac{i}{2}(\alpha\cdot\mathit{N}) +\tilde{\mathit{C}^{z}_{\S}})[h], \\
       \langle \tilde{\mathit{C}^{z}_{\S}}[h],g\rangle_{\mathit{H}^{-1/2}(\S)^4,\mathit{H}^{1/2}(\S)^4}&=\langle h,\mathit{C}^{\overline{z}}_{\S}[g]\rangle_{\mathit{H}^{-1/2}(\S)^4,\mathit{H}^{1/2}(\S)^4},\label{dual}
    \end{align}
    for any $g\in\mathit{H}^{1/2}(\S)^4$ and  $ h\in\mathit{H}^{-1/2}(\S)^4$, where $\Phi^{z}_{\O_\pm}:\mathit{H}^{-1/2}(\S)^4  \longrightarrow \mathit{L}^2(\O_\pm)^4 $ is defined by  $\Phi^z_{\O_\pm}[h](x)=\Phi^z[h](x)$, for $h\in\mathit{H}^{-1/2}(\S)^4$  and $x\in\O_\pm$.
\end{itemize}
\end{proposition}
 
  All the problems that we are going to consider share the following properties when the parameter $\t$ is critical:
 \begin{itemize}
  \item[(P1)] $A_\t$ and $\tilde{A}_\t$ admit continuous extensions from $\mathit{H}^{-1/2}(\S)^4$ into itself, which we still denote by $A_\t$ and $\tilde{A}_\t$. Moreover, $ \mathrm{sgn}(\t)\mathit{I}_4= \tilde{A}_\t  A_\t= A_\t\tilde{A}_\t$, holds in $\mathit{H}^{-1/2}(\S)^4$.
    \item[(P2)] The operator  $\tilde{\Lambda}^{z}_{\t,+}A_\t\tilde{\Lambda}^{z}_{\t,-}$ is bounded from $\mathit{H}^{-1/2}(\S)^4$ to $\mathit{H}^{1/2}(\S)^4$, where $\tilde{\Lambda}^{z}_{\t,\pm}$ is the continuous extension of $\Lambda^{z}_{\t,\pm}$ defined from $\mathit{H}^{-1/2}(\S)^4$ onto itself.
\end{itemize}
 Now, following the same arguments as in \cite[Section 3]{BB}, one can easily show that $H_\t$ is closable. Moreover, the adjoint operator $H^{\ast}_\t$ acts in the sense of distribution as $H^{\ast}_\t(u+\Phi[g])=H(u)$, on the domains
   \begin{align}\label{adjoint}
  \mathrm{dom}(H^{\ast}_{\t})=\left\{ \varphi=u+\Phi[g]: u\in\mathit{H}^1(\rr^3)^4, g\in\mathit{H}^{-1/2}(\S)^4, t_{\S}u=-\tilde{\Lambda}_{\t,+}[g]\right\}.
  \end{align}
  Then we have the following theorem.
  \begin{theorem}\label{thincritic} Let $H_\t$ be as in \eqref{dom1} with a critical parameter  $\t\in\rr^n$. If $(\mathrm{P1})$ and $(\mathrm{P2})$ hold, then $H_\t$ is essentially self-adjoint and we have 
   \begin{align}\label{the closure}
  \mathrm{dom}(\overline{H_{\t}})=\left\{ \varphi=u+\Phi[g]: u\in\mathit{H}^1(\rr^3)^4, g\in\mathit{H}^{-1/2}(\S)^4, t_{\S}u=-\tilde{\Lambda}_{\t,+}[g]\right\}.
  \end{align}
  Moreover, $\mathrm{dom}(\overline{H_{\t}})\not\subset  \mathrm{dom}(H_{\t})$.
  \end{theorem}
  \textbf{Proof.} Since $H_\t$ is closable,  it is sufficient to show the inclusion $H^{\ast}_{\t}\subset\overline{H_{\t}}$. To this end, fix  $\varphi=u+\Phi[g]\in  \mathrm{dom}(H^{\ast}_{\t}) $ and let $(h_j)_{j\in\mathbb{N}}\subset \mathit{H}^{1/2}(\S)^4$ be a sequence of functions that converges to $g$ in $\mathit{H}^{-1/2}(\S)^4$. Set 
\begin{align}\label{construcsubseq}
   g_j:= g +\frac{1}{2} A_{\t} \tilde{\Lambda}_{\t,-}[h_j -g],\quad\forall j\in\mathbb{N}.
   \end{align} 
   Then
   \begin{align}\label{argCV}
   \begin{split}
    (g_j)_{j\in\mathbb{N}}&\subset \mathit{H}^{1/2}(\S)^4 \, \text{ and }\,  g_j \xrightarrow[j\to\infty]{} g   \text{ in }  \mathit{H}^{-1/2}(\S)^4,\\
    (\Lambda_{\t,+}[g_j])_{j\in\mathbb{N}}&\subset \mathit{H}^{1/2}(\S)^4  \, \text{ and }\, \Lambda_{\t,+}[g_j] \xrightarrow[j\to\infty]{} \tilde{\Lambda}_{\t,+}[g] ,  \text{ in } \mathit{H}^{1/2}(\S)^4.
    \end{split}
    \end{align}
   Indeed, observe that 
    \begin{align}
   -\frac{1}{2}A_\t\tilde{\Lambda}_{\t,-}[g]&= -\frac{1}{2}A_\t\left(\frac{1}{\mathrm{sgn}(\t)}\tilde{A}_\t- \tilde{\mathit{C}_{\S}} \right)[g]\\
   &=-g +\frac{1}{2}A_\t\tilde{\Lambda}_{\t,+}[g].
   \end{align}
   where the property $(\mathrm{P1})$ was used in the last equality. Hence, we obtain that  
   \begin{align}\label{}
   g_j:= \frac{1}{2} A_{\t} \left( \tilde{\Lambda}_{\t,+}[g]+ \Lambda_{\t,-}[h_j ] \right), \\
 \tilde{\Lambda}_{\t,+}[  g_j- g]= \frac{1}{2}\tilde{\Lambda}_{\t,+} A_{\t} \tilde{\Lambda}_{\t,-}[h_j -g].
   \end{align}
   Therefore, \eqref{argCV} follows by \eqref{construcsubseq}, the continuity of $\Lambda_{\t,-}$  in  $\mathit{H}^{1/2}(\S)^4$, and the property $(\mathrm{P2})$.
Now, define 
\begin{align}
v_j= \frac{1}{4}\left( E_{\O_+}( \tilde{\Lambda}_{\t,+} A_{\t}\tilde{\Lambda}_{\t,-}[h_j- g]) +E_{\O_-}( \tilde{\Lambda}_{\t,+} A_{\t}\tilde{\Lambda}_{\t,-}[h_j- g]) \right), \text{ for all } j\in\mathbb{N}.
\end{align}
Clearly,  $v_j\in\mathit{H}^1(\rr^3)^4 $ and  $v_j  \xrightarrow[j\to\infty]{} 0$ in $\mathit{H}^1(\rr^3)^4 $. Now  set  $\varphi_j:= u_j +\Phi[g_j]$, where $u_j=u-v_j$,  for all  $j\in\mathbb{N}$. Note that  
\begin{align*}
t_\S u_j= t_\S u- \frac{1}{2}\tilde{\Lambda}_{\t,+} A_{\t} \tilde{\Lambda}_{\t,-}[h_j -g]=-\tilde{\Lambda}_{\t,+}[g_j]+ ( t_\S u +\tilde{\Lambda}_{\t,+} [g])=\Lambda_{\t,+}[g_j].
\end{align*}
Thus, $t_\S u_j= -\Lambda_{\t,+}[g_j]$ holds in $ \mathit{H}^{1/2}(\S)^4$, and hence $(\varphi_j)_{j\in\mathbb{N}}\subset \mathrm{dom}(H_{\t})$.  Using that $H_\t(\varphi_j)=H(u)-H(v_j)$,  and the fact that $(g_j)_{j\in\mathbb{N}}$ converges to $g$ in $\mathit{H}^{-1/2}(\S)^4$ when $j\longrightarrow \infty$,  we obtain that  
\begin{align}
(\varphi_j,H_{\t}\varphi_j) \xrightarrow[j\to\infty]{}(\varphi, H^{\ast}_{\t}\varphi) \text{ in } \mathit{L}^2(\rr^3)^4\times \mathit{L}^2(\rr^3)^4 .
\end{align}
 Therefore $ H^{\ast}_{\t}\subset\overline{H_{\t}}$.  Now the fact that $\mathrm{dom}(\overline{H_{\t}})\not\subset  \mathrm{dom}(H_{\t})$ follows essentially from the above arguments. Indeed, given any function $0\neq g\in\mathit{H}^{-1/2}(\S)^4\setminus\mathit{L}^{2}(\S)^4$,  set 
 \begin{align}
\varphi= \frac{1}{2}\left( E_{\O_+}( \tilde{\Lambda}_{\t,+} A_{\t}\tilde{\Lambda}_{\t,-}[g]) +E_{\O_-}( \tilde{\Lambda}_{\t,+} A_{\t}\tilde{\Lambda}_{\t,-}[g]) \right)-\Phi[A_{\t}\tilde{\Lambda}_{\t,-}[g]].
\end{align}
Then, it is clear that $\varphi\notin \mathrm{dom}(H_{\t})$ and $\varphi\in \mathrm{dom}(\overline{H_{\t}})$. This finishes the proof of the theorem.\qed
\newline

The next proposition gives us another way to define the Dirac operator $\overline{H_{\t}}$.
  \begin{proposition}\label{defcritic} Let $\overline{H_{\t}}$ be as in Theorem \ref{thincritic}. Then we have 
   \begin{align*}
  \mathrm{dom}(\overline{H_{\t}})=\bigg\{ (\varphi_+,\varphi_-)&\in\mathit{L}^2(\Omega_+)^4\oplus\mathit{L}^2(\Omega_-)^4:  (\alpha\cdot\nabla) \varphi_\pm\in\mathit{L}^2(\Omega_{\pm})^4 \text{ and }\\
  & \left(\frac{1}{2} + iA^{-1}_\t(\aN)\right) t_\S\varphi_{+}=- \left(\frac{1}{2} + iA^{-1}_\t(\aN)\right) t_\S\varphi_{-}\bigg\}, 
  \end{align*}
  where the transmission condition holds in $\mathit{H}^{-1/2}(\S)^4$.
  \end{proposition}
 \textbf{Proof.} Given  $\varphi=(u+\Phi[g])\in\mathrm{dom}(H_{\t})$, set $\varphi_{\pm} := \varphi|_{\Omega_\pm}$. Then, 
 a simple computation in the sens of distributions yields
\begin{align*}
 (\overline{H_\t}+A_\t\delta_{\S})\varphi=& (-i\alpha\cdot\nabla +m\beta)\varphi +\frac{1}{2}A_\t(t_{\S}\varphi_+ +t_{\S}\varphi_-)\delta_{\S},\\
=& (-i\alpha\cdot\nabla +m\beta)\varphi_+\oplus(-i\alpha\cdot\nabla +m\beta)\varphi_- +i\alpha\cdot\mathit{N}(t_{\S}\varphi_+ -t_{\S}\varphi_-)\delta_{\S}\\
&+\frac{1}{2}A_\t(t_{\S}\varphi_+ +t_{\S}\varphi_-)\delta_{\S} .
\end{align*}
 Thus, the proposition follows from this, the definition of $\tilde{\Lambda}_{\t,+}$ and Proposition \ref{extension}.\qed

    \subsection{On the confinement}\label{Subb3.3} In this part, we briefly discuss the case where the operator $H_\t$ generates confinement. In such a case we shall always restrict ourselves to Lipschitz domains.  We mention that, for some Dirac operators (for example the coupling of the free Dirac operator with the Lorentz scalar $\delta$-potential), one can prove that they generate the confinement for UR domains, provided that $\mathit{B}^2_{1/2}(\S)$ is the trace to $\mathit{H}^1(\O_\pm)^4$. For instance  this is possible if $\O_+$ is a two-sided NTA domain with an ADR boundary by combing Proposition \ref{Les traces} with Jones's results \cite[Theorem 1 and Theorem 2]{PJ},  for more details we refer to \cite{PH}.  As it was observed in \cite{BB} (see also \cite{CLMT} for the two-dimensional case) the anomalous magnetic $\delta$-potential generates confinement with a critical parameter. So, here we are going to show how to deal with both situations, i.e confinement with critical or non-critical parameters. 
    
Recall the definition of $\mathrm{dom}(H_{\t})$ from \eqref{dom1}.  Let $\Phi_{\O_\pm}:\mathit{L}^2(\S)^4  \longrightarrow \mathit{L}^2(\O_\pm)^4 $  be the operators defined by  $\Phi_{\O_\pm}[g](x)=\Phi[g](x)$, for $g\in\mathit{L}^2(\S)^4$  and $x\in\O_\pm$.  Given any $\varphi=(u+\Phi[g])\in\mathrm{dom}(H_{\t})$, we set
\begin{align}\label{decompfunc}
\varphi_{\pm} := \varphi|_{\Omega_\pm}=u|_{\Omega_\pm}+\Phi_{\O_\pm}[g].
\end{align} 
For simplicity, we denote by $\lim\limits_{ \mathrm{nt}}\varphi_{\pm}$ the nontangential limit of $\varphi_{\pm}$. By definition it holds that
\begin{align}\label{calculconf}
\begin{split}
 t_{\S}u=-\Lambda_{\t,+}[g]  &\Longleftrightarrow  t_{\S}u +\mathit{C}_\S[g] =-A^{-1}_\t [g]\\
 &\Longleftrightarrow  \frac{1}{2}(  \lim\limits_{ \mathrm{nt}}\varphi_{+}+  \lim\limits_{ \mathrm{nt}}\varphi_{-}) = -iA^{-1}_\t(\alpha\cdot \mathit{N}))(  \lim\limits_{ \mathrm{nt}}\varphi_{+}- \lim\limits_{ \mathrm{nt}}\varphi_{-})\\
 &\Longleftrightarrow  \left(\frac{1}{2} +iA^{-1}_\t(\aN)\right) \lim\limits_{ \mathrm{nt}}\varphi_{+} = - \left(\frac{1}{2} -iA^{-1}_\t(\aN)\right)  \lim\limits_{ \mathrm{nt}}
 \varphi_{-}.
 \end{split}
\end{align} 
From this, we get the following properties:
\begin{itemize}
  \item[(P3)]  $\left(1/2 \pm iA^{-1}_\t(\aN)\right)$ are projectors in $\mathit{L}^2(\S)^4$.
   \item[(P4)]  $\mathrm{sgn}(\t)=-4$ and $\tilde{A}_\t(\aN)=(\aN)A_\t$. 
\end{itemize} 
Then, the following proposition illustrates the phenomenon of confinement for non-critical parameters.
\begin{proposition}\label{confinementgenerale} Let $H_\t$ be as in \eqref{dom1} with a non-critical parameter  $\t\in\rr^n$.  If $(\mathrm{P3})$ or $(\mathrm{P4})$ holds, then  $H_{\t}$ generates confinement and we have 
$$H_{\t}\varphi=H_\t^{\Omega_+}\varphi_+\oplus H_\t^{\Omega_-}\varphi_-=\left(-i\alpha\cdot \nabla+m\beta\right)\varphi_+\oplus\left(-i\alpha\cdot \nabla+m\beta\right)\varphi_-,$$
    where $H_\t^{\Omega_\pm}$ are the self-adjoint Dirac operators defined on
\begin{align*}
\dom(H_\t^{\Omega_\pm})=\bigg\{ u_{\Omega_\pm}+\Phi_{\O_\pm}[g]&: u_{\Omega_\pm}\in\mathit{H}^1(\Omega_\pm)^4, g\in\mathit{L}^2(\S)^4 \text{ and }\\
& \left(\frac{1}{2} \pm iA^{-1}_\t(\aN)\right)(t_\S u_{\Omega_\pm}+\mathit{C}_\pm[g])=0 \bigg\}.
\end{align*}

\end{proposition}
\textbf{Proof.}   If $(\mathrm{P3})$ holds true, then the proof  follows directly from \eqref{calculconf}. Assume that $(\mathrm{P4})$ holds true, then a simple computation yields that 
\begin{align*}
\label{}
    \left(\frac{1}{2} \pm iA^{-1}_\t(\aN)\right)A_\t \left(\frac{1}{2} \pm iA^{-1}_\t(\aN)\right)&=-4\left(\frac{1}{2} \pm iA^{-1}_\t(\aN)\right)\\
      \left(\frac{1}{2} \pm iA^{-1}_\t(\aN)\right)A_\t \left(\frac{1}{2} \mp iA^{-1}_\t(\aN)\right)&=0,
\end{align*}
Again, using \eqref{calculconf} we get the desired result.\qed
\newline

Now, in the case of a critical parameter, one need to replace $(\mathrm{P3})$ and $(\mathrm{P4})$ by the following properties:
 \begin{itemize}
  \item[($\mathrm{P}^{\prime}3$)]  $\left(1/2 \pm iA^{-1}_\t(\aN)\right)$ are projectors in $\mathit{H}^{-1/2}(\S)^4$.
   \item[($\mathrm{P}^{\prime}4$)]  $\mathrm{sgn}(\t)=-4$ and $\tilde{A}_\t(\aN)=(\aN)A_\t$. 
\end{itemize} 
 Here $A_\t$ (respectively  $\tilde{A}_\t$) is the extension given by the property $(\mathrm{P1})$. Then, using Proposition \ref{defcritic} and  following essentially the same arguments as Proposition \ref{confinementgenerale}, we get the following result for the confinement in this case.
 \begin{proposition}\label{confinementgenerale1} Let $\overline{H_{\t}}$ be as in \eqref{dom1} with a critical parameter  $\t\in\rr^n$.  If $(\mathrm{P}^{\prime}3)$ or $(\mathrm{P}^{\prime}4)$  holds true, then  $\overline{H_{\t}}$ generates confinement and we have 
$$\overline{H_{\t}}\varphi=H_\t^{\Omega_+}\varphi_+\oplus H_\t^{\Omega_-}\varphi_-=\left(-i\alpha\cdot \nabla+m\beta\right)\varphi_+\oplus\left(-i\alpha\cdot \nabla+m\beta\right)\varphi_-,$$
    where $H_\t^{\Omega_\pm}$ are the self-adjoint Dirac operators defined on
\begin{align*}
\dom(H_\t^{\Omega_\pm})=\bigg\{ \varphi_\pm&\in\mathit{L}^2(\Omega_\pm)^4:  (\alpha\cdot\nabla) \varphi_\pm\in\mathit{L}^2(\Omega_{\pm})^4  \text{ and } \left(\frac{1}{2} \pm iA^{-1}_\t(\aN)\right)t_\S\varphi_\pm =0 \bigg\}.
\end{align*}
where the boundary conditions holds in $\mathit{H}^{-1/2}(\S)^4$.
\end{proposition}
\subsection{$\beta$ and $\g$ transformations of the electrostatic and the magnetic $\delta$-potentials}\label{subb3.4}
In this part,  we introduce the $\beta$ and the $\g$ transformations of the electrostatic and the magnetic $\delta$-potentials.
We insist that we are not formulating a theory here, but rather we are describing facts which are based on our observation.

 Let $\ep,\e\in\rr$, recall that the electrostatic and the magnetic $\delta$-shell interactions of strength $\ep$ and $\e$,  respectively, supported on $\S$ are defined by
 \begin{align}
 V_\ep:=\ep I_4\delta_{\S} \,\text{ and }\,V_{\e}= \eta(\alpha\cdot N) \delta_{\S}.
 \end{align} 
Then, the $\beta$ transformation, which we denote by $\G_\beta$ is the multiplication operator by $\beta$ which preserves the symmetry of the above $\delta$-potentiels, that is 
 \begin{align}
\G_\beta( V_\ep):=\ep\beta\delta_{\S} \,\text{ and }\,\G_\beta(V_{\e})= i\eta\beta(\alpha\cdot N) \delta_{\S}.
 \end{align} 
 Thus, the $\beta$ transformation of the electrostatic $\delta$-potential gives the Lorentz scalar $\delta$-potential, and the $\beta$ transformation of the magnetic $\delta$-potential gives what it was called by the author in \cite{BB} (and independently in \cite{CLMT} in the two-dimensional setting) the \textit{anomalous magnetic} $\delta$-potential.
 
Similarly, the $\g$ transformation denoted by $\G_\g$, is the multiplication operator by $\g_5$ which preserves the symmetry, thus we get  
 \begin{align}
\G_\g( V_\ep):=\ep\g_5 I_4\delta_{\S} \,\text{ and }\,\G_\g(V_{\e})= \eta\g_5 (\alpha\cdot N) \delta_{\S}.
 \end{align} 
Again, the $\g$ transformation of the electrostatic $\delta$-potential  was already considered by the author in \cite{BB}. In there, it was shown that  $\pm2\g_5 I_4\delta_{\S}$ coincided with the electrostatic $\delta$-potential of constant strength $\mp2$, cf. \cite[Remark 5.2]{BB}, so we called it here the \textit{modified electrostatic} $\delta$-potential. Also we call $ \eta\g_5 (\alpha\cdot N) \delta_{\S}$ the  \textit{modified magentic} $\delta$-potential. 

Finally, we have the composition of the  $\beta$ transformation and the $\g$ transformation which gives us the following potentials
  \begin{align}
\G_\beta\G_\g( V_\ep):=i\ep\g_5 \beta\delta_{\S}\,  \text{ and }\, \G_\beta\G_\g(V_{\e})= i\eta\g_5\beta (\alpha\cdot N) \delta_{\S},
 \end{align} 
and we call them respectively the \textit{modified Lorentz scalar} $\delta$-potential and the \textit{modified anomalous magnetic} $\delta$-potential. 

To put things in order we use the following notations:
\begin{align} \label{Poto1}
\begin{split}
V_{\ept}= \ept\g_5\delta_{\S},\quad V_{\mu}=\mu\beta\delta_{\S},\quad V_{\mt}=i\mt\g_5\beta\delta_{\S},  \quad  \ept,\m,\mt\in\rr.
\end{split}
\end{align}
\begin{align} \label{Poto2}
\begin{split}
V_{\et}=  \et\g_5(\aN)\delta_{\S},\quad V_{\u}=  i\u\beta(\aN)\delta_{\S} ,\quad V_{\ut}=i\ut\g_5\beta(\aN)\delta_{\S},  \quad \et,\u,\ut\in\rr.
\end{split}
\end{align}
The reader may wonder why we introduced such transformations, and what is the interest behind that.  In fact, the answer is:
\begin{itemize}
  \item The $\g$ transformation preserves the characteristics of the potentials given by \eqref{Poto1}  and  \eqref{Poto2}. That is, the characteristics of the potentials are stable under the $\g$ transformation in the sense that,  for any potential $V_{\bullet}$ from \eqref{Poto1} or \eqref{Poto2}, if the parameter is critical then it remains critical after applying the $\g$ transformation, and the same holds true to the case of the confinement.
  \item   The characteristics of the potentials  given by \eqref{Poto1}  and  \eqref{Poto2} are not stable under the $\beta$ transformation, in the sense that, if $V_{\bullet}$ has a critical parameter then $\G_\beta(V_{\bullet})$ does not  have a critical parameter and vice versa. Similarly,   if $V_{\bullet}$ does not generates confinement then $\G_\beta(V_{\bullet})$ generates confinement and vice versa.
\end{itemize}

As a simple example, one can consider the electrostatic $\delta$-potential. As is well-known,  (\cite{OBV,BH}) $V_\ep$ has a critical parameter which is $\ep=\pm2$, and it does not generates confinement. Now that $V_{\mu}$ generates confinement is known (\cite{AMV2}), moreover by \cite[Section 5]{BB} we also know that $V_{\ept}$ has a critical parameter which is $\ept=\pm2$. More generally, one can prove the above facts using directly  Definition \ref{defcritic}, $(\mathrm{P3})$, $(\mathrm{P4})$,  $(\mathrm{P}^{\prime}3)$ and $(\mathrm{P}^{\prime}4)$. Thus we conclude that 

\begin{itemize}
  \item  $V_{\mu}$ and  $V_{\mt}$ generate confinement (for $\m=\mt=\pm2$) without critical parameters. 
  \item $V_{\u}$ and  $V_{\ut}$ generate confinement for the critical parameters  $\u=\ut=\pm2$.
\end{itemize}
We are going to  show this in detail in the next section.

\section{Delta interactions involving  $\beta$ and $\g$ transformations of the electrostatic and the magnetic $\delta$-potentials}\label{sec5}
\setcounter{equation}{0}
As the title of this section indicates, here we focus on  the spectral  study of the Dirac operator $ H_\t$ when  $V_\t$ is a 
combination of the $\delta$-potentials given by \eqref{Poto1}  and  \eqref{Poto2}. We first consider the following Dirac operator
\begin{align}
H_{\k} =  H + V_{\k}=H+(\ep I_4 +\mu\beta + \eta(\alpha\cdot N))\delta_{\S},\quad  \k:=(\ep,\m,\e)\in\rr^3,
\end{align}
which was already study by the author in \cite{BB}, for a critical and non-critical parameter, when $\O_+$ is a $\mathit{C}^2$-smooth domain. Thus, we only focus  on the spectral properties of $H_\k$ for non-critical parameter in the case of UR domains. 

For the convenience of the reader, we begin our study with the subclass of bounded Lipschitz domains with $\mathrm{VMO}$ normals, where  we can discuss the spectral properties of $H_\k$, for all  $\mathrm{sgn}(\k)\neq0,4$. Subsequently, we discuss  the Sobolev regularity of $\dom(H_\k)$ in the case of bounded $\mathit{C}^{1,\o}$-smooth domains. Finally, we study separately the couplings $(H+ V_\ep+V_\m)$ and $(H+ V_\e)$ for general UR domains.

 Afterwards,  we consider in Subsection \ref{sub4.4} the following Dirac operators 
\begin{align}\label{autresoperator}
\begin{split}
H_{\mt} &=  H + V_{\ut}=H+i\mt\g_5\beta\delta_{\S},\quad \mt \in\rr,\\
H_{\ut} &=  H + V_{\ut}=H + i\ut\g_5\beta(\aN)\delta_{\S},  \quad \ut\in\rr,
\end{split}
\end{align}
which deserve to be analysed in detail. Since $V_{\ept}$ and  $V_{\et}$ can be treated in a similar way  as  $V_\ep$ and  $V_\e$  respectively, and the potential $V_\u$ have been  studied in \cite{BB}, so in order to avoid repetition we will only make a few remarks about the latter potentials.

\subsection{$\delta$-interactions supported on the boundary of a Lipschitz domain with a normal in $\mathrm{VMO}$ }\label{sub5.1}
 In what follows, unless otherwise specified, we always suppose that $\S$ satisfies the following property:  
\begin{itemize}
  \item[(H1)]  $\S=\partial\Omega_+$ with $\Omega_+$  a bounded Lipschitz  domain with a normal $\mathit{N}\in\mathrm{VMO}(\partial\O,\mathrm{dS})^3$. 
  \end{itemize}
  Roughly speaking, the above assumption implies  smallness of the Lipschitz constant of $\Omega_+$.  Another way to reformulate  the assumption $(\mathrm{H1})$ is to say that $\Omega_+$ belongs to the intersection of the class of bounded Lipschitz domains and the class of regular SKT domains, see the proof of Proposition \ref{Steve}.
 
Recall that $\mathrm{sgn}(\k)=\ep^2-\m^2-\e^2$,  and the operator $\lambda^{z}_{\k,+}$ is given by
\begin{align}
\Lambda^{z}_{\k,+}=\frac{1}{\mathrm{sgn}(\k)}(\ep I_4 -(\mu\beta+ \eta(\alpha\cdot N)) )+\mathit{C}^{z}_{\S},\quad\forall z\in\cc\setminus\left((-\infty,-m]\cup[m,\infty)\right).
\end{align} 
For technical reasons, we define the operator $\lambda^{z}_{\k,-}$ as follows 
\begin{align}
\Lambda^{z}_{\k,-}=\frac{1}{\mathrm{sgn}(\k)}(\ep I_4 +\mu\beta+ \eta(\alpha\cdot N) )-\mathit{C}^{z}_{\S},\quad\forall z\in\cc\setminus\left((-\infty,-m]\cup[m,\infty)\right).
\end{align} 
 
  Now we can state the main result about the spectral properties of the Dirac operator $ H_{\k}$.
\begin{theorem}\label{main4}  Let $\k\in\rr^3$ such that  $\mathrm{sgn}(\k)\neq0,4$, and assume that $\S$ satisfies $(\mathrm{H1}) $, and let $H_{\k}$ be as in \eqref{dom1}. Then $H_{\k}$ is self-adjoint. Moreover, the following hold true:
\begin{itemize}
 \item [(i)]For all $z\in\cc\setminus\rr$, it holds that 
 \begin{align}\label{resolvent2}
    (H_{\k} -z)^{-1}= (H -z)^{-1} - \Phi^{z}(\Lambda^{z}_{\k,+} )^{-1}  t_\S (H -z)^{-1}.
\end{align}
\item [(ii)]  $\mathrm{Sp}_{\mathrm{ess}}(H_{\k})=(-\infty,-m]\cup [m,+\infty)$. 
\item [(iii)] $a\in\mathrm{Sp}_{\mathrm{p}}(H_{\k})$ if and only if $-a\in\mathrm{Sp}_{\mathrm{p}}(H_{\tilde{\k}})$, where  $\tilde{\k}$ is given by
$$\tilde{\k}=\left(-\frac{4\ep}{\mathrm{sgn}(\k)},-\frac{4\m}{\mathrm{sgn}(\k)},\frac{4\e}{\mathrm{sgn}(k)}\right).$$ 
\end{itemize}
\end{theorem}

Before proving this result, we first give a characterization of the assumption $(\mathrm{H1}) $ via the compactness of the anticommutator $\lbrace \alpha\cdot \mathit{N}, \mathit{C}^{z}_{\S}\rbrace$. 

 Let  $g\in\mathit{L}^2(\S)$, then the harmonic double layer $K$ and the  Riesz transforms $(R_k)_{1\leqslant k \leqslant3}$ on $\S$ are defined by
 \begin{align}\label{pv double layer2}
 \begin{split}
K[g](x) &=\lim\limits_{\rho\searrow 0}\int_{|x-y|>\rho}\frac{\mathit{N}(y)\cdot(x-y)}{4\pi|x-y|^3}g(y)\mathrm{dS}(y),\\
 R_k[g](x) &=\lim\limits_{\rho\searrow 0}\int_{|x-y|>\rho}\frac{x_k-y_k}{4\pi|x-y|^3}g(y)\mathrm{dS}(y).
  \end{split}
  \end{align}
 Also,  for $z\in\cc\setminus\left((-\infty,-m]\cup[m,\infty)\right)$, the trace of the single-layer potential associated to $(\Delta +m^2-z^2)I_4$ is given by 
 \begin{align}\label{SL}
S^{z}[g](x)=\int_{\S} \psi^z(x-y)g(y)\mathrm{dS}(y), \quad \forall x\in\S \text{ and } g\in\mathit{L}^{2}(\S)^4,
\end{align}
where 
\begin{align}
\label{}
 \psi^z(x)=\frac{e^{i\sqrt{z^2-m^2}|x|}}{4\pi|x|}\mathit{I}_4,\quad\text{ for } x\in\rr^3.
\end{align}
 If $z=0$, we simply write $S:=S^0$. 
 
  The following proposition is implicitly contained in \cite{HTM}, but we state and prove it here for the sake of completeness.  
  \begin{proposition}\label{Steve} Assume that $\S$ satisfies $(\mathrm{H1}) $. Then, the harmonic double layer $K$ and the commutators $\left[ \mathit{N}_j, R_k\right]$, $1\leqslant j,k \leqslant3$, are compact operators on $\mathit{L}^2(\S)$.
  \end{proposition}
  Before giving the proof, we need to introduce the notion of  bounded regular Semmes-Kenig-Toro domains (regular SKT domains for short) developed by S. Hofmann, M. Mitrea and M. Taylor in \cite{HTM}. 
  \begin{definition}[regular SKT Domains] \label{SKTdomain}We say that a bounded open set $\O\subset\rr^3$ is a regular  Semmes-Kenig-Toro domain, or briefly regular SKT domain,   provided $\O$ is two-sided NTA domain, $\partial\O$ is ADR and   whose geometric measure theoretic outward unit normal $\nu\in\mathrm{VMO}(\partial\O,\mathrm{dS})^3$. 
\end{definition}
\begin{remark}We mention that Definition \ref{SKTdomain} is rather a characterization of regular SKT domains, for the precise definition we refer to \cite[Definition 4.8]{HTM}.
\end{remark}
   \textbf{Proof of Proposition \ref{Steve}.} The result follows from the fact that $\Omega_+$ is regular SKT domain. To see this indeed, note that bi-Lipschitz mappings preserve the class  of two-sided NTA domains with Ahlfors regular boundaries, and the class of regular SKT domains is invariant under continuously differentiable diffeomorphisms, see \cite{HTM2}. Now, by definition $\Omega_+$  is  locally the region above the graph of a Lipschitz function $\phi: \rr^2\longrightarrow \rr$. Therefore, one may (and do) assume (via a partition of unity and a local flattening of
the boundary) that 
   \begin{align*}
  \O_+ =\{x=(\overline{x},x_3)\in\rr^2\times\rr: x_3>\phi(\overline{x})\}.
   \end{align*}
  Let  $F: \rr^3\longrightarrow \rr^3$ be defined for all $(\overline{x},x_3)\in\rr^2\times\rr$ as $F(\overline{x},x_3):=(\overline{x},x_3+\phi(\overline{x} ))$. Then it easily follows that $F$ is a bijective function with inverse  $F^{-1}: \rr^3\longrightarrow \rr^3$ given by $F^{-1} (\overline{y},y_3):=(\overline{y},y_3-\phi(\overline{y} ))$ for all $(\overline{y},y_3)\in\rr^2\times\rr$. Moreover, $F$ and $F^{-1}$ are both Lipschitz functions with constants  $L_F, L_{F^{-1}}\leqslant(1+ \left|\left| \nabla \phi\right|\right|_{\mathit{L}^\infty})$.  It is clear that $\O_+$ (respectively $ \O_-$) is the image of  $\rr^3_+$ (respectively $\rr^3_-$) under the bi-Lipschitz homeomorphism $F$,  which also maps $\rr^2\times\{0\}$ onto $\S$. From this,  it follows that   $\Omega_+$ is a two-sided NTA domain and $\S$ is  Ahlfors regular (because  $\rr^3_+$ is a two sided NTA domain and $\rr^2\times\{0\}$ is Ahlfors regular). Since  $\mathit{N}\in \mathrm{VMO}(\S)$ by assumption, thanks to  \cite[Theorem 4.21]{HTM}, we know that  $\Omega_+$ is a regular SKT domain. Therefore the claimed result follows by \cite[Theorem 4.47]{HTM}.\qed
  \newline

  \begin{lemma}\label{L1} Assume that $\S$ satisfies $(\mathrm{H1}) $. Then,  $\lbrace \alpha\cdot \mathit{N}, \mathit{C}^{z}_{\S}\rbrace$ is a compact operator on $\mathit{L}^{2}(\S)^4$, for all $z\in\cc\setminus\left((-\infty,-m]\cup[m,\infty)\right)$.
  \end{lemma}
  \textbf{Proof. } Given  $g\in\mathit{L}^2(\S)^4$, then a  straightforward   computation shows that
  \begin{align}\label{eq5.5}
  \begin{split}
    \lbrace \alpha\cdot \mathit{N}, \mathit{C}^{z}_{\S}\rbrace[g](x)&= T_{K_1}[g](x)  + T_{K_2}[g](x),
    \end{split}
    \end{align}
where the kernels $K_j$, $j=1,2$,  are given by 
\begin{align*}
 K_1(x,y)=& \frac{e^{i\sqrt{z^2-m^2}|x-y|}}{4\pi |x-y|}\left(\alpha\cdot\mathit{N}(x)\right)\left(z+m\beta+\sqrt{z^2-m^2}\left( \alpha\cdot\frac{x-y}{|x-y|}\right)\right)\\ 
 &+   \frac{e^{i\sqrt{z^2-m^2}|x-y|}}{4\pi |x-y|}\left(z+m\beta+\sqrt{z^2-m^2}\left( \alpha\cdot\frac{x-y}{|x-y|}\right)\right)(\alpha\cdot\mathit{N}(y))\\
 &+ \frac{e^{i\sqrt{z^2-m^2}|x-y|} -1}{4\pi |x-y|^3}\left[(\alpha\cdot\mathit{N}(x))(i \alpha\cdot(x-y)) +   \left(i \alpha\cdot(x-y)\right)\left(\alpha\cdot\mathit{N}(y)\right)\right].
 \end{align*}
 \begin{align*}
 K_2(x,y)= \frac{i}{4\pi |x-y|^3}\left( (\mathit{N}(x))( \alpha\cdot(x-y)) + \alpha\cdot(x-y))(\mathit{N}(y)) \right).
 \end{align*}
Using the estimate
\begin{align}\label{majoration}
\left| e^{i\sqrt{z^2-m^2}|x|}-1\right|\leqslant \left|\sqrt{z^2-m^2}\right||x|,
\end{align}
it easily follows that 
\begin{align}\label{estimation2}
\sup_{1\leqslant k,j\leqslant 4}\left| K_{1}(x-y)\right| &=\mathcal{O}(|x-y|^{-1})\quad \text{ when }  |x-y|\longrightarrow 0.
\end{align}
Once \eqref{estimation2} has been established, working component by component and using \cite[Lemma 3.11]{Fo},  one can show that  $T_{K_1}$ is a compact operator in $\mathit{L}^{2}(\S)^4$. Now it is straightforward to check that 
 \begin{align}\label{lasomme}
 T_{K_2}[g](x)&= \tilde{K}[g](x) + \tilde{K}^{\ast}[g](x) +\sum_{j=1}^{3}\sum_{\substack{k={1}\\ k\neq j} }^{3}\alpha_j\alpha_k \left[ \mathit{N}_j, \mathcal{R}_k\right][g](x).
 \end{align}
 where   $\tilde{K}$ denotes the matrix valued harmonic double layer,  $\tilde{K}^{\ast}$ is the associated adjoint operator, and $\mathcal{R}_k$ are the matrix versions of the Riesz transforms. That is, for $x\in\S$ and $g\in\mathit{L}^2(\S)^4$, we have
 \begin{align}\label{pv double layer}
 \begin{split}
  \tilde{K}[g](x) &=\lim\limits_{\rho\searrow 0}\int_{|x-y|>\rho}\frac{\mathit{N}(y)\cdot(x-y)}{4\pi|x-y|^3}\mathit{I}_4g(y)\mathrm{dS}(y),\\
   \tilde{K}^{\ast}[g](x) &=\lim\limits_{\rho\searrow 0}\int_{|x-y|>\rho}\frac{\mathit{N}(x)\cdot(x-y)}{4\pi|x-y|^3}\mathit{I}_4g(y)\mathrm{dS}(y),\\
   \mathcal{R}_k[g](x) &=\lim\limits_{\rho\searrow 0}\int_{|x-y|>\rho}\frac{x_k-y_k}{4\pi|x-y|^3}\mathit{I}_4g(y)\mathrm{dS}(y).
  \end{split}
  \end{align}
Since the adjoint of a compact operator is a compact operator and $\alpha_j$'s are constants matrices, using Proposition \ref{Steve}  and working component by component, we get that  $T_{K_2}$ is a  compact operators in $\mathit{L}^{2}(\S)^4$. Therefore $\lbrace \alpha\cdot \mathit{N}, \mathit{C}^{z}_{\S}\rbrace$ is a compact operator in $\mathit{L}^{2}(\S)^4$ and this finishes the proof of the lemma.\qed
\newline

 Note that Lemma \ref{L1} is not valid for general Lipschitz surfaces. In fact, it turns out that assuming $(\mathrm{H1}) $ means that we are excluding the special class of corner domains. Indeed, from Proposition \ref{Steve} we know that any bounded Lipschitz domain $\O_+$ is an NTA domain and $\S$ is ADR. However,  the presence of any angle $\theta\neq0$, implies that
 $$\mathrm{dist}( \mathit{N},\mathrm{VMO}(\S,\mathrm{dS})^3)>0,$$ 
 where the distance is taken in $\mathrm{BMO}(\S,\mathrm{dS})^3$, cf. \cite[Proposition  4.38]{HTM} and the discussion that precedes it. Hence, $\O_+$ is not a regular SKT domain and then by \cite[Theorem 4.47]{HTM},  the principale value of the harmonic double layer $K$ and the commutators $\left[ \mathit{N}_j, R_k\right]$, $1\leqslant j,k \leqslant3$, are not compact  on $\mathit{L}^2(\S)$. So, $  \lbrace \alpha\cdot \mathit{N}, \mathit{C}_{\S}\rbrace$ is not a compact operator on $\mathit{L}^2(\S)^4$, and thus the assumption $(\mathrm{H1}) $ is sharp. To make this clearer, we have the following result.
\begin{theorem}\label{qestion3} Let $\O_+$  be a bounded Lipschitz domain, such that the decomposition $\rr^3=\O_+\cup\S\cup\O_-$ holds, where $\partial\O_+=\S$. Then,  $\S$ satisfies $(\mathrm{H1}) $ if and only if $\lbrace \alpha\cdot \mathit{N}, \mathit{C}_{\S}\rbrace$ is  compact in $\mathit{L}^{2}(\S)^4$.
\end{theorem}
\textbf{Proof.} The first implication follows from Lemma \ref{L1} by taking $z=0$. Let us prove the reverse implication, so assume that $\lbrace \alpha\cdot \mathit{N}, \mathit{C}_{\S}\rbrace$ is compact in $\mathit{L}^{2}(\S)^4$. Recall the definition of the operator  $W$ from \eqref{opeWW}. Then,  from \eqref{eq5.5} it holds that 
\begin{align}
\lbrace \alpha\cdot \mathit{N}, \mathit{C}_{\S}\rbrace  =  T_{K_1}+ \begin{pmatrix}
\lbrace \sigma\cdot \mathit{N}, W \rbrace   & 0\\
 0& \lbrace \sigma\cdot \mathit{N},  W \rbrace
\end{pmatrix},
\end{align}
where $T_{K_1}$ is a compact operator in  $\mathit{L}^2(\S)^4$. Using this, it follows that 
\begin{align}\label{WetC}
\lbrace \alpha\cdot \mathit{N}, \mathit{C}_{\S}\rbrace \text{ is compact in } \mathit{L}^2(\S)^4 \Longleftrightarrow
\lbrace \sigma\cdot \mathit{N}, W \rbrace \text{ is compact in } \mathit{L}^2(\S)^2.
  \end{align}
  Hence, it remains to show that 
  \begin{align}
\lbrace \sigma\cdot \mathit{N}, W \rbrace \text{ is compact in } \mathit{L}^2(\S)^2  \Longrightarrow \S \text{ satisfies }(\mathrm{H3}). 
  \end{align}
  For this, note that from Proposition \ref{Steve} we know that $\O_+$ is a two-sided NTA domain and $\S$ is  ADR. So $\O_+$ satisfies the two-sided corkscrew condition with an  ADR boundary.  Hence, $\O_+$ is a uniformly rectifiable by  \cite[Corollary 3.9]{HTM}. Next, we claim that there exists $C>0$, depending only on the dimension, the uniform rectifiability and the ADR constants of $\S$, such that 
  \begin{align}\label{Wope}
  \mathrm{dist}\left( \mathit{N},\mathrm{VMO}(\partial\O,\mathrm{dS})^3\right) \leqslant C \mathrm{dist}\left(\lbrace  \sigma\cdot \mathit{N}, W \rbrace , \mathcal{K}( \mathit{L}^2(\S)^2)\right), 
  \end{align}
 where the distance in the right-hand side is measured in $\mathcal{B}(\mathit{L}^2(\S)^2)$.  Now, assume for instance that \eqref{Wope} holds true.  Since $\lbrace \sigma\cdot \mathit{N}, W \rbrace $ is compact in $ \mathit{L}^2(\S)^2$,  by \eqref{Wope}, it holds that  $\mathit{N}\in\mathrm{VMO}(\partial\O,\mathrm{dS})^3$. Therefore, $\S$ satisfies  $(\mathrm{H1}) $, which proves the theorem. Now, let us come back to  the proof of  \eqref{Wope}. Given $x,y\in\rr^3$, we define the following multiplication operator
 \begin{align}\label{Diracmulti}
x\odot y:= (\sigma\cdot x)(-\sigma\cdot y).
\end{align}
 Using the anticommutation properties of the Pauli matrices, it is easy to check that:
 \begin{align*}
x\odot x:= -|x|^2,\quad x\odot y+ y\odot x =-2 (x\cdot y)\mathit{I}_2,\quad \forall x,y\in\rr^3.
\end{align*}
 Now, we make the observation that the multiplication operator defined by \eqref{Diracmulti} has the same properties as the multiplication operator in the Clifford algebra $\mathcal{C}l_3$ (see \cite[Section 4.6]{HTM} for the precise definition). Moreover, $W(\sigma\cdot\mathit{N})$ plays the same role as the Cauchy-Clifford operator defined on $\mathit{L}^2(\S)\otimes\mathcal{C}l_3$\footnote{Note that  the algebra generated by the Pauli spin matrices $\sigma=(\sigma_1,\sigma_2,\sigma_3)$ (as an algebra on the real field) is isomorphic to $\mathcal{C}l_3$. Thus, $W(\sigma\cdot\mathit{N})$ can be viewed as  the restriction of the Cauchy-Clifford operator on $\mathit{L}^2(\S)\otimes\cc^2$.} (i.e. it acts on $\mathcal{C}l_3$-valued functions), cf. \cite[Section 4.6]{HTM}. Thus, one can adapt the same arguments of \cite[Theorem 4.46]{HTM} and show that the claim \eqref{Wope} holds true, we leave the details for the reader. This completes the proof of the theorem.\qed
 \newline

As it was done in \cite[Theorem 4.47]{HTM},  one can also characterize the class of bounded regular SKT domains via  the compactness of the anticommutators  $\lbrace \sigma\cdot \mathit{N}, W \rbrace $ in $\mathit{L}^2(\S)^2$, or equivalently via the compactness of the anticommutator $\lbrace \alpha\cdot \mathit{N}, \mathit{C}_{\S} \rbrace $ in $\mathit{L}^2(\S)^4$.  This is the purpose of the following proposition.
\begin{proposition}\label{remarkSKT} Let $\O\subset\rr^3$ be a bounded two-sided NTA domain with a compact,  ADR boundary. Then the following statements are equivalent:
\begin{itemize}
 \item [(i)]  $\O$ is a regular SKT domain.
 \item [(ii)] The harmonic double layer $K$ and the commutators $\left[ \mathit{N}_j, R_k\right]$, $1\leqslant j,k \leqslant3$, are compact operators on $\mathit{L}^2(\partial\O)$.
\item [(iii)]  $\lbrace \alpha\cdot \mathit{N}, \mathit{C}_{\partial\O} \rbrace $ is a compact operator on $\mathit{L}^2(\partial\O)^4$.
\item [(iv)] $\lbrace \sigma\cdot \mathit{N}, W \rbrace $ is a compact operator on $\mathit{L}^2(\partial\O)^2$.
\end{itemize}
\end{proposition}
\textbf{Proof.} $ (\mathrm{i})\Rightarrow (\mathrm{ii})$ is a consequence of  \cite[Theorem 4.21]{HTM} and \cite[Theorem 4.47]{HTM}. $ (\mathrm{ii})\Rightarrow (\mathrm{iii})$  readily follows from \eqref{eq5.5} and  \eqref{lasomme}. $ (\mathrm{iii})\Rightarrow (\mathrm{iv})$ is an immediate consequence of \eqref{WetC}. Finally, $ (\mathrm{iv})\Rightarrow (\mathrm{i})$ follows from \eqref{Wope} and  \cite[Theorem 4.21]{HTM}. \qed
\newline

\begin{corollary}\label{Fredy} Let $z\in\cc\setminus\left((-\infty,-m]\cup[m,\infty)\right)$, then $\Lambda_{\k,\pm}^{z}$ is a Fredholm operator on $\mathit{L}^{2}(\S)^4$.
\end{corollary}
\textbf{Proof.} Fix $z\in\cc\setminus\left((-\infty,-m]\cup[m,\infty)\right)$.  Observe that 
 \begin{align}\label{psi}
 \{ \beta,\mathit{C}^{a}_{\S}\}[g](x)=2(m \mathit{I}_4+z\beta) S^{a}[g](x).
\end{align}
Using this, it follows that 
\begin{align}
\begin{split}
\Lambda^z_{\k,\mp}\Lambda^z_{\k,\pm}=&\frac{1}{\mathrm{sgn}(\k)}-\frac{1}{4}-\mathit{C}^{z}_{\S} (\alpha\cdot \mathit{N})\lbrace \alpha\cdot \mathit{N}, \mathit{C}^z_{\S}\rbrace +\frac{2\mu}{\mathrm{sgn}(\k)}( m\mathit{I}_4 +z\beta)S^z \\
&+\frac{\e}{\mathrm{sgn}(\k)} \lbrace \alpha\cdot \mathit{N}, \mathit{C}^z_{\S}\rbrace\label{multi111},
\end{split}
\end{align}
 As  $S^z$  is bounded from $\mathit{L}^{2}(\S)^4$ to $\mathit{H}^{1/2}(\S)^4$ (see \cite[Theorem 6.11]{Mc} for example), and  the injection $\mathit{H}^{1/2}(\S)^4$ into $\mathit{L}^{2}(\S)^4$ is compact,  it follows that $S^z$  is a compact operator in $\mathit{L}^{2}(\S)^4$.  Now, using that  $\mathit{C}^{z}_{\S} (\alpha\cdot \mathit{N})$ is bounded in $\mathit{L}^{2}(\S)^4$ and that $\lbrace \alpha\cdot \mathit{N}, \mathit{C}^z_{\S}\rbrace$ is a compact operator on $\mathit{L}^{2}(\S)^4$ by Lemma \ref{L1},  we thus obtain that $\mathit{C}^{z}_{\S} (\alpha\cdot \mathit{N})\lbrace \alpha\cdot \mathit{N}, \mathit{C}^z_{\S}\rbrace$ is a compact operator on $\mathit{L}^{2}(\S)^4$.  Hence  $\Lambda^z_{\k,\mp}\Lambda^z_{\k,\pm}$ is Fredholm operator and therefore $\Lambda^z_{\k,\pm}$  is Fredholm operator by \cite[Theorem 1.46 ($\mathrm{iii}$)]{Aiena}. This finishes  the proof of the corollary.  \qed
 \newline

 Now we are in position to complete the proof of Theorem \ref{main4}.
 
\textbf{Proof of Theorem \ref{main4} } Since $\Lambda^z_{\k,+}$ is Fredholm for all $z\in\cc\setminus\left((-\infty,-m]\cup[m,\infty)\right)$, a direct application of Theorem \ref{luis} yields that $H_\k$ is self-adjoint which proves the first statement of the theorem. Assertions $(\mathrm{i})$ and $(\mathrm{ii})$ are consequences of  Proposition \ref{surlespectre}. Thus it remains to show $(\mathrm{iii})$. For this, we let 
$$\tilde{\k}=\left(-\frac{4\ep}{\mathrm{sgn}(\k)},-\frac{4\m}{\mathrm{sgn}(\k)},\frac{4\e}{\mathrm{sgn}(k)}\right),$$ 
and observe that 
 \begin{align}
\Lambda_{\tilde{\k},+}^a=\frac{1}{4}(-\ep I_4 +\mu\beta-\e(\aN))+\mathit{C}^{a}_{\S}.
\end{align}
Using Proposition \ref{surlespectre}-$(\mathrm{i})$ and Lemma \ref{lemme 2.1} it follows that 
\begin{align*}
 0\in\mathrm{Sp}_{\mathrm{disc}}(\Lambda_{\k,+}^a)&\Longleftrightarrow \text{ there is } 0\neq g\in\mathit{L}^2(\S)^4: -\frac{1}{\mathrm{sgn}(\k)}(\ep I_4 -\mu\beta-\e(\aN))g= \mathit{C}^{a}_{\S}[g]\\
 &\Longleftrightarrow  \frac{4}{\mathrm{sgn}(\k)}(\ep I_4 -\mu\beta-\e(\aN))((\alpha\cdot\mathit{N})\mathit{C}^{a}_{\S})^2[g]= \mathit{C}^{a}_{\S}[g]\\
  &\Longleftrightarrow \mathit{C}^{a}_{\S}((\alpha\cdot\mathit{N})\mathit{C}^{a}_{\S})[g]=   \frac{1}{4}(\ep I_4 -\mu\beta+\e(\aN)) (\alpha\cdot\mathit{N})\mathit{C}^{a}_{\S}[g]\\
  &\Longleftrightarrow\text{ there is } 0\neq f=((\alpha\cdot\mathit{N})\mathit{C}^{a}_{\S})[g]\in\mathit{L}^2(\S)^4: \Lambda_{\tilde{\k},+}^a[f]=0\\
   &\Longleftrightarrow  a\in\mathrm{Sp}_{\mathrm{disc}}(H_{\tilde{\k}})\cap(-m,m).
 \end{align*}
Therefore, $a\in\mathrm{Sp}_{\mathrm{disc}}(H_{\k})$ if and only if $a\in\mathrm{Sp}_{\mathrm{disc}}(H_{\tilde{\k}})$, which proves $(\mathrm{iii})$. This completes the proof of the theorem.
\qed
\newline

The reader interested on confinement may wonder if the Dirac operator $H_{\k}$ generates this phenomenon under the assumption that $\mathrm{sgn}(\k)=-4$.  To clarify and provide an answer to this question, given $\varphi= u+ \Phi[g]\in\mathrm{dom}(H_{\k})$ and  recall the decomposition \ref{decompfunc}. Then
 \begin{align}\label{DC1}
 \begin{split}
  \lim\limits_{ \mathrm{nt}}\varphi_{\pm}&= t_{\S}u+  \lim\limits_{\mathrm{nt}}\Phi_{\O_\pm}[g]=  t_{\S}u+(\mathit{C}_\S \mp \frac{i}{2}(\alpha\cdot \mathit{N}))[g]\\
  &= \left(\frac{1}{4}(\ep -\m\beta -\e(\alpha\cdot \mathit{N})) \mp \frac{i}{2}(\alpha\cdot \mathit{N}) \right)g,
  \end{split}
 \end{align}
 where in the last equality we used that $t_{\S}u=-\Lambda_+[g]$. Now, multiplying the identity  \eqref{DC1} by 
 $$\left(\frac{1}{2}(\ep +\m\beta +\e(\alpha\cdot \mathit{N}))\pm i(\alpha\cdot \mathit{N})\right),$$
  we get 
   \begin{align}\label{DC2}
 \begin{split}
 \left(\frac{1}{2}(\ep +\m\beta +\e(\alpha\cdot \mathit{N})) \pm i(\alpha\cdot \mathit{N}) \right) \lim\limits_{ \mathrm{nt}}\varphi_{\pm}= \mp i\e g.
  \end{split}
 \end{align}
As consequence, if $\mathrm{sgn}(\k)=-4$ and $\e\neq0$, then $H_{\k}$ cannot generate confinement. Hence $\S$ is penetrable. Clearly, if we set $\e=0$ in \eqref{DC2}, then  $H_{\k}$ generates confinement, but we postpone this case to  subsection \ref{sub5.2}, where we establish that for Lipschitz surfaces.

\subsection{Sobolev regularity of $\dom(H_{\k})$ for $\delta$-interactions supported on the boundary of a $\mathit{C}^{1,\o}$-domain}\label{sub5.3}
In this part, we discuss  how  the smoothness of the surface supporting the singular perturbation influences  the Sobolev regularity of $\dom(H_{\k})$ in the non-critical case.  As it was shown by the author in \cite[ Section 3]{BB},   when $\S$ is a  $\mathit{C}^{2}$-smooth compact surface, we know that the functions in $\mathrm{dom}(H_\k)$ are indeed in $\mathit{H}^1(\rr^3\setminus\S)^4$. However, such a result can fails if $\S$ is less regular. Indeed, there are two obstacles which prevent us from obtaining such a result. The first one is that $(\alpha\cdot\mathit{N})\Lambda_{\k,+}[g]$ should belong to  $\mathit{H}^{1/2}(\S)^4$,  which clearly fails if for example $\S$ is $\mathit{C}^{1}$-smooth. The second reason is that we also need to extend the anticommutator  $\lbrace \alpha\cdot \mathit{N}, \mathit{C}_{\S}\rbrace$ to a bounded operator from $\mathit{L}^{2}(\S)^4$ to $\mathit{H}^{1/2}(\S)^4$. Although again, we know that behind this operator there are components of the Riesz transforms as well as the principale value of the harmonic double layer operator and its adjoint, which do not have this property, even if $\S$ is $\mathit{C}^{1,\o}$-smooth with $\o<1/2$.

 In the following, we assume that $\Omega_+$ is a bounded $\mathit{C}^{1,\o}$-smooth domain with $\g\in(0,1)$. The main result of this subsection reads as follows:
\begin{theorem}\label{self gamma} Let $\k\in\rr^3$ such that  $\mathrm{sgn}(\k)\neq0,4$ and let $H_{\k}$ be as in Theorem \ref{main4}. Then $H_{\k}$ is self adjoint and the following hold:
\begin{itemize}
  \item[(i)] If $\o\leqslant 1/2$, then  for all $s<\o$ we have 
  \begin{align*}
\mathrm{dom}(H_{\k})&\subset \left\{ u+\Phi[g]: u\in\mathit{H}^1(\rr^3)^4, g\in\mathit{H}^{s}(\S)^4, t_{\S}u=-\Lambda_{\k,+}[g]\right\}\subset\mathit{H}^{1/2+s}(\rr^3\setminus\S)^4.
\end{align*}
  \item[(ii)] If $\o>1/2$, then 
   \begin{align*}
  \mathrm{dom}(H_{\k})&=\left\{ u+\Phi[g]: u\in\mathit{H}^1(\rr^3)^4, g\in\mathit{H}^{1/2}(\S)^4, t_{\S}u=-\Lambda_{\k,+}[g]\right\}\subset\mathit{H}^1(\rr^3\setminus\S)^4.
\end{align*}
\end{itemize}
\end{theorem}

\begin{proposition}\label{normal}  There is a constant $C>0$ such that for all $x,y\in\S$, it hold that 
\begin{align}
\left| N(x)\cdot(x-y)\right| \leqslant C |x-y|^{1+\o}. 
\end{align}
\end{proposition} 

\textbf{Proof. } Since  $\S$ is $\mathit{C}^{1,\o}$-smooth, it suffies to prove the statement for $ |x-y|<1$. Without loss of generality (after translation and rotation if necessary), we may assume that $x=0$ and $N(x)=(0,0,1)$. There is a $\mathit{C}^{1,\o}$-smooth function  $\phi: B(0,1)\subset\rr^2\longrightarrow\rr$ such that $\phi(0)=0$, $|\nabla\phi(0)|=0$ and  
$$B(0,1)\cap\S =\{x=(x_1,x_2,x_3): x_3=\phi(x_1,x_2)\}.$$
Then we get 
\begin{align}
\left| \mathit{N}(x)\cdot(x-y)\right|=|y_3| =|\phi(y_1,y_2)| \leqslant C |y|^{1+\o}. 
\end{align}
Therefore the statement is proven since $\S$ is compact. \qed
\newline

In the following proposition, we prove that the anticommutator of the Cauchy operator $\mathit{C}_{\S}$ and the multiplication operator  by $(\alpha\cdot \mathit{N})$ is bounded from $\mathit{L}^{2}(\S)^4$ to $\mathit{H}^{s}(\S)^4$, for all $s\in(0,\o)$. This result should be compared to  \cite[Proposition 3.10]{BHM}, where the authors showed that  for $\S$ a $\mathit{C}^{2}$-smooth compact surface, the commutator of the Cauchy operator $\mathit{C}_{\S}$ with a H\"older continuous function of order $a\in(0,1)$ is bounded from $\mathit{L}^{2}(\S)^4$ to $\mathit{H}^{s}(\S)^4$, for all $s\in(0,a)$. In fact, both results are identical modulo a slight change of the assumptions. 

\begin{lemma}\label{encore les commutateurs} Suppose that $\S$ is  $\mathit{C}^{1,\o}$. Then,  for all $s\in(0,\o)$,  the anticommutator  $\lbrace \alpha\cdot \mathit{N}, \mathit{C}_{\S}\rbrace$ is a bounded operator from $\mathit{L}^{2}(\S)^4$ to $\mathit{H}^{s}(\S)^4$ .
\end{lemma}
\textbf{Proof. }  Let $g\in\mathit{L}^2(\S)^4$, in the same manner as in the proof of Lemma \ref{L1}, one can check that 
  \begin{align}\label{dec4}
  \begin{split}
    \lbrace \alpha\cdot \mathit{N}, \mathit{C}_{\S}\rbrace[g](x)&=\int_{y\in\S}K^{\prime}(x,y)g(y)\mathrm{dS}(y)+ \tilde{K}^{\ast}[g](x) := T_{K^{\prime}}[g](x) + \tilde{K}^{\ast}[g](x),
    \end{split}
    \end{align}
where $\tilde{K}^{\ast}$  is defined by \eqref{pv double layer} and the kernel $K^{\prime}$ is given by 
\begin{align*}
 K^{\prime}[g](x)=& \phi(x-y)(\alpha\cdot(\mathit{N}(y)-\mathit{N}(x)) -m\frac{e^{-m|x-y|}}{2i\pi |x-y|^2}(\mathit{N}(x)\cdot(x-y))\mathit{I}_4\\
 & -\frac{e^{-m|x-y|}-1}{2i\pi |x-y|^3}(\mathit{N}(x)\cdot(x-y))\mathit{I}_4.
 \end{align*}
As $\S$ is $\mathit{C}^{1,\o}$-smooth, there is a constant $C>0$ such that $\left| \mathit{N}(x)-\mathit{N}(y)\right| \leqslant C |x-y|$. Using this, the estimate \eqref{majoration}  and Proposition \ref{normal},  we obtain that $| K^{\prime}(x,y)| \leqslant C |x-y|^{-1}$. Hence the integral operator  $T_{K^{\prime}}$ is not singular, and since $ \mathit{N}$ is in the H\"older class $\mathit{C}^{0,\o}(\S)^4$,  we then can adapt the proof of \cite[Proposition 2.8]{OBV} (see also \cite[Proposition 3.10]{BHM}) and show that $T_{K^{\prime}}$ is bounded from $\mathit{L}^{2}(\S)^4$  to $\mathit{H}^{s}(\S)^4$,  for all $s\in(0,\o)$.  Finally, the fact that $\tilde{K}^{\ast}$ is bounded from $\mathit{L}^{2}(\S)^4$ onto  $\mathit{H}^{s}(\S)^4$, for all $s\in(0,\o)$,  follows by \cite[p. 165]{Tay}. \qed
\newline

We are now in a position to complete the proof of  Theorem \ref{self gamma}:

\textbf{Proof of Theorem \ref{self gamma}. } The first statement is a direct consequence of Theorem \ref{main4}. The second statement follows by the same method as in \cite[Theorem 3.1]{BB} . Indeed, fix $\g\in(0,1)$ and assume that $\S$ is $\mathit{C}^{1,\o}$, and let $g\in\mathit{L}^{2}(\S)^4$ such that  $\Lambda_{\k,+}[g]\in\mathit{H}^{1/2}(\S)^4$. Note that multiplication by $\mathit{N}$ is bounded in $\mathit{H}^{s}(\S)^4 $  for all $s\in[0,\g)$ (cf. \cite[Lemma A.2]{BHM}) and $\mathit{C}_{\S}$ is bounded from $\mathit{H}^{1/2}(\S)^4$ into itself.  Therefore, we obtain that $ \Lambda_{\k,-} \Lambda_{\k,+}[g]\in\mathit{H}^{s}(\S)^4 $, for all $s\in[0,\g)$. Now, note that   
\begin{align*}
\Lambda_{\k,-}\Lambda_{\k,+}= \frac{1}{\mathrm{sgn}(\k)}-\frac{1}{4}-\lbrace \alpha\cdot \mathit{N}, \mathit{C}_{\S}\rbrace (\alpha\cdot \mathit{N})\mathit{C}_{\S}  +\frac{2m\mu}{\mathrm{sgn}(\k)}S +\frac{\e}{\mathrm{sgn}(\k)}\lbrace \alpha\cdot \mathit{N}, \mathit{C}_{\S}\rbrace ,
\end{align*}
here we used the fact that $\mathit{C}_{\S} (\alpha\cdot \mathit{N})\lbrace \alpha\cdot \mathit{N}, \mathit{C}_{\S}\rbrace =\lbrace \alpha\cdot \mathit{N}, \mathit{C}_{\S}\rbrace (\alpha\cdot \mathit{N})\mathit{C}_{\S}$. Thus we get 
\begin{align}\label{une autre iden}
\begin{split}
g=&\frac{4(\mathrm{sgn}(\k))}{4-\mathrm{sgn}(\k)}\bigg(\Lambda_{\k,-}\Lambda_{\k,+}+\lbrace \alpha\cdot \mathit{N}, \mathit{C}_{\S}\rbrace(\alpha\cdot \mathit{N}) \mathit{C}_{\S} \\
&-\frac{\e}{\mathrm{sgn}(\k)}\lbrace \alpha\cdot \mathit{N}, \mathit{C}_{\S}\rbrace -\frac{2m\mu}{\mathrm{sgn}(\k)}S\bigg)[g],
\end{split}
\end{align}
As $\mathit{C}_{\S}(\alpha\cdot \mathit{N})$ is bounded from $\mathit{L}^{2}(\S)^4$ into itself and $S$ is bounded from $\mathit{L}^{2}(\S)^4$ to $\mathit{H}^{1/2}(\S)^4$, using Lemma \ref{encore les commutateurs},  from  \eqref{une autre iden}  it follows that $g\in\mathit{H}^{s}(\S)^4$, for all $s\in[0,\g)$. Since for all $s\in[0,1/2]$ the operator  $\Phi$ gives rise to a bounded operator $\Phi:\mathit{H}^{s}(\S)^4  \longrightarrow \mathit{H}^{1/2+s}(\rr^3\setminus\S)^4 $ (this follows by adapting the same arguments as \cite[Proposition 3.6]{BHM}),  we then get  the inclusions  in $(\mathrm{i})$. In particular, if $\o>1/2$, we then obtain that $g\in\mathit{H}^{1/2}(\S)^4$ and therefore $\Phi[g]\in\mathit{H}^{1}(\rr^3\setminus\S)^4$. This gives the equality in $(\mathrm{ii})$ and the theorem is shown.\qed
\newline

\begin{remark}\label{gammaregu} Note that if $\mathrm{sgn}(\k)\notin\{0,4\}$,  and $\S$ is $\mathit{C}^{0,\o}$-smooth with $\g\in(1/2,1)$, then using the same technique as in \cite[Theorem 3.1]{BB}, one can show that $H_{\k}$ is self-adjoint. In fact, as  $\lbrace \alpha\cdot \mathit{N}, \mathit{C}_{\S}\rbrace$ is self-adjoint, and bounded  from $\mathit{L}^{2}(\S)^4$ to $\mathit{H}^{1/2}(\S)^4$, by duality, we can extend it to a bounded operator from  $\mathit{H}^{-1/2}(\S)^4$ to $\mathit{L}^{2}(\S)^4$. Hence, by iterating twice the same argument to those of the proof of  \cite[Theorem 3.1]{BB}, we then get that 
  \begin{align*}
  \mathrm{dom}(H^{\ast}_{\k})=\left\{ u+\Phi[g]: u\in\mathit{H}^1(\rr^3)^4, g\in\mathit{H}^{1/2}(\S)^4, t_{\S}u=-\Lambda_{\k,+}[g]\right\},
  \end{align*}
  which proves the self-adjointness of $H_{\k}$ in this case.
\end{remark}

\subsection{$\delta$-interactions supported on  the boundary of a bounded uniformly rectifiable domain}\label{sub5.2} Here we discuss  special cases where we can show the self-adjointness of $H_{\k}$,  when $\O_+$ is bounded uniformly rectifiable  and $\e=0$. The idea is to identify some situations where the operator $\Lambda_{\k,+}$  gives rise to a Fredholm operator, and thereby use Theorem \ref{luis} to get the self-adjointness of $\mathcal{H}_{\k}$.  So in this subsection the domain  $\O_+$ is UR unless stated otherwise, and  we suppose that $\e=0$. Thus, $H_{\k}$ coincide with $H_{\ep,\m}$,  the Dirac operator with electrostatic and Lorentz scalar $\delta$-shell interactions supported on $\S$. The first main result on the spectral properties of  the Dirac operator $H_{\ep,\m}$ reads as follows:
\begin{theorem}\label{th5.6} Let $\ep,\m\in\rr$ such that $0<|\ep^2-\m^2|<1/\| \mathit{C}_{\S}\|^2_{\mathit{L}^2(\S)^4\rightarrow\mathit{L}^2(\S)^4} $, then $H_{\ep,\m}$ is self-adjoint. In particular, if $\O_+$ is Lipschitz and we assume  that there is  $z_0\in\cc\setminus\rr$, such that $|\ep^2-\m^2|<1/\| \mathit{C}^{z_0}_{\S}\|^2_{\mathit{L}^2(\S)^4\rightarrow\mathit{L}^2(\S)^4} $, then it holds that  
\begin{align}
\mathrm{Sp}_{\mathrm{ess}}(H_{\ep,\m})=(-\infty,-m]\cup [m,+\infty).
\end{align}
\end{theorem}
\textbf{Proof}. Fix $\ep,\m\in\rr$ such that 
$$0<|\ep^2-\m^2|<1/\| \mathit{C}^{z}_{\S}\|^2_{\mathit{L}^2(\S)^4\rightarrow\mathit{L}^2(\S)^4},$$ 
holds for some $z\in\cc\setminus\left((-\infty,-m]\cup[m,\infty)\right)$. Then, from the proof of Corollary \ref{Fredy} we have 
\begin{align}\label{encore la com}
\begin{split}
\Lambda^{z}_{(\ep,\m),\mp}\Lambda^{z}_{(\ep,\m),\pm}=\frac{1}{\ep^2-\m^2}-(\mathit{C}^{z}_{\S})^2  +\frac{2\mu}{\ep^2-\m^2}( m\mathit{I}_4 +z\beta)S^z,
\end{split}
\end{align}
Recall that $\mathit{C}^{z}_{\S}$ is bounded in $\mathit{L}^2(\S)^4$. Using  Neumann's lemma, it follows that   $$M_z:=\left( \mathit{I} -(\ep^2-\m^2)(\mathit{C}^{z}_{\S})^2\right),$$
 is a bounded invertible operator in $\mathit{L}^2(\S)^4$. Now, since $( m\mathit{I}_4 +z\beta)$ is bounded and $S$ is compact on $\mathit{L}^2(\S)^4$,  we therefore get that $K_z:=\frac{2\mu}{\ep^2-\m^2}( m\mathit{I}_4 +z\beta)S^z$ is compact on $\mathit{L}^2(\S)^4$. Combining this with  \eqref{encore la com}, we obtain  that 
\begin{align}
\begin{split}
\mathit{I}-(\ep^2-\m^2)M_z^{-1}\Lambda^z_{(\ep,\m),-}\Lambda^z_{(\ep,\m),+}&=-(\ep^2-\m^2)M_z^{-1}K_z,\\
\mathit{I}-(\ep^2-\m^2)\Lambda^z_{(\ep,\m),+}\Lambda^z_{(\ep,\m),-}M_z^{-1}&=-(\ep^2-\m^2)K_zM_z^{-1}.
\end{split}
\end{align}
As $M_z^{-1}\Lambda^z_{(\ep,\m),-}$ and $\Lambda^z_{(\ep,\m),-}M_z^{-1}$ are bounded operators on $\mathit{L}^2(\S)^4$, $M_z^{-1}K_z$ and $K_zM_z^{-1}$ are compact on $\mathit{L}^2(\S)^4$, then \cite[Theorem 1.50 and Theorem 1.51]{Aiena}  yields that $\Lambda^{z}_{(\ep,\m),+}$ is Fredholm. Hence, the first statement is a direct consequence of Theorem \ref{luis} and the fact that $\Lambda_{(\ep,\m),+}$ is a self-adjoint,  Fredholm operator on $\mathit{L}^2(\S)^4$. The last statement readily  follows  by Proposition \ref{surlespectre} and the fact that $\Lambda^{z_0}_{(\ep,\m),+}$ is Fredholm by assumption. \qed

\begin{remark} From Lemma  \ref{lemme 2.1}($\mathrm{iv}$) it easily follows that $\| \mathit{C}^{z}_{\S}\| \geqslant 1/2$ (cf. \cite[Remark 3.5]{AMV2}), which implies that  $|\ep^2-\m^2|<4$. Hence, the combination of coupling constants $\ep$ and $\m$ is not  critical. Of course, we already know that the above result is false in the case  $\ep^2-\m^2=4$. Note that Theorem \ref{th5.6} remains valid if one control the norm of the Cauchy operator instead of controlling the combination of interactions. However, this can influence the geometrical characterization of $\S$ which can imply an increase in terms of regularity.
\end{remark}
As it was mentioned in the introduction, for  $m=0$ and $\m\in(-2,2)$, in \cite{PV} it was shown the existence of a unique self-adjoint realization of the two dimensional Dirac operator with pure Lorentz scalar $\delta$-shell interactions, where  $\S$ is a closed curve with finitely many corners. It seems that their assumption (i.e the restriction on $\m$ to be in $(-2,2)$) is related to the assumption that we imposed in the Theorem \ref{th5.6}.

 Although Theorem \ref{th5.6} gives an upper bound for $|\ep^2-\m^2|$ so that $H_{\ep,\m}$ is self-adjoint, this is not  satisfactory in the sense that this bound involves $\| \mathit{C}_{\S}\|^2_{\mathit{L}^2(\S)^4\rightarrow\mathit{L}^2(\S)^4} $, which is not easy to quantify. In what follow, we are going to remove this restriction by imposing a sign hypothesis on the coupling constants,  and give a better quantitative assumption than the one of Theorem \ref{th5.6}. The next theorem makes this more precise. 
 
\begin{theorem}\label{th5.7}Let $\ep,\m\in\rr$ such that $|\ep|\neq|\m|$, and let  $(H_{\ep,\m},\mathrm{dom}(H_{\ep,\m}))$ be as above.  Assume that $\ep$ and $\m$ satisfy one of the the following assumptions:\begin{itemize}
  \item [(a)] $\m^2>\ep^2$.
   \item [(b)]  $\ep^2>\m^2$  and $16\| W\|^2_{\mathit{L}^2(\S)^2\rightarrow\mathit{L}^2(\S)^2}<\ep^2-\m^2<1/\| W\|^2_{\mathit{L}^2(\S)^2\rightarrow\mathit{L}^2(\S)^2} $.
\end{itemize}
Then $H_{\ep,\m}$ is self-adjoint. In particular, if $\O_+$ is Lipschitz, then the following statements hold true:
\begin{itemize}
  \item[(i)]  Given $a\in(-m,m)$, then $\mathrm{Kr}(H_{\ep,\m}-a)\neq0$ $\Longleftrightarrow$ $\mathrm{Kr}( \Lambda^{a}_{(\ep,\m),+})\neq0$.
\item[(ii)] For all $z\in\cc\setminus\rr$, it holds that 
 \begin{align*}
    (H_{\ep,\m} -z)^{-1}= (H -z)^{-1} - \Phi^{z}(\Lambda^{z}_{(\ep,\m),+} )^{-1} t_\S(H-z)^{-1}.
\end{align*}
  \item[(iii)] $\mathrm{Sp}_{\mathrm{ess}}(H_{\ep,\m})=(-\infty,-m]\cup [m,+\infty)$. 
  \item [(iv)]$ a\in\mathrm{Sp}_{\mathrm{p}}(H_{\ep,\m})$ if and only if $a\in\mathrm{Sp}_{\mathrm{p}}(H_{\frac{-4\ep}{\ep^2-\m^2},\frac{-4\m}{\ep^2-\m^2}})$, for all $a\in(-m,m)$.
    \item [(v)] $C_0:=\sup_{a\in[-m,m]}\Vert  \mathit{C}^a_\S\Vert <\infty$. Moreover, $\mathrm{Sp}_{\mathrm{disc}}(H_{\ep,\m})\cap(-m,m)=\emptyset$  either if $| \ep-\m|<1/C_0$ and $| \ep+\m|<1/C_0$, or if $| \ep-\m|>4C_0$ and $| \ep+\m|>4C_0$.
\end{itemize}
\end{theorem}
\textbf{Proof.} To prove the theorem, in both situations, we show that $\Lambda^{z}_{(\ep,\m),+}$ is Fredholm  for all $z\in\cc\setminus\left((-\infty,-m]\cup[m,\infty)\right)$. Once this is shown, we use the fact that $\Lambda_{(\ep,\m),+}$ is a bounded self-adjoint operator, and we conclude by using Theorem \ref{luis} to obtain the first statement of the theorem.  So, fix $\l,\m\in\rr$ such that $|\ep|\neq|\m|$, and let $z\in\cc\setminus\left((-\infty,-m]\cup[m,\infty)\right)$. Then, from the definition of $\mathit{C}^{z}_{\S}$ it follows that 
\begin{align}\label{zCau}
\mathit{C}^{z}_{\S}= T^{z}_{K}+ \begin{pmatrix}
 0  & W\\
 W& 0
\end{pmatrix}:= T^{z}_{K} +\tilde{W},
\end{align}
where the kernel $K$ satisfies  
\begin{align}\label{estimation8}
\sup_{1\leqslant k,j\leqslant 4}\left| K(x-y)\right| &=\mathcal{O}(|x-y|^{-1})\quad \text{ when }  |x-y|\longrightarrow 0.
\end{align}
Hence, $T^{z}_{K}$ is compact in $\mathit{L}^2(\S)^4$. Therefore, in the same way as in \eqref{encore la com} we get that 
\begin{align}\label{encore la com2}
\begin{split}
\Lambda^{z}_{(\ep,\m),\mp}\Lambda^{z}_{(\ep,\m),\pm}&=\frac{1}{\ep^2-\m^2}- \tilde{W}^2 +  (T^{z}_{K})^2 +\{T^{z}_{K},\tilde{W}\} +\frac{2\mu}{\ep^2-\m^2}( m\mathit{I}_4 +z\beta)S^z\\
&:=\frac{1}{\ep^2-\m^2}- \tilde{W}^2 +T_K,
\end{split}
\end{align}
where $T_K$ is compact in $\mathit{L}^2(\S)^4$. Observe that 
\begin{align*}
\tilde{W}^2=  \begin{pmatrix}
 W^2  & 0\\
 & W^2
\end{pmatrix}.
\end{align*}
Now, if $\ep^2<\m^2$ holds, then using that $W$ is a bounded self-adjoint operator in $\mathit{L}^2(\S)^2$, it follows that   $\tilde{W}^2$ is a nonnegative, self-adjoint operator on $\mathit{L}^2(\S)^4$. From this, it follows that  $1/(\ep^2-\m^2)$ belongs to the resolvent set of $\tilde{W}^2$ and hence $\mathit{I}_4- (\ep^2-\m^2)\tilde{W}^2$ in invertible on $\mathit{L}^2(\S)^4$ .  In another hand, assume that  $\ep^2>\m^2$ and $\ep^2-\m^2<1/\| W\|^2_{\mathit{L}^2(\S)^2\rightarrow\mathit{L}^2(\S)^2} $ hold, then Neumann's lemma yields that $\mathit{I}_4- (\ep^2-\m^2)\tilde{W}^2$ is invertible on $\mathit{L}^2(\S)^4$. In both cases, similar arguments to those of the proof of Theorem \ref{th5.6} yield that $\Lambda^z_{(\ep,\m),+}$ is Fredholm, which proves the first statement of the theorem for this two cases. Now we deal with the case $16\| W\|^2_{\mathit{L}^2(\S)^2\rightarrow\mathit{L}^2(\S)^2}<\ep^2-\m^2$.  From Lemma \ref{lemme 2.1}, we know  that $W$ is invertible on $\mathit{L}^2(\S)^2$ and $W^{-1}= -4(\sigma\cdot\mathit{N})  W (\sigma\cdot\mathit{N})$.  Thus, from \eqref{encore la com2} it follows that
\begin{align}\label{encore la com3}
\begin{split}
(\ep^2-\m^2)(\tilde{W}^{-1})^2\Lambda^{z}_{(\ep,\m),-}\Lambda^{z}_{(\ep,\m),+}&=(\tilde{W}^{-1})^2-(\ep^2-\m^2)\mathit{I}_4  +(\ep^2-\m^2)(\tilde{W}^{-1})^2 T_K,\\
(\ep^2-\m^2)\Lambda^{z}_{(\ep,\m),+}\Lambda^{z}_{(\ep,\m),-}(\tilde{W}^{-1})^2&=(\tilde{W}^{-1})^2-(\ep^2-\m^2)\mathit{I}_4  +(\ep^2-\m^2)T_K(\tilde{W}^{-1})^2.
\end{split}
\end{align}
As $\| \tilde{W}\|_{\mathit{L}^2(\S)^4\rightarrow\mathit{L}^2(\S)^4} = \| W\|_{\mathit{L}^2(\S)^2\rightarrow\mathit{L}^2(\S)^2} $,   using again Lemma \ref{lemme 2.1}  we get that 
\begin{align}\label{une norme}
\| \tilde{W}^{-1}\|_{\mathit{L}^2(\S)^4\rightarrow\mathit{L}^2(\S)^4}  \leqslant 4\| \tilde{W}\|_{\mathit{L}^2(\S)^4\rightarrow\mathit{L}^2(\S)^4}=4\| W\|_{\mathit{L}^2(\S)^2\rightarrow\mathit{L}^2(\S)^2}
\end{align} 
Hence, if   $\ep^2-\m^2>4\| W\|_{\mathit{L}^2(\S)^2\rightarrow\mathit{L}^2(\S)^2}$, then $\ep^2-\m^2>\| \tilde{W}^{-1}\|^2_{\mathit{L}^2(\S)^4\rightarrow\mathit{L}^2(\S)^4}$. Thus $\ep^2-\m^2$ is not in the spectrum of $(\tilde{W}^{-1})^2$.  Thereby  $ \tilde{W}^{-1}-(\ep^2-\m^2)\mathit{I}_4 $ is invertible on $\mathit{L}^2(\S)^4$. Now, from \eqref{encore la com3} it follows that 
\begin{align}\label{encore la com4}
\begin{split}
\mathit{I}_4  - (\ep^2-\m^2)\left((\tilde{W}^{-1})^2-(\ep^2-\m^2)\mathit{I}_4\right)^{-1}(\tilde{W}^{-1})^2\Lambda^{z}_{(\ep,\m),-}\Lambda^{z}_{(\ep,\m),+}&= T_{K_1},\\
\mathit{I}_4  - (\ep^2-\m^2)\Lambda^{z}_{(\ep,\m),+}\Lambda^{z}_{(\ep,\m),-}(\tilde{W}^{-1})^2\left((\tilde{W}^{-1})^2-(\ep^2-\m^2)\mathit{I}_4\right)^{-1}&= T_{K_2}.
\end{split}
\end{align}
where $ T_{K_1},T_{K_2}\in\mathcal{K}(\mathit{L}^2(\S)^4)$. Thereby,  \cite[Theorem 1.50 and Theorem 1.51]{Aiena}  yields that $\Lambda^{z}_{(\ep,\m),+}$ is Fredholm, and this finishes the proof of the first and the second statements. Item $(\mathrm{i})$ is consequence of Proposition \ref{surlespectre}. The proof of items $(\mathrm{ii})$, $(\mathrm{iii})$ and $(\mathrm{iv})$  runs as in the proof of Theorem \ref{main4}.  Now we turn to the proof of item $(\mathrm{v})$, the first claim of statement  is contained in \cite[Lemma 3.2]{AMV2} and \cite[Proposition 3.5]{BEHL1}, for a $\mathit{C}^2$-compact surface $\S$, and the same arguments hold true in the Lipschitz case. To see the last claim of $(\mathrm{v})$, note that for all $a\in(-m,m)$, we have
\begin{align}\label{compiir}
0\in \mathrm{Sp}_{\mathrm{disc}}(\Lambda^a_{(\ep,\m),+}) \Longleftrightarrow  -1\in \mathrm{Sp}_{\mathrm{disc}}( (\ep\mathit{I}_4+\m\beta) \mathit{C}^a_\S).
\end{align}
Using the first claim, it follows that  if $| \ep-\m|<1/C_0$ and $| \ep+\m|<1/C_0$, then 
\begin{align*}
\| (\ep\mathit{I}_4+\m\beta) \mathit{C}^a_\S\|_{\mathit{L}^2(\S)^4\rightarrow\mathit{L}^2(\S)^4}<1.
\end{align*}
Therefore, $-1\notin \mathrm{Sp}_{\mathrm{disc}}( (\ep\mathit{I}_4+\m\beta) \mathit{C}^a_\S)$. Hence,  \eqref{compiir} and $(\mathrm{i})$ yield that $\mathrm{Sp}_{\mathrm{disc}}(H_{\ep,\m})\cap(-m,m)=\emptyset$. Using the equivalence given by $(\mathrm{iv})$, and  iterating  the previous arguments we easily recover the case $| \ep-\m|>4C_0$ and $| \ep+\m|>4C_0$,  which gives $(\mathrm{v})$.   This completes the proof of theorem.\qed
\newline

\begin{remark} Assume that $\O_+$ is Lipschitz. Then, following essentially the same arguments as in Theorem \ref{th5.7},  one can show that, if for all $a\in(-m,m)$ the following holds:
$$16\| \mathit{C}^{a}_{\S}\|^2_{\mathit{L}^2(\S)^4\rightarrow\mathit{L}^2(\S)^4}<\ep^2-\m^2<1/\| \mathit{C}^{a}_{\S}\|^2_{\mathit{L}^2(\S)^4\rightarrow\mathit{L}^2(\S)^4}, $$
then $H_{\ep,\m}$ is self-adjoint.  Moreover, if $\m=0$, then from Theorem  \ref{th5.7}-$(\mathrm{v})$ it follows that 
\begin{align*}
\mathrm{Sp}_{\mathrm{disc}}(H_{\ep,0})\cap(-m,m)=\emptyset,
\end{align*}
 see also \cite[Theorem 3.3]{AMV2} and \cite[Theorem 4.4]{BEHL1} for a similar result.
\end{remark}
\begin{remark} Assume that $\O_+$ is UR. Then, using exactly the same technique as in the proof of Theorem \ref{th5.7}, one can show that the coupling $H+\ept\g_5\delta_{\S}$ is self-adjoint under the assumption  $(\mathrm{b})$,  with $\m=0$.
\end{remark}
\begin{remark} Note that in Theorem \ref{th5.7} $(\mathrm{b})$, the combination of the coupling constants $\ep$ and $\m$ is not critical. Moreover, there is an interval $J\subset\rr_{+}$, such that we have no information on the self-adjointness character of $H_{\ep,\m}$, if $\ep^2-\m^2\in J$. 
\end{remark}


Next, we discuss the particular case $\ep^2-\m^2=-4$.  Assume that $\O_+$ is Lipschitz, then it is clear that $(\mathrm{P4})$ holds true. Thus, if we let 
$$P_\pm= \left( 1\mp i\frac{i}{2}(\ep +\m\beta )(\alpha\cdot \mathit{N}) \right),$$
 as consequence of proposition \ref{confinementgenerale}, we then have the following result.
\begin{proposition}\label{Un confinement}Assume that $\O_+$ is Lipschitz. Let $\ep,\m\in\rr$ such that $\ep^2-\m^2=-4$,  and let $ H_{\ep,\m}$ be as in Theorem \ref{th5.7}. Then $\S$ is impenetrable and it holds that
\begin{align}
H_{\ep,\m}\varphi=H_{\ep,\m}^{\Omega_+}(\varphi_+)\oplus H_{\ep,\m}^{\Omega_-}(\varphi_-)=\left(-i\alpha\cdot \nabla+m\beta\right)\varphi_+\oplus\left(-i\alpha\cdot \nabla+m\beta\right)\varphi_-
\end{align}
where $H_{\ep,\m}^{\Omega_\pm}$ are the self-adjoint Dirac operators defined on
\begin{align*}
\mathrm{dom}(H_{\ep,\m}^{\Omega_\pm})&=\left\{ \varphi_{\pm} :=u_{\Omega_\pm}+\Phi_{\O_\pm}[g], u_{\Omega_\pm}\in\mathit{H}^1(\Omega_\pm)^4, g\in\mathit{L}^2(\S)^4: P_\pm \lim\limits_{ \mathrm{nt}}\varphi_{\pm}=0 \right\}.
\end{align*}
\end{proposition}

\begin{remark} By Taking $\ep=0$ in Theorem \ref{th5.7} $(\mathrm{a})$, we conclude that if $\O_+$ is UR, then  $H_{0,\m}$  is always self-adjoint.  Moreover, $H_{0,\m}$ generates confinement when $\m=\pm2$, for any compact Lipschitz surface $\S$. 
\end{remark}

 The reason we assumed that $ \eta = 0 $ is purely  technical. The  following proposition is about the self-adjointness of   the coupling $H+\e(\alpha\cdot\mathit{N})\delta_\S$.
\begin{proposition}\label{eta tout } Assume that $\O_+$ is UR. Let $\e\in\rr\setminus\{0\}$, set $\k=(0,0,\e)$ and let $ H_{\k}$ be as above. Then $ H_{\k}$ is self-adjoint and we have 
\begin{align*}
  \mathrm{dom}(H_{\k})=\left\{ u+\Phi[-4\e^2(\e+4)^{-1}(\alpha\cdot\mathit{N})\Lambda_{\k,-}(\alpha\cdot\mathit{N})[T_{\S}u]]: u\in\mathit{H}^1(\rr^3)^4 \right\}.
  \end{align*}
Moreover, If $\O_+$ is Lipschitz, then the spectrum of $ H_{\k}$ is given by
\begin{align}\label{speceta}
\mathrm{Sp}(H_{\k})=\mathrm{Sp}_{\mathrm{ess}}(H_{\k})=(-\infty,-m]\cup [m,+\infty).
\end{align}
\end{proposition}
\textbf{Proof.} Assume that $\e\in\rr\setminus\{0\}$  and fix $z\in\cc\setminus\left((-\infty,-m]\cup[m,\infty)\right)$. Recall that $\Lambda^{z}_{\k,\pm}$ are given by
\begin{align*}
\Lambda^{z}_{\k,\pm}=\frac{1}{\e}(\alpha\cdot N) \pm\mathit{C}^{z}_{\S}.
\end{align*} 
Now, using Lemma \ref{lemme 2.1}, a simple computation yields
\begin{align*}
(\e(\alpha\cdot\mathit{N}))\Lambda^{z}_{\k,-}(\e(\alpha\cdot\mathit{N}))\Lambda^{z}_{\k,+}=\Lambda^{z}_{\k,+}(\e(\alpha\cdot\mathit{N}))\Lambda^{z}_{\k,-}(\e(\alpha\cdot\mathit{N}))=1+ \frac{\e^2}{4}.
\end{align*}
Therefore, $\Lambda^{z}_{\k,+}$ is invertible with $(\Lambda^{z}_{\k,+})^{-1}=4\e^2(\e+4)^{-1}(\alpha\cdot\mathit{N})\Lambda^{z}_{\k,-}(\alpha\cdot\mathit{N})$. In particular, $\Lambda^{z}_{\k,+}$ is Fredholm, for all $z\in\cc\setminus\left((-\infty,-m]\cup[m,\infty)\right)$. As $\Lambda_{\k,+}$ is bounded invertible, and self-adjoint in $\mathit{L}^2(\S)^4$, using  Theorem \ref{luis}  we then  get the first statement.  That $\mathrm{Sp}(H_{\k})$ is characterized by  \eqref{speceta} is a consequence of Proposition \ref{surlespectre}. This completes the proof of the proposition.\qed
\newline

To finish this part, we briefly discuss  the particular case $\mu=\pm\ep$. Assume that $\O_+$ is Lipschitz and  given $\ep\in\rr\setminus\{0\}$, recall that $H_{\ep,\pm\ep}$ is defined on the domain
\begin{align}
\mathrm{dom}(H_{\ep,\pm\ep})=\left\{ u+\Phi[g]: u\in\mathit{H}^1(\rr^3)^4, g\in P_\pm\mathit{L}^{2}(\S)^4 \text{ and }P_\pm t_{\S}u=-P_\pm \Lambda_{+}[g]\right\},
\end{align}
where $\Lambda_{\ep,+}$ are given by
\begin{align}\label{lambda2}
\Lambda_{\ep,+}=P_\pm\left(1/2\ep+\mathit{C}_{\S}\right)P_\pm \quad \text{and } \Lambda_{\ep,-}=P_\pm\left(1/2\ep -\mathit{C}_{\S}\right)P_\pm, \quad\forall \ep \neq0.
\end{align}
see \cite[Proposition 3.5]{BB} for example.  Observe that 
 \begin{align}
\Lambda_{\ep,-}\Lambda_{\ep,+}=  \frac{1}{4\ep^2}P_\pm -\frac{m^2}{2\ep}(S)^2P_\pm.
\end{align}
Since $S$ is bounded from $\mathit{H}^{-1/2}(\S)^4$ to $\mathit{H}^{1/2}(\S)^4$ and hence compact on $\mathit{L}^{2}(\S)^4$ (see \cite[ Theorem 6.11]{Mc} for example), we then get the following result.
\begin{proposition} Let  $\ep\in\rr\setminus\{0\}$ and assume that $\S$ is a compact Lipschitz surface. Then $ H_{\ep,\pm\ep}$ is self-adjoint and 
 \begin{align*}
  \mathrm{dom}(H_{\ep,\pm\ep})=\left\{ u+\Phi[g]: u\in\mathit{H}^1(\rr^3)^4, g\in P_\pm\mathit{H}^{1/2}(\S)^4, P_\pm t_{\S}u=-P_\pm\Lambda_+[g]\right\}.
  \end{align*}
\end{proposition} 
\textbf{Proof.} This readily follows from the compactness and regularization property of the operator $S$ and Theorem \ref{luis}. \qed

 \subsection{Spectral properties of $H_{\mt}$ and $H_{\ut}$}\label{sub4.4} In this part, we briefly discuss  the spectral properties of the Dirac operators  $H_{\mt}$ and $H_{\ut}$ defined by \eqref{autresoperator}. Recall that $\Lambda^{z}_{\mt,+}$ and $\Lambda^{z}_{\ut,+}$ are given by
 \begin{align}
\label{}
    \Lambda^{z}_{\mt,\pm}=  \frac{i}{\mt}\g_5\beta\pm \CS^{z}\, \text{ and }\,\Lambda^{z}_{\ut,\pm}=  \frac{1}{\ut}\g_5\beta(\aN)\pm \CS^{z},
\end{align}
  that is  $A_{\bullet}=\bullet\tilde{A}_{\bullet}$, for $\bullet=\mt$ or $\ut$. 
  
  The following two propositions gather the most important spectral properties of $H_{\mt}$. We remark that $H_{\mt}$ has almost the same properties as the coupling $(H+\m\beta\delta_{\S})$, and this is the reason why we have called the potential $V_{\mt}$  the \textit{ modified Lorentz scalar} $\delta$-potential. 
 \begin{proposition}\label{Th4.8} Let $\mt\in\rr\setminus\{0\}$ and assume that $\O_+$ is UR, then,  $(H_{\mt},\dom(H_{\mt}))$ is self-adjoint. In particular, if $\O_+$ is Lipschitz, then the following hold:
\begin{itemize}
\item[(i)] $\mathrm{Sp}_{\mathrm{ess}}(H_{\mt})=(-\infty,-m]\cup [m,+\infty)$.
  \item[(ii)] If $\mt^2=4$, then $H_{\mt}$  generates confinement and we have 
$$H_{\mt}\varphi=H_{\mt}^{\Omega_+}\varphi_{+}\oplus H_{\mt}^{\Omega_-}\varphi_{-}=\left(-i\alpha\cdot \nabla+m\beta\right)\varphi_{+}\oplus\left(-i\alpha\cdot \nabla+m\beta\right)\varphi_{-},$$
    where $H_{\mt}^{\Omega_\pm}$ are the self-adjoint Dirac operators defined on
\begin{align*}
\dom(H_{\mt}^{\Omega_\pm})=\bigg\{ u_{\Omega_\pm}+\Phi_{\O_\pm}[g]&: u_{\Omega_\pm}\in\mathit{H}^1(\Omega_\pm)^4, g\in\mathit{L}^2(\S)^4 \text{ and }\\
& \left(\frac{1}{2} \mp \frac{1}{\mt}\g_5\beta(\aN)\right)(t_\S u_{\Omega_\pm}+\mathit{C}_\pm[g])=0 \bigg\}.
\end{align*}
\end{itemize} 
 \end{proposition}
 \textbf{Proof.} Fix  $z\in\cc\setminus\left((-\infty,-m]\cup[m,\infty)\right)$, and observe that 
 \begin{align}
\label{}
 \{\g_5\beta,\CS^z\}= 2z\g_5\beta S^{z}.
\end{align}
 Now, using \eqref{zCau}, similar arguments   as in the proof of Theorem \ref{th5.7} yield that  
 \begin{align}\label{}
\begin{split}
(\Lambda^{z}_{\mt,+})^2=\frac{1}{\mt^2}+ \tilde{W}^2 +  \frac{i}{\mt}\{\g_5\beta,\CS^z\} + T_z,
\end{split}
\end{align}
 where $T_z\in\mathcal{K}(\mathit{L}^2(\S)^4)$. Thus $(\Lambda^{z}_{\mt,+})$ is Fredholm, for all $z\in\cc\setminus\left((-\infty,-m]\cup[m,\infty)\right)$. Therefore, $H_{\mt}$ is self-adjoint  by Theorem \ref{luis}. Assertion $(\mathrm{i})$ follows by Proposition \ref{surlespectre}. Now it is easy to check that $(1/2 \mp \frac{1}{\mt}\g_5\beta(\aN))$ are projectors, thus the property $(\mathrm{P3})$ holds true. Therefore,   $(\mathrm{ii})$ is a direct consequence of Proposition \ref{confinementgenerale}. This completes the proof of the proposition. \qed
 \newline

 The following proposition gives us more information about the spectrum of $H_{\mt}$ in the case of $\mathit{C}^{1,\o}$ domains. The arguments of the proof are rather standard, so we are not going to give a complete
proof.
 \begin{proposition} Assume that $\O_+$ is $\mathit{C}^{1,\o}$-smooth with $\g>1/2$, and let $H_{\mt}$ be as in Proposition \ref{Th4.8}. Then, the following is true:
  \begin{itemize}
\item[(i)] The spectrum of $H_{\mt}$ is symmetric with respect to $0$.
\item[(ii)] $\mathrm{Sp}_{\mathrm{disc}}(H_{\mt})\cap(-m,m)$ is finite, and each  eigenvalue of $H_{\mt}$ has an even multiplicity.
  \item[(iii)] $H_{\mt}$ is unitarily equivalent to $H_{-\mt}$.
  \item [(iv)]$ a\in\mathrm{Sp}_{\mathrm{p}}(H_{\mt})$ if and only if $a\in\mathrm{Sp}_{\mathrm{p}}(H_{\frac{-4}{\mt}})$, for all $a\in(-m,m)$.
  \item[(v)]  There is $C_0>0$ such that  $\mathrm{Sp}_{\mathrm{disc}}(H_{\mt})\cap(-m,m)=\emptyset$  either if $| \mt|<1/C_0$ or if $| \mt|>4C_0$.
\end{itemize} 
 \end{proposition}
 \textbf{Proof.}  First, observe that for all $\mt\in\rr\setminus\{0\}$, $\dom(H_{\mt})\subset\mathit{H}^1(\rr^3\setminus\S)^4$ (this follows in the same way as in Theorem \ref{self gamma}). Moreover, $H_{\mt}$ acts in the sense of distributions as  
 \begin{align}
H_{\mt}\varphi= \left(-i\nabla\cdot\alpha+m\beta\right)\varphi_+\oplus\left(-i\nabla\cdot\alpha+m\beta\right)\varphi_-,
 \end{align}
 on the domain
   \begin{align*}
 \begin{split}
\dom(H_{\mt})=\bigg\{ \varphi=(\varphi_+,\varphi_-)&\in\mathit{H}^1(\Omega_+)^4\oplus\mathit{H}^1(\Omega_-)^4 : \\
& \left(\frac{1}{2} - \frac{1}{\mt}\g_5\beta(\aN)\right)t_{\S}\varphi_+ =- \left(\frac{1}{2} + \frac{1}{\mt}\g_5\beta(\aN)\right)t_{\S}\varphi_- \bigg\}.
\end{split}
\end{align*}  
 Now, assertions $(\mathrm{i})$ and the fact that each  eigenvalue of $H_{\mt}$ has an even multiplicity can be proved as much as \cite[Theorem 2.3]{HOP}. Also, the fact that $\mathrm{Sp}_{\mathrm{disc}}(H_{\mt})\cap(-m,m)$ is finite can  be deduced by applying the same arguments as \cite[Theorem 4.1]{BEHL2}. Let us prove $(\mathrm{iii})$, for that given $\psi\in\mathit{L}^{2}(\rr^3)^4$ and define the operator
\begin{align}\label{matt}
 \quad T(\psi)=\gamma_5\beta \psi.
\end{align}
Then, a simple computation yields that $T^2_2(\psi)=-\psi$ and $T(H(\psi))=-H(T(\psi))$.  Moreover, one can check easily  that 
\begin{align*}
  \left(\frac{1}{2} \mp \frac{1}{\mt}\g_5\beta(\aN)\right)\left(\g_5\beta t_{\S}\varphi_+\right)&=-\g_5\beta \left(\frac{1}{2} \pm \frac{1}{\mt}\g_5\beta(\aN)\right)t_{\S}\varphi_\pm.
 \end{align*}
Hence, we conclude that  $\varphi\in\dom(H_{\mt})$ if and only if $ T(\varphi)\in\dom(H_{-\mt})$, which proves $(\mathrm{iii})$.  Assertions  $(\mathrm{iv})$ and $(\mathrm{v})$ can be proved in the same way as Theorem \ref{th5.7}$(\mathrm{v})$-$(\mathrm{vi})$.\qed
\newline

 We now move on to the spectral study of the operator $H_{\ut}$. Again, we observe that $H_{\ut}$  has almost the same spectral properties as $(H+i\u\beta(\aN)\delta_{\S})$. In the following proposition, we are interested  only  on the self-adjointness character of $H_{\ut}$, we omit the other specific spectral properties  because they can be deduced  from \cite[Theorem 5.1 and Theorem 5.2]{BB}.
\begin{proposition}\label{Propo4.9} Let $\ut\in\rr\setminus\{0\}$ and let  $(H_{\ut},\dom(H_{\mt}))$ be as in \eqref{dom1}. Then following hold true:
\begin{itemize}
\item[(i)]  If $\O_+$ is a regular SKT domain and $\mt^2\neq4$, then $(H_\l, \dom(H_\l))$ is self-adjoint. 
\item[(ii)] If $\O_+$ is a $\mathit{C}^2$-smooth domain and $\mt^2=4$, then $(H_{\ut}, \dom(H_{\ut}))$ is essentially self-adjoint and generates confinement and we have 
$$\overline{H_{\ut}}\varphi=H_{\ut}^{\Omega_+}\varphi_+\oplus H_{\ut}^{\Omega_-}\varphi_-=\left(-i\alpha\cdot \nabla+m\beta\right)\varphi_+\oplus\left(-i\alpha\cdot \nabla+m\beta\right)\varphi_-,$$
    where $H_{\ut}^{\Omega_\pm}$ are the self-adjoint Dirac operators defined on
\begin{align*}
\dom(H_{\ut}^{\Omega_\pm})=\bigg\{ \varphi_\pm&\in\mathit{L}^2(\Omega_\pm)^4:  (\alpha\cdot\nabla) \varphi_\pm\in\mathit{L}^2(\Omega_{\pm})^4  \text{ and } \left(\frac{1}{2} \pm \frac{i}{\mt}\g_5\beta\right)t_\S\varphi_\pm =0 \bigg\}.
\end{align*}
\end{itemize} 
\end{proposition}
\textbf{Proof.} Let $z\in\cc\setminus\left((-\infty,-m]\cup[m,\infty)\right)$, then 

  \begin{align}\label{}
\begin{split}
\Lambda^{z}_{\ut,\mp}\Lambda^{z}_{\ut,\pm}=\frac{1}{\ut^2}- (\CS^z)^2 \pm  \frac{1}{\ut}[\g_5\beta(\aN),\CS^z],
\end{split}
\end{align}
Now, observe that 
 \begin{align}\label{}
\begin{split}
 \frac{1}{\ut}[\g_5\beta(\aN),\CS^z]=m[\g_5(\aN),S^z] +\g_5\beta (T_z+ \{\aN, \tilde{W}\}),
\end{split}
\end{align}
where $\tilde{W}$ is given by \eqref{zCau},  and $T_z$ is an integral operator with  kernel $K_z$ given by:
\begin{align*}
 K_z(x,y)=&\sqrt{z^2-m^2} \frac{e^{i\sqrt{z^2-m^2}|x-y|}}{4\pi |x-y|^2}\left((\alpha\cdot\mathit{N}(x))\alpha\cdot(x-y) + \alpha\cdot(x-y)(\aN(y))\right)\\
 &+ \frac{e^{i\sqrt{z^2-m^2}|x-y|} -1}{4\pi |x-y|^3}\left[(\alpha\cdot\mathit{N}(x))(i \alpha\cdot(x-y)) +   \left(i \alpha\cdot(x-y)\right)\left(\alpha\cdot\mathit{N}(y)\right)\right].
 \end{align*}
Clearly, $T_z, [\g_5(\aN),S^z]\in\mathcal{K}(\mathit{L}^{2}(\S)^4)$. Thus, using  Proposition \ref{remarkSKT} we get that 
  \begin{align}\label{}
\begin{split}
\Lambda^{z}_{\ut,-}\Lambda^{z}_{\ut,+}=\frac{4-\ut^2}{4\ut^2}- \tilde{T}_z,
\end{split}
\end{align}
with $\tilde{T}_z\in\mathcal{K}(\mathit{L}^{2}(\S)^4)$. Thus, if $\O_+$ is a regular SKT domain and $\mt^2\neq4$, then similar arguments to those of the proof of Theorem  \ref{th5.7} yield that $\Lambda^{z}_{\ut,+}$ is Fredholm, for all $z\in\cc\setminus\left((-\infty,-m]\cup[m,\infty)\right)$. Therefore, $(\mathrm{i})$ follows by Theorem \ref{luis}. 

Now we are going to prove $(\mathrm{ii})$, we only consider the case $\ut=2$, since the case $\ut=-2$ can be treated analogously. So assume that $\O_+$ is a $\mathit{C}^2$-smooth, then it is clear that  $\ut=2$ is a critical parameter, and $(\mathrm{P1})$ holds true.  Thus, $\Lambda_{2,\pm}$ extends to a bounded operator from $\mathit{H}^{-1/2}(\S)^4$ onto itself by Proposition \ref{extension}.  Now observe that 
  \begin{align}\label{}
\tilde{\Lambda}_{2,+}(\frac{1}{2}\g_5\beta(\aN))\tilde{\Lambda}_{2,-}= \Lambda_{2,+}\tilde{\Lambda}_{2,+}\tilde{\Lambda}_{2,-}+ \tilde{\Lambda}_{2,+}\tilde{\Lambda}_{2,-}\tilde{\Lambda}_{2,-}.
\end{align}
Then, as in \cite[Lemma 3.1]{BB} one can show that $\tilde{\Lambda}_{2,+}\tilde{\Lambda}_{2,-}$ and  $\tilde{\Lambda}_{2,+}\tilde{\Lambda}_{2,-}$ are bounded from $\mathit{H}^{-1/2}(\S)^4$ to $\mathit{H}^{1/2}(\S)^4$. Hence, $(\mathrm{P2})$ also holds true and thus $H_{2}$ is essentially self-adjoint and generates confinement by Theorem \ref{thincritic} and Proposition \ref{confinementgenerale1}, because $(\mathrm{P}^{\prime}3)$ also holds true. Which proves $(\mathrm{ii})$ and complets the proof of the proposition.\qed


 \section{ On the confinement induced by  delta interactions involving the Cauchy operator}\label{Sec 8}
\setcounter{equation}{0}
In this section, we are interested in the families of Dirac operators given by
\begin{align}
 (-m,m)&\ni a\longmapsto H_{a,\l}= H+ \l \CS^{a}\delta_{\S},\quad\l\in\rr\setminus\{0\},\\
  (-m,m)&\ni a\longmapsto H_{a,\ll}= H+ \ll( \aN) \mathit{C}^{a}_{\S}( \aN)\delta_{\S},\quad \ll\in\rr\setminus\{0\}.
\end{align}

As we have already mentioned in Subsection \ref{Subb2.2}, the above families of Dirac operators involve the Calder\'on projectors
\begin{align}
 \left(\frac{1}{2}\mp i(\aN)\CS^a\right)\, \text{ and }\,  \left(\frac{1}{2}\mp i\CS^a(\aN)\right),
\end{align}
for $\l,\ll\in\{-4,4\}$, and hence they induce confinement. So throughout this section we focus only on those two cases.

First, we study the Dirac operators  $ H_{a,\l}$. As usual we let 
  \begin{align}
\Lambda^{z}_{\l,\pm}=-\frac{4}{\l}(\aN)\CS^a(\aN)\pm\CS^z.
\end{align}
and 
 \begin{align}\label{Hl}
\mathrm{dom}(H_{a,\l})=\left\{ u+\Phi[g]: u\in\mathit{H}^1(\rr^3)^4, g\in\mathit{L}^2(\S)^4 \text{ and } T_{\S}u=-\Lambda_{\l,+}[g]\right\}.
\end{align}
The following proposition is about the basic spectral properties of $ H_{a,\l}$.

\begin{proposition}\label{th8.1} Let $H_\l$ be as in \eqref{Hl}. The following hold true:
\begin{itemize}
  \item[(i)] If $\O_+$ is a UR domain and $\l=4$, then $(H_{a,\l}, \dom(H_{a,\l}))$ is self-adjoint for all $a\in (-m,m)$.  Moreover, if $\O_+$  is Lipschitz, then $a\notin\mathrm{Sp}(H_{a,\l})$, $\S$ is impenetrable and the following hold:
  \begin{itemize}
  \item[(a)]   $H_{a,\l}\varphi=H_{a,\l}^{\Omega_+}\varphi_+\oplus H_{a,\l}^{\Omega_-}\varphi_-=\left(-i\alpha\cdot \nabla+m\beta\right)\varphi_+\oplus\left(-i\alpha\cdot \nabla+m\beta\right)\varphi_-$,
    where $H_{a,\l}^{\Omega_\pm}$ are the self-adjoint Dirac operators defined on
\begin{align*}
\dom(H_{a,\l}^{\Omega_\pm})=\bigg\{ u_{\Omega_\pm}+\Phi_{\O_\pm}[g]:& u_{\Omega_\pm}\in\mathit{H}^1(\Omega_\pm)^4, g\in\mathit{L}^2(\S)^4 \text{ and }\\
& \left(\frac{1}{2}\mp i(\aN)\CS^a\right)(t_\S u_{\Omega_\pm}+\mathit{C}_\pm[g])=0 \bigg\}.
\end{align*}
   \item[(b)] $\mathrm{Sp}_{\mathrm{ess}}(H_{a,\l})=(-\infty,-m]\cup [m,+\infty)$.
  \end{itemize}
  \item[(ii)]   If $\O_+$ is a $\mathit{C}^2$-smooth domain and $\l=-4$, then $(H_{a,\l}, \dom(H_{a,\l}))$ is essentially self-adjoint. Moreover, we have 
$$\overline{H_{a,\l}}\varphi=H_{a,\l}^{\Omega_+}\varphi_+\oplus H_{a,\l}^{\Omega_-}\varphi_-=\left(-i\alpha\cdot \nabla+m\beta\right)\varphi_+\oplus\left(-i\alpha\cdot \nabla+m\beta\right)\varphi_-,$$
 where $H_{a,\l}^{\Omega_\pm}$ are the self-adjoint Dirac operators defined on
\begin{align*}
\dom(H_{a,\l}^{\Omega_\pm})&=\left\{ \varphi_{\pm}\in\mathit{L}^2(\Omega_\pm)^4:  (\alpha\cdot\nabla) \varphi_\pm\in\mathit{L}^2(\Omega_{\pm})^4 \text{ and } \left(\frac{1}{2}\pm i(\aN)\tilde{\CS^a}\right) t_\S\varphi_{\pm}=0 \right\}.
\end{align*}
\end{itemize}
\end{proposition}
 
\textbf{Proof.}  First, we prove $(\mathrm{i})$, so assume that $\O_+$ is a UR domain and $\l=4$. Fix $a\in(-m,m)$, then using the decomposition \eqref{zCau} we obtain that  
\begin{align}\label{5.1.1}
\begin{split}
\Lambda^{z}_{\l,+}= -(\aN)\tilde{W}(\aN)+ \tilde{W}+T^{z}_{K}=-4(\aN)\tilde{W}(\aN)(\frac{1}{4}+ \tilde{W}^2)+T^{z}_{K}
\end{split}
\end{align}
where $T^{z}_{K}\in\mathcal{K}(\mathit{L}^2(\S)^4)$. Since $(\aN)\tilde{W}(\aN)$ and $(1/4+ \tilde{W}^2)$ are invertible in $\mathit{L}^2(\S)^4$, by  \cite[Theorem 1.50 and Theorem 1.51]{Aiena}  it follows that $\Lambda_{\l,+}$ is a Fredholm operator. As $\Lambda_{\l,+}$ is self-adjoint in $\mathit{L}^2(\S)^4$,  by Theorem \ref{luis} we conclude that $(H_{a,\l}, \dom(H_{a,\l}))$ is self-adjoint for all $a\in (-m,m)$, which proves the first statement of $(\mathrm{i})$. Now  assume that $\O_+$ is Lipschitz, and observe that 
\begin{align}\label{5.1.2}
\begin{split}
\Lambda^{a}_{\l,+}= -(\aN)\CS^a(\aN) + \CS^a=-4(\aN)\CS^a(\aN)(\frac{1}{4}+(\CS^a)^2).
\end{split}
\end{align}
Then, the same reasoning as before yields that $\Lambda^{a}_{\l,+}$ is invertible for all $a\in (-m,m)$. Therefore,  by Proposition \ref{surlespectre}-$(\mathrm{i})$ it follows that $a\notin\mathrm{Sp}(H_{a,\l})$. Item $(\mathrm{a})$ and $(\mathrm{b})$ are  consequences  of Proposition \ref{Calderonprojectors} and Proposition \ref{confinementgenerale}, respectively. This finishes the proof of $(\mathrm{i})$.

Now, we prove $(\mathrm{ii})$, so assume that $\O_+$ is $\mathit{C}^2$-smooth and $\l=-4$. It is clear that $\Lambda^{a}_{\l,+}\in\mathcal{K}(\mathit{L}^2(\S)^4)$, therefore $\l=-4$ is a critical parameter. Also, observe that the properties $(\mathrm{P1})$ and $(\mathrm{P}^{\prime}3)$  hold true by Proposition \ref{extension}. Thus, the only thing left to check is the property $(\mathrm{P2})$.  To this end,  recall again the decomposition \eqref{zCau}, then  we make the observation that in order to prove that $\{\aN, \CS^a\}$ extends to a  bounded operator from $\mathit{H}^{-1/2}(\S)^4$ to $\mathit{H}^{1/2}(\S)^4$ (see \cite[Proposition 2.8]{OBV} and \cite[Lemma 3.1 and Remark 3.1]{BB}),  the most delicate part is to show that $\{\aN, \tilde{W}\}$   extends to a  bounded operator from $\mathit{H}^{-1/2}(\S)^4$ to $\mathit{H}^{1/2}(\S)^4$, because the kernel of $T^a_K$ behaves locally as $|x-y|^{-1}$ and thus $T^a_K$ extends to a bounded operator from $\mathit{H}^{-1/2}(\S)^4$ to $\mathit{H}^{1/2}(\S)^4$, even if $\S$ is Lipschitz. Now, a straightforward computation shows that 

\begin{align}
\begin{split}
\tilde{\Lambda}_{\l,+}\tilde{\CS^a} \tilde{\Lambda}_{\l,-}=& \left((\aN)T^a_K(\aN) +T^0_K +\left( (\aN)\tilde{W}(\aN) + \tilde{W}\right)\right)\tilde{\CS^a} \tilde{\Lambda}_{\l,-}\\
=& \left((\aN)T^a_K(\aN) +T^0_K +(\aN)\{\aN,\tilde{W}\}\right)\tilde{\CS^a} \tilde{\Lambda}_{\l,-}.
\end{split}
\end{align}
Combining this with the above observation, and taking into account the fact that $\mathit{N}$ is $\mathit{C}^1$-smooth we then get the property $(\mathrm{P2})$. Therefore, item $(\mathrm{ii})$ follows by Theorem \ref{thincritic} and Proposition \ref{confinementgenerale1}. This achieves the proof of the proposition.\qed
\newline

Now we turn to the analysis of the Dirac operator $H_{a,\ll}$. We recall that 
\begin{align}\label{5.99}
\Lambda^z_{\l^{\prime},+}=-\frac{4}{\ll}\CS^a+\CS,\quad \text{ for all } (a,\ll)\neq(0,4).
\end{align}
The case $(a,\ll)=(0,4)$ is special, we will discuss it separately in the end of this section. The main result about the self-adjointness of $H_{a,\ll}$ reads as follows.

\begin{proposition}\label{th8.2} Let $(H_{\ll}, \dom(H_{\ll}))$ be as above. The following hold:
\begin{itemize}
  \item[(i)]   If $\O_+$ is a UR domain and  $\ll=-4$ , then $(H_{a,{\ll}}, \dom(H_{a,{\ll}}))$ is self-adjoint.  Moreover, if $\O_+$  is Lipschitz, then $a\notin\mathrm{Sp}(H_{a,\ll})$, $\S$ is impenetrable and we have 
    \begin{align}
H_{a,\ll}(\varphi)=H_{a,\ll}^{\Omega_+}(\varphi_+)\oplus H_{a,\ll}^{\Omega_-}(\varphi_-)=\left(-i\alpha\cdot \nabla+m\beta\right)\varphi_+\oplus\left(-i\alpha\cdot \nabla+m\beta\right)\varphi_-
\end{align}
where $H_{a,\ll}^{\Omega_\pm}$ are the self-adjoint Dirac operators defined on
\begin{align*}
\dom(H_{a,\ll}^{\Omega_\pm})=\bigg\{u_{\Omega_\pm}+\Phi_{\O_\pm}[g]:& u_{\Omega_\pm}\in\mathit{H}^1(\Omega_\pm)^4, g\in\mathit{L}^2(\S)^4 \text{ and }\\
&\left(\frac{1}{2}\pm i\CS^a(\aN)\right) (t_\S u_{\Omega_\pm}+ \mathit{C}_\pm[g])=0 \bigg\}.
\end{align*}
  \item[(ii)] If $\O_+$ is a $\mathit{C}^2$-smooth domain, $\ll=4$ and $a\neq0$, then $(H_{a,\ll}, \dom(H_{a,\ll}))$ is essentially self-adjoint. furthermore  we have 
 $$\overline{H_{a,\ll}}(\varphi)=H_{a,\ll}^{\Omega_+}(\varphi_+)\oplus H_{a,\ll}^{\Omega_-}(\varphi_-)=\left(-i\alpha\cdot \nabla+m\beta\right)\varphi_+\oplus\left(-i\alpha\cdot \nabla+m\beta\right)\varphi_-,$$
where $H_{a,\ll}^{\Omega_\pm}$ are the self-adjoint Dirac operators defined on
\begin{align*}
\dom(H_{a,\ll}^{\Omega_\pm})&=\left\{ \varphi_{\pm}\in\mathit{L}^2(\Omega_\pm)^4:  (\alpha\cdot\nabla) \varphi_\pm\in\mathit{L}^2(\Omega_{\pm})^4 \text{ and }\left(\frac{1}{2}\mp i\tilde{\CS^a}(\aN)\right) t_\S\varphi_{\pm}=0 \right\}.
\end{align*}
\end{itemize}
\end{proposition}
We omit the proof of this proposition, since it is easier than, and can easily be extracted from the proof of Proposition \ref{th8.1}.
\newline

In the following, we describe the property of the confinement induced by $H_{a,\ll}$ when $(a,\ll)=(0,4)$. From \eqref{5.99}, we notice that $\Lambda_{\l^{\prime},+}$ vanishes in this case. Indeed, let $\varphi=u+\Phi[g]$ with $u\in\mathit{H}^1(\rr^3)^4$ and $g\in\mathit{L}^2(\S)^4$, then  
\begin{align}
\label{}
H_{a,\ll}(\varphi)= H(u)+ (\aN)\CS(\aN)t_\S u,
\end{align}
holds in the sense of distribution. Thus, we need to assume that $u\in\mathit{H}^1_0(\rr^3\setminus\S)^4$ in order to ensure that 
$H_{a,\ll}(\varphi)\in\mathit{L}^2(\rr^3)^4$. Hence, we cannot define $H_{a,\ll}$ as we did before. To work around this problem, note that for $\varphi=(\varphi_+,\varphi_-)\in\mathit{H}^1(\Omega_+)^4\oplus\mathit{H}^1(\Omega_-)^4$, a simple computation in the sense of distributions yields
\begin{align*}
 H_{0,\ll}(\varphi)=& H(\varphi) +2(\aN)\CS(\aN)(t_{\S}\varphi_+ +t_{\S}\varphi_-)\delta_{\S}\\
=& H(\varphi_+)\oplus H(\varphi_-) +i\alpha\cdot\mathit{N}(t_{\S}\varphi_+ -t_{\S}\varphi_-)\delta_{\S}\\
&+2(\aN)\CS(\aN)(t_{\S}\varphi_+ +t_{\S}\varphi_-)\delta_{\S}.
\end{align*}
Thus, if we let 
\begin{align}\label{transf}
\left(\frac{1}{2}-i\CS(\aN)\right)t_{\S}\varphi_+ = \left(\frac{1}{2}+i\CS(\aN)\right)t_{\S},\varphi_-
\end{align}
 then $H_{0,\ll}(\varphi)\in\mathit{L}^2(\rr^3)^4$. In particular, this leads us to define $H_{0,\ll}$ as follows:
 \begin{align}\label{domconf}
\mathrm{dom}(H_{0,\ll})=\left\{ \varphi=(\varphi_+,\varphi_-)\in\mathit{H}^1(\Omega_+)^4\oplus\mathit{H}^1(\Omega_-)^4 :  \eqref{transf} \text{ holds in } \mathit{H}^{1/2}(\S)^4 \right \}.
\end{align}
Clearly,  $H_{0,\ll}$ is symmetric.  Moreover, if we let $P_\mp=1/2\mp i\CS(\aN)$, then it is straightforward to check that 
\begin{align}\label{danslapr}
H_{0,\ll}(\varphi)=H_{\Omega_+}(\varphi_+)\oplus H_{\Omega_-}(\varphi_-)=\left(-i\alpha\cdot \nabla+m\beta\right)\varphi_+\oplus\left(-i\alpha\cdot \nabla+m\beta\right)\varphi_-,
\end{align}
where $H_{\Omega_\pm}$ are the symmetric Dirac operators defined on
\begin{align*}
\dom(H_{\Omega_\pm})&=\left\{ \varphi_{\pm}\in\mathit{H}^1(\Omega_\pm)^4: \left(\frac{1}{2}\mp i\CS(\aN)\right) t_\S\varphi_{\pm}=0 \right\}.
\end{align*}
Then, we have the following theorem about the self-adjointness of   $H_{a,\ll}$ when $(a,\ll)=(0,4)$.  In the proof, we use the notation $\langle\, , \, \rangle_{\mathit{H}^{-1/2},\mathit{H}^{1/2}}$ for the duality pairing between $\mathit{H}^{-1/2}(\S)^4$ and $\mathit{H}^{1/2}(\S)^4$.
\begin{theorem} Assume that $(a,\ll)=(0,4)$ and let $H_{a,\ll}$ be as in \eqref{domconf}. Then $H_{a,\ll}$ is self-adjoint, the restriction  $H_{a,\ll}\downharpoonright \mathit{H}^1(\rr^3\setminus\S)^4$ is essentially self-adjoint. Moreover, $\S$ is impenetrable and we have $\overline{H_{a,\ll}}=\overline{H_{\Omega_+}}\oplus \overline{H_{\Omega_-}}$, with 
 \begin{align}\label{adjoint22}
\dom(\overline{H_{\Omega_\pm}})&=\left\{  \psi\in\mathit{L}^2(\Omega_+)^4:  (\alpha\cdot\nabla) \psi\in\mathit{L}^2(\Omega_{+})^4 \text{ and } P_\mp t_{\S}\psi=0 \right\},
\end{align}
where the boundary condition has to be understood as an equality in $\mathit{H}^{-1/2}(\S)^4$.
\end{theorem} 
\textbf{Proof.} The proof is standard and follows essentially the same idea as in \cite[Theorem 3.2 and Theorem 4.2]{OBV}. Indeed, due to the decomposition \eqref{danslapr}, it is sufficient to prove that both $H_{\Omega_+}$ and $H_{\Omega_-}$ are self-adjoint. In what follows we deal with the self-adjointness of  $H_{\Omega_+}$ only, since $H_{\Omega_-}$  can be treated analogously. For the convenience of the reader,   the proof  is divided into two streps  as follows:
\begin{itemize}
  \item[(a)] The domain of $H^{\ast}_{\Omega_+}$ is given by
  \begin{align}\label{adjointt}
\dom(H^{\ast}_{\Omega_+})&=\left\{  \psi\in\mathit{L}^2(\Omega_+)^4:  (\alpha\cdot\nabla) \psi\in\mathit{L}^2(\Omega_{+})^4 \text{ and } P_-t_{\S}\psi=0 \right\},
\end{align}
where the boundary condition has to be understood as an equality in $\mathit{H}^{-1/2}(\S)^4$.
  \item [(b)] The inclusion $\overline{H_{\O_+}}\subset H^{\ast}_{\Omega_+}$ holds.
\end{itemize}
Once $(\mathrm{a})$ and $(\mathrm{b})$ are proved, we use the fact that $H_{\O_+}$ is symmetric and then we conclude that  $\overline{H_{\O_+}}= H^{\ast}_{\Omega_+}$.

 \textbf{Proof of $(\mathrm{a})$.} Denote by $D$ be the set on the right-hand of \eqref{adjointt}. First, we show the inclusion $D\subset\dom(H^{\ast}_{\Omega_+})$, for this let $\varphi\in\dom(H_{\Omega_+}) $ and $\psi\in D$. Then, using \cite[Corollary
2.15]{OBV} it follows that 
\begin{align}\label{5.16}
 \langle H(\psi), \varphi\rangle_{\mathit{L}^2(\O_+)^4}= \langle\psi, H(\varphi)\rangle_{\mathit{L}^2(\O_+)^4} + \langle -i(\aN) t_{\S}\psi, t_{\S}\varphi\rangle_{\mathit{H}^{-1/2},\mathit{H}^{1/2}}.
 \end{align} 
Now, using \eqref{dual} and the fact that  $-i(\aN)t_{\S}\psi=2(\aN)\tilde{\CS}(\aN) t_{\S}\psi$ , it follows that 
\begin{align}
\begin{split}
  \langle -i(\aN) t_{\S}\psi, t_{\S}\varphi\rangle_{\mathit{H}^{-1/2},\mathit{H}^{1/2}}&= \langle 2(\aN)\tilde{\CS}(\aN) t_{\S}\psi, t_{\S}\varphi\rangle_{\mathit{H}^{-1/2},\mathit{H}^{1/2}}\\
  &= \langle -i(\aN) t_{\S}\psi, -2i\CS(\aN)t_{\S}\varphi\rangle_{\mathit{H}^{-1/2},\mathit{H}^{1/2}}.
  \end{split}
 \end{align} 
 Similarly,  using  that  $t_{\S}\varphi=2i\CS(\aN) t_{\S}\varphi$ , we get that
 \begin{align}
 \begin{split}
 \langle -i(\aN) t_{\S}\psi, t_{\S}\varphi\rangle_{\mathit{H}^{-1/2},\mathit{H}^{1/2}}= -\langle -i(\aN) t_{\S}\psi, -2i\CS(\aN)t_{\S}\varphi\rangle_{\mathit{H}^{-1/2},\mathit{H}^{1/2}}.
   \end{split}
 \end{align} 
From this, we conclude that 
\begin{align}\label{5.19}
 \begin{split}
 \langle -i(\aN) t_{\S}\psi, t_{\S}\varphi\rangle_{\mathit{H}^{-1/2},\mathit{H}^{1/2}}= 0.
   \end{split}
     \end{align} 
Therefore, we obtain 
 \begin{align*}
 \langle H(\psi), \varphi\rangle_{\mathit{L}^2(\O_+)^4}= \langle\psi, H(\varphi)\rangle_{\mathit{L}^2(\O_+)^4},
 \end{align*} 
which yields the inclusion $D\subset\dom(H^{\ast}_{\Omega_+})$.  We now prove the reverse inclusion. Given $\varphi\in\mathcal{D}(\O_+)^4$ and let $\psi\in\dom(H^{\ast}_{\Omega_+})$. Then, by definition there is $\chi\in\mathit{L}^2(\O_+)^4$ such that 
\begin{align*}
 \langle H(\psi), \varphi\rangle_{\mathcal{D}^{\prime}(\O_+)^4,\mathcal{D}(\O_+)^4}&= \langle\psi, H(\varphi)\rangle_{\mathcal{D}^{\prime}(\O_+)^4,\mathcal{D}(\O_+)^4}\\
 &= \langle\psi, H(\varphi)\rangle_{\mathit{L}^2(\O_+)^4}= \langle\psi, \chi\rangle_{\mathcal{D}^{\prime}(\O_+)^4,\mathcal{D}(\O_+)^4}.
 \end{align*} 
 Thus, we get that $H(\psi)=\chi$ in $\mathcal{D}^{\prime}(\O_+)^4$ and then in $\mathit{L}^2(\O_+)^4$. Hence, $\psi,(\alpha\cdot\nabla) \psi\in\mathit{L}^2(\Omega_{+})^4$, so it remains to show that $P_-t_{\S}\psi=0$ in $\mathit{H}^{-1/2}(\S)^4$. To this end, recall the definition of the extension operator $E_{\O_+}$ from Subsection \ref{BesovSobolev}.  Observe that $E_{\O_+}(P_+g)\in\dom(H_{\Omega_+})$, for all $g\in\mathit{H}^{1/2}(\S)^4$. Hence, from \eqref{5.16} and \eqref{5.19} it follows that 
 \begin{align}
 \begin{split}
 \langle -i(\aN) t_{\S}\psi, P_+g\rangle_{\mathit{H}^{-1/2},\mathit{H}^{1/2}}= 0.
 \end{split}
 \end{align} 
Thus, we get 
 \begin{align}
 \begin{split}
 \langle -i(\aN) t_{\S}\psi, g\rangle_{\mathit{H}^{-1/2},\mathit{H}^{1/2}}&=  \langle -i(\aN) t_{\S}\psi, P_-g\rangle_{\mathit{H}^{-1/2},\mathit{H}^{1/2}}.
 \end{split}
 \end{align}  
  Now, using \eqref{dual} and the identity $-4(\CS(\aN))^2=\mathit{I}_4$, a simple computation yields
 \begin{align}
 \begin{split}
 \langle -i(\aN) t_{\S}\psi, g\rangle_{\mathit{H}^{-1/2},\mathit{H}^{1/2}}&=  \langle -i(\aN) t_{\S}\psi, -4P_-(\CS(\aN))^2g\rangle_{\mathit{H}^{-1/2},\mathit{H}^{1/2}}\\
 &=  \langle 2i\tilde{\CS}(\aN) t_{\S}\psi, i(\aN)P_-g\rangle_{\mathit{H}^{-1/2},\mathit{H}^{1/2}}.
 \end{split}
 \end{align}  
 Therefore, we get 
 \begin{align}
 \begin{split}
 \langle \left(\frac{1}{2}-i\tilde{\CS}(\aN) \right) t_{\S}\psi, g\rangle_{\mathit{H}^{-1/2},\mathit{H}^{1/2}}&= 0.
  \end{split}
 \end{align}  
 Since this is true for all $g\in\mathit{H}^{1/2}(\S)^4$, it follows that $\psi\in D$. Hence $\dom(H^{\ast}_{\Omega_+})\subset D$, which proves $(\mathrm{a})$. 
 
\textbf{Proof of $(\mathrm{b})$.}  Fix $\psi\in\dom(H^{\ast}_{\Omega_+})$ and let $(g_j)_{j\in\mathbb{N}}=(P_+h_j)_{j\in\mathbb{N}}\subset \mathit{H}^{1/2}(\S)^4$ be a sequence of functions that converges to $t_{\S}\psi$ in $\mathit{H}^{-1/2}(\S)^4$. Set 
   \begin{align}\label{suite}
   \begin{split}
   \psi_j&= \psi +\Phi_{\O_+}[i(\aN)(g_j - t_{\S}\psi])+ E_{\O_+}\left(\{\aN,\tilde{\CS}\}(\aN)(g_j- t_{\S}\psi)\right)\\
   &:= \psi +F_1+ F_2.
   \end{split}
   \end{align}
 Clearly, $\psi_j,(\alpha\cdot\nabla) \psi_j\in\mathit{L}^2(\Omega_{+})^4$, for all $j\in\mathbb{N}$. Now observe that 
  \begin{align}
 \begin{split}
t_{\S}F_1=P_+(g_j- t_{\S}\psi) \, \text{ and }\, t_{\S}F_2=\{\aN,\tilde{\CS}\}(\aN)(g_j- t_{\S}\psi).
  \end{split}
 \end{align}  
Hence we get 
  \begin{align}\label{tracede la fon}
 \begin{split}
t_{\S}\psi_j=g_j + \{\aN,\tilde{\CS}\}(\aN)(g_j- t_{\S}\psi).
  \end{split}
 \end{align} 
 Since  $\{\aN,\tilde{\CS}\}$ is bounded from $\mathit{H}^{-1/2}(\S)^4$ into $\mathit{H}^{1/2}(\S)^4$, it follows that $t_{\S}\psi_j\in\mathit{H}^{1/2}(\S)^4$. Therefore, $\psi_j\in\mathit{H}^1(\Omega_+)^4$ holds by \cite[Proposition 2.16]{OBV}. As $P_-g_j=0=P_-t_{\S}\psi$,  we get that
  \begin{align*}
2i\CS(\aN)\{\aN,\tilde{\CS}\}(\aN)(g_j- t_{\S}\psi)=\{\aN,\tilde{\CS}\}(\aN)(g_j- t_{\S}\psi).
 \end{align*} 
 Using this and the fact that $g_j=P_+g_j$, from \eqref{tracede la fon} it follows that $P_-t_{\S}\psi_j=0$. Thus, $\psi_j\in\dom(H_{\O_+})$, for all $j\in\mathbb{N}$. Now, by Proposition \ref{extension} $(\mathrm{i})$-$(\mathrm{ii})$, we obtain that 
   \begin{align}\label{cv}
\psi_j \xrightarrow[j\to\infty]{} \psi \,\text{ in } \mathit{L}^{2}(\O_+)^4.
\end{align}
 Next, by \cite[Theorem 2.2]{OBV} and the trace theorem, there is $C>0$ such that
   \begin{align}
  \Vert H(\psi_j -\psi)\Vert_{\mathit{L}^2(\O_+)^4}^2\leqslant    \Vert t_{\S}\psi_j -t_{\S}\psi\Vert_{\mathit{H}^{-1/2}(\S)^4}^2.
\end{align}
Thus
  \begin{align}\label{cv}
H(\psi_j) \xrightarrow[j\to\infty]{} H(\psi) \,\text{ in } \mathit{L}^{2}(\O_+)^4.
\end{align}
Summing up, we have proved that  $(\psi_j, H_{\Omega_+}(\psi_j))$ convergences to $(\psi, H^{\ast}_{\Omega_+}(\psi))$ when $j$ tends to infinity. Therefore,  $\overline{H_{\O_+}}\subset H^{\ast}_{\Omega_+}$ and this completes the proof of $(\mathrm{b})$.

Finally, it remains to prove that $\overline{H_{\Omega_+}}\not\subset H_{\Omega_+}$. Pick  $g\in \mathit{H}^{-1/2}(\S)^4\setminus \mathit{L}^{2}(\S)^4$ and set $\psi=E_{\O_+}( P_+g)$. Then $\psi\in\dom(\overline{H_{\Omega_+}})$  and $\psi\notin \dom(H_{\Omega_+})$,  as otherwise $g\in \mathit{H}^{1/2}(\S)^4$ by \cite[Proposition 2.7]{OBV}. This achieves the proof of the theorem. \qed
\newline

\begin{remark} It should be noted that the reason why  $H_{a,\l}\downharpoonright \mathit{H}^1(\rr^3\setminus\S)^4$ and $H_{a,\ll}\downharpoonright \mathit{H}^1(\rr^3\setminus\S)^4$ are not self-adjoint for critical parameters, is due to the fact that we are projecting in the wrong direction. In other words, we have forced the the terms (more precisely the Calder\'on projectors associated to each problems) that allow us to regularize the functions in $\dom(\overline{H_{a,\l}})$ (respectively  in $\dom(\overline{H_{a,\ll}})$) to be zero; see e.g. \cite[Proposition 2.7]{OBV}.
\end{remark}

\section*{Acknowledgement} I would like to warmly thank my PhD supervisors Vincent Bruneau and Luis Vega for their encouragement, and for their precious discussions and advices during the preparation of this paper. I would also like to thank Mihalis Mourgoglou for various fruitful discussions on topics related to Analysis of PDEs in rough domains. I am particularly grateful to him for suggesting several references related to Section  \ref{sec5} and for his careful reading of a preliminary version of that section.  This project is based upon work supported by the Government of the Basque Country under Grant PIFG$18/06$ "ERC Grant: Harmonic Analysis and Differential Equations", researcher in charge: Luis Vega.


\end{document}